\documentclass[11pt]{amsart}
\usepackage{amsmath}
\usepackage{amsxtra}
\usepackage{amscd}
\usepackage{amsthm}
\usepackage{amsfonts}
\usepackage{amssymb}
\usepackage{psfig}
\usepackage{eucal}
\newcommand{\cal}{\mathcal}
\newtheorem{thm}{Theorem}
\newtheorem{lem}{Lemma} \newtheorem{prop}{Proposition}
 
\newtheorem{cor}{Corollary}
\newtheorem{conj}{Conjecture}

\newcommand{\on}{\operatorname} 
\newcommand{\proj}{\mathbf P}
\newcommand{\psdraw}[3]{\begin{array}{c} \hspace{-1mm}
\raisebox{-4pt}{\psfig{figure=#1.ps,width=#2,height=#3}}
\hspace{-1mm}\end{array}}
\numberwithin{thm}{section} 
\numberwithin{lem}{section}
\numberwithin{cor}{section} 
\numberwithin{prop}{section} 
\numberwithin{equation}{section}
\numberwithin{conj}{section}
\title{Moduli of Trigonal Curves}
\author{Zvezdelina E. Stankova-Frenkel}
\begin{document}

\thispagestyle{empty}
\maketitle

\begin{abstract} We study the moduli of trigonal curves.
We establish the exact upper bound of ${36(g+1)}/(5g+1)$
for the slope of trigonal fibrations. Here, 
the slope of any fibration $X\rightarrow B$ of stable curves with
smooth general member is the ratio $\delta_B/\lambda_B$ of the
restrictions of the boundary class $\delta$ and the Hodge class $\lambda$
on the moduli space $\overline{\mathfrak{M}}_g$ to the base $B$.
We associate to a trigonal family $X$ a canonical rank two vector 
bundle $V$, and show that for Bogomolov-semistable
$V$ the slope satisfies the stronger inequality
${\delta_B}/{\lambda_B}\leq 7+{6}/{g}$.
We further describe the rational Picard group
of the {trigonal} locus
$\overline{\mathfrak T}_g$ in the moduli space $\overline{\mathfrak{M}}_g$
of genus $g$ curves. In the even genus case, we interpret the above
Bogomolov semistability condition in
terms of the so-called Maroni divisor in $\overline{\mathfrak T}_g$. 
\end{abstract}

\tableofcontents

\section*{\hspace*{1.9mm}1. Introduction}

\setcounter{section}{1} 
\setcounter{subsection}{0}
\setcounter{lem}{0}
\setcounter{thm}{0}
\setcounter{prop}{0}
\setcounter{defn}{0} 
\setcounter{cor}{0} 
\setcounter{conj}{0} 
\setcounter{claim}{0} 
\setcounter{remark}{0}
\setcounter{equation}{0}
\label{introduction}

In this paper $\overline{\mathfrak M}_g$ denotes the
Deligne-Mumford compactification of the moduli space
of smooth curves over $\mathbb{C}$ of genus $g\geq 2$. 
The boundary locus $\Delta$ of
$\overline{\mathfrak M}_g$ consists of nodal curves with finite automorphism
groups, which together with the smooth curves are referred to as {\it
stable} curves. The locus of hyperelliptic curves will be denoted by
$\overline{\mathfrak{I}}_g$, and the closure of the locus of trigonal curves
will be denoted by $\overline{\mathfrak{T}}_g$.

\medskip  The main objects of our study will be families  
of genus $g$ stable curves, whose general members
are smooth. Associated to any such {\it flat and proper}
family $f\!:\!X\!\rightarrow\! B$ 
 are three basic invariants $\lambda|_B$, $\delta|_B$ and $\kappa|_B$.
We define these in Section~\ref{definition}
as divisors on $B$, but
for most purposes one can think of them as integers by considering their
respective degrees. The invariant $\delta|_B$ counts, with appropriate 
multiplicities, the number of singular fibers of $X$.
The self-intersection of the relative dualizing sheaf $\omega_f$ on
$X$ defines $\kappa|_B$, and its pushforward to $B$ is a  rank $g$ 
locally free sheaf, whose determinant is $\lambda|_B$.

\smallskip
The basic relation $12\lambda|_B=\delta|_B+\kappa|_B$ and the positivity of
the three invariants for non-isotrivial families force the {\it slope}
$\displaystyle{\frac{\delta|_B}{\lambda|_B}}$ of $X/\!_{\displaystyle{B}}$
to fall into the interval $[0,12)$ (cf.~Sect.~\ref{slope-non-isotrivial}). 
In fact, Cornalba-Harris and Xiao establish for this slope 
an exact upper bound of $8+4/g$, which is achieved only for certain
hyperelliptic families (cf.~Theorem~\ref{CHX}). 
However, if the base curve $B$ passes through
a {\it general}
 point of $\overline{\mathfrak{M}}_g$, Mumford-Harris-Eisenbud 
give the better bound of $6+\on{o}(1/g)$ (cf.~Theorem~\ref{generalbound1}). 
The families violating this inequality are 
entirely contained in the closure $\overline{\cal{D}}_k$ 
of the locus of $k$-sheeted covers of ${\proj}^1$, for
a suitably chosen $k$. In particular, for $k=2$ we recover the 
{hyperelliptic} locus $\overline{\mathfrak{I}}_g$, 
for $k=3$ - the {trigonal} locus $\overline{\mathfrak{T}}_g$, etc. 
Therefore, higher than the above
``generic'' ratio can be obtained only for families with special
linear series, such as $g^1_2$, $g^1_3$, etc. 
These observations clearly raise the following

\medskip
\noindent{\bf Question.} {\it According to the possession of special 
linear series, is there a stratification of $\overline{\mathfrak M}_g$ 
which would give successively smaller slopes $\delta/\lambda$? What 
would be the successive upper bounds with respect to such a stratification?}

\smallskip 
The following result, whose proof will be given in the paper,
answers this question for an exact upper bound for
families with linear series $g^1_3$.
 
\medskip
\noindent{\bf Theorem I.} 
{\it If $f\!:\!X{\rightarrow} B$ is a trigonal nonisotrivial family 
with smooth general member, then the slope of $X/\!_{\displaystyle{B}}$ 
satisfies:
\begin{equation*}
\frac{\delta|_B}{\lambda|_B}\leq \frac{36(g+1)}{5g+1}\cdot
\end{equation*}
Equality is achieved if and only if all fibers are irreducible,
$X$ is a triple cover of a ruled surface $Y$ over $B$,
 and a certain divisor class $\eta$ on $X$ is numerically zero.}

\smallskip
To understand the importance of this result and the above question, 
consider Mumford's alternative description 
of the basic invariants (cf.~ Sect.~\ref{linebundles}):
$\lambda|_B$, $\delta|_B$ and $\kappa|_B$ are restrictions of certain 
rational divisor classes $\lambda, \delta, \kappa\in
\on{Pic}_{\mathbb{Q}}\overline{\mathfrak{M}}_g$. Specifically,  
$\delta=\delta_0+\cdots+\delta_{[g/2]}$, where $\delta_i$ the class
of the boundary divisor $\Delta_i$ of $\overline{\mathfrak{M}}_g$,
and $\on{Pic}_{\mathbb{Q}}\overline{\mathfrak{M}}_g$ 
is freely generated by the Hodge class $\lambda$  and the boundary
classes $\delta_i$ for $g\geq 3$ (cf.~\cite{Ha2}).
Thus, our question about a stratification of $\overline{\mathfrak{M}}_g$ 
translates into a question about 
the relations among the fundamental classes of various
subvarieties defined by geometric conditions in $\overline{\mathfrak{M}}_g$ .
Moreover, such a stratification would provide a link between
the {\it global} invariant $\lambda$ (the degree of the Hodge bundle on
$\overline{\mathfrak M}_g$) and the {\it locally defined} invariant $\delta$ of
the singularities of our families. In the process of estimating the 
ratio $\delta / \lambda$ we hope to understand the geometry of interesting
loci in $\overline{\mathfrak M}_g$, and describe their rational Picard groups.

\smallskip
Such a program for the hyperelliptic locus $\overline{\mathfrak{I}}_g$
is completed by Cornalba-Harris (cf.~Theorem~\ref{theoremCHPic}), 
who exhibit generators and relations 
for $\on{Pic}_{\mathbb{Q}}{\overline{\mathfrak{I}}_g}$.
The typical examples of families with maximal ratio of $8+{4}/{g}$ are 
constructed as blow-ups of pencils of hyperelliptic curves, 
embedded in the same ruled surface. 

\smallskip
Similar examples for
{trigonal families} yield the slope $7+{6}/{g}$, but as Theorem I
shows, this ratio is {\it not} an upper bound. This happens because of
an ``extra'' {\it Maroni} locus in $\overline{\mathfrak{T}}_g$ (cf.~
Sect.~12).
While a general trigonal curve embeds in ${\mathbf F}_0={\proj}^1\times
{\proj}^1$ or in the blow-up ${\mathbf F}_1$ 
of ${{\proj}^2}$ at a point, the remaining
trigonal curves embed in other rational 
ruled surfaces and comprise a closed subset in
$\overline{\mathfrak{T}}_g$, called the Maroni locus. The proof of
Theorem II, stated below, implies that 
all trigonal families achieving the maximal bound lie entirely in the
Maroni locus, and moreover, their members 
are embedded in ruled surfaces ``as far as possible
from the generic'' ruled surfaces ${\mathbf F}_0$ and ${\mathbf F}_1$.

\smallskip
The ratio $7+{6}/{g}$, though not the ``correct'' maximum,
plays a significant role in understanding the geometry of the trigonal
locus, and in describing its rational Picard group. In particular, in 
 a linear relation is established between the Hodge class,
the boundary classes on $\overline{\mathfrak{T}}_g$, and 
a canonically defined vector bundle $V$ of rank 2 on a ruled surface 
$\widehat{Y}$ (cf.~Sect.~9):

\medskip
\noindent{\bf Theorem II.}
{\it  Let $\delta_0$ denote the boundary class in $\overline
{\mathfrak{T}}_g$ corresponding to irreducible singular curves, and let 
$\delta_{k,i}$ be the remaining boundary classes.
For any trigonal non-isotrivial family with general smooth member, we have
\begin{equation*}
(7g+6)\lambda|_B=g\delta_0|_B+\sum_{k,i} \widetilde{c}_{k,i}\delta_{k,i}|_B
+\frac{g-3}{2}(4c_2(V)-c_1^2(V)),
\end{equation*}
where $\widetilde{c}_{k,i}$ 
is a quadratic polynomial in $i$ with linear coefficients
in $g$, and it is determined by the geometry of $\delta_{k,i}$.}

\smallskip
For example, $\widetilde{c}_{1,i}=3(i+2)(g-i)/2$
corresponds to the boundary divisor $\Delta{\mathfrak{T}}_{1,i}$,
whose general member is the join in three points of two trigonal curves
of genera $i$ and $g-i-2$, respectively (cf.~Fig.~\ref{Delta-k,i}).

\medskip
Recall that the vector bundle $V$ is called {\it Bogomolov semistable}
if its Chern classes satisfy $4c_2(V)\geq c_1^2(V)$ (cf.~\cite{Bo,Re}). 
We show in Section~9 the following

\medskip
\noindent{\bf Theorem III.} 
{\it For any trigonal nonisotrivial family $X\rightarrow B$ 
with general smooth member, if $V$ is Bogomolov semistable, then
\begin{equation*}
 \frac{\delta|_B}{\lambda|_B}\leq 7+\frac{6}{g}\cdot
\end{equation*}}

\medskip
In the even genus case, the Maroni locus is in fact a divisor on
$\overline{\mathfrak{T}}_g$, whose class we denote by $\mu$. 
We further recognize the ``Bogomolov quantity''
$4c_2(V)-c_1^2(V)$ as counting roughly four times the number of
Maroni fibers in $X$, and deduce

\medskip
\noindent{\bf Theorem IV.}
{\it For even $g$, $\on{Pic}_{\mathbb{Q}}\overline{\mathfrak{T}}_g$ is
freely generated by all boundary divisors $\delta_0$ and $\delta_{k,i}$, and
the Maroni divisor $\mu$. The class of the Hodge bundle on
$\overline{\mathfrak{T}}_g$ is expressed in terms of these generators as
the following linear combination:
\begin{equation*}
(7g+6)\lambda|_{\overline{\mathfrak{T}}_g}=g\delta_0+
\sum_{k,i}\widehat{c}_{k,j}\delta_{k,i}+2(g-3){\mu}.
\end{equation*}}
Consequently, the condition $\eta\equiv 0$ in Theorem I can be interpreted
as a relation among the number of irreducible singular curves and the
``Maroni'' fibers in our family: $(g+2)\delta_0|_B=-72(g+1)\mu|_B$, and hence
maximal slope families are entirely contained in the Maroni locus of 
$\overline{\mathfrak{T}}_g$ (cf.~Theorem~\ref{maximalmaroni}).
The stated theorems complete the program for the trigonal locus 
$\overline{\mathfrak{T}}_g$, which was outlined earlier in this section.

\smallskip
An important interpretation of these results can be traced
back to \cite{MHE}, where it is shown that the moduli space 
$\mathfrak{M}_g$ is of {\it general type}. The $k$-gonal locus 
$\overline{\cal{D}}_k$ is realized in terms of the generating classes as:
$[\overline{\cal{D}}_k]=a\lambda-b\delta-\cal{E}$ for some $a,b>0$,
and an effective boundary combination $\cal{E}$. 
Restricting to a general curve $B\subset \overline{\mathfrak{M}}_g$, we have 
$\overline{\cal{D}}_k|_B> 0$, and hence
$a\lambda|_B-b\delta|_B>0$. Because of 
Seshadri's criterion for ampleness of line bundles,
in effect, we are asking for all positive numbers $a$ and $b$ such that
the linear combination $a\lambda-b\delta$ is ample on $\overline{\mathfrak{M}}_g$.
The smaller the ratio $a/b$ is, the stronger result we obtain.
In other words, we are aiming at a maximal bound of $\delta /\lambda$,
when we think of these classes as restricted to any curve $B\subset
\overline{\mathfrak{M}}_g$.
In view of this, part of this paper can be described as looking for all {\it
ample} divisors on $\overline{\mathfrak{T}}_g$ of the form
$a\lambda-b\delta$ with $a,b>0$. Theorem I then gives the necessary 
condition $(5g+1)a\geq 36(g+1)b$ (compare with \cite{M2,MHE,CH}).
Some of the results can be applied to the study of 
the discriminant loci of a certain type of triple
covers of surfaces. 

\smallskip
The methods and ideas for the trigonal case are
in principal extendable to more general families of $k$-gonal curves.
We refer the reader to  Sect.~13
for a general maximal bound for tetragonal curves (for
$g$ odd), and conjectures for the maximal
and general bounds for any $d$-gonal and other families of stable curves. 

\bigskip
\begin{center}{\sc Acknowledgments}\end{center}

\medskip
This paper is based on my Ph.D. thesis at Harvard University.
Joe Harris, my advisor, introduced me to the problem of finding
a stratification of the moduli space with respect to a descending sequence
of slopes of one--parameter families. I am very grateful to him for
his advice and support throughout my work on the present thesis.
I would like to thank Fedor Bogomolov, 
David Eisenbud, Benedict Gross, Brendan Hassett,
David Mumford, Tony Pantev and Emma Previato for the helpful discussions
that I have had with them at different stages of the project, as well as
Kazuhiro Konno for providing me with his recent results on trigonal
families. A source of inspiration and endless moral support has been my
husband, Edward Frenkel, to whom goes my gratitude and love.

\bigskip
\section*{\hspace*{1.9mm}2. Preliminaries}

\setcounter{section}{2}
\setcounter{subsection}{0}
\setcounter{lem}{0}
\setcounter{thm}{0}
\setcounter{prop}{0}
\setcounter{defn}{0} 
\setcounter{cor}{0} 
\setcounter{conj}{0} 
\setcounter{claim}{0} 
\setcounter{remark}{0}
\setcounter{equation}{0}
\label{preliminaries}

\subsection{Definition of $\lambda|_B,\,\,\delta|_B$ and
$\kappa|_B$ in Pic$B$} 
\label{definition} 
\begin{figure}[h]
$$\psdraw{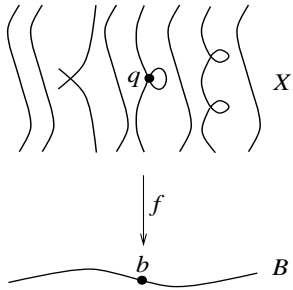}{1.5in}{1.5in}$$
\hspace*{3mm}
\caption{Trigonal family $f\!:\!X\!\rightarrow\! B$}
\label{family}
\end{figure}

Let $f:X\rightarrow B$ be a flat proper
one-parameter family of stable curves of genus $g$, where
$B$ is a smooth projective curve. Assume in addition that the
general member of $X$ is {\it smooth} (cf.~Fig.~\ref{family}). 

\medskip

Let $\omega_f=\omega_X\otimes f^*\omega_B^{-1}$ be the relative dualizing 
sheaf of $f$. Its pushforward $f_*(\omega_f)$ is a locally free sheaf 
on $B$ of rank $g$, and we set
\[\lambda|_B=\lambda_X:=\wedge^gf_*(\omega_f)\in \on{Pic}B\]
to be its determinant. The sheaf $f_*(\omega_f)$ is known as the ``Hodge
bundle'' on $B$, and $\lambda|_B$ - as the ``Hodge class'' of $B$.
In a similar way, we set $\kappa|_B$ to be the self-intersection of $\omega_f$:
\[\kappa|_B=\kappa_X:=f_*(c_1^2(\omega_f))\in\on{Pic}B.\]

\bigskip
The definition of $\delta|_B$, on the other hand, is local and requires
some notation. Let $q$ be any singular point of a fiber $X_b$,
$b\in B$. Since the general fiber of $X$ is smooth, the total space
of $X$ near $q$ is given locally analytically by $xy=t^{m_q}$, 
where $x$ and $y$  are local parameters on $X_b$, $t$ is a local parameter
on $B$ near $b$, and $m_q\geq 1$.
(This follows from the one-dimensional versal
deformation space  of the nodal singularity at $q$.) For any other
point $q$ of $X$ we set $m_q=0$.
Now we can define
\[\delta|_B=\delta_X:=f_*(\sum_{q\in X}m_qq)\in \on{Pic}B.\]

By abuse of notation, we shall use the same letters for the line bundles
$\lambda|_B, \,\,\kappa|_B$ and $\delta|_B$ and for their respective degrees,
e.g. $\lambda|_B=\on{deg}\lambda|_B$.

\medskip\noindent{\bf Remark 2.1.} It is possible to define the three basic
invariants for a wider variety of families. In particular, dropping the
assumption of smoothness of the general fiber roughly means that the base
curve $B$ is contained entirely in the boundary locus of 
$\overline{\mathfrak M}_g$. Since such families are not discussed in
our paper, we shall not give here these definitions. The existence,
however, of such invariants for any one-parameter family of stable curves
will follow from the description of $\lambda,\,\,\delta$ and $\kappa$ as 
``global'' classes  in $\on{Pic}_{\mathbb Q}\overline{\mathfrak M}_g$ (cf.
Sect.~\ref{linebundles}).

\medskip\noindent{\bf Remark 2.2.} It is also possible to consider families
whose special members are not stable, e.g. cuspidal, tacnodal and other
types of singular curves. One reduces to the above cases by applying 
{\it semistable reduction} (cf.~\cite{FM}).

\subsection{The line bundles $\lambda,\,\,\delta$ and $\kappa$ in
Pic$_{\mathbb Q}\overline{\mathfrak M}_g$}
\label{linebundles} Another way to
interpret the classes $\lambda|_B,\,\,\delta|_B$ and $\kappa|_B$ is to
think of them as rational divisor classes on $\overline{\mathfrak{M}}_g$. In fact,
Mumford shows that such invariants, defined for {\it
any} proper flat family $X\rightarrow S$ and natural under base change,
induce line bundles in Pic$_{\mathbb Q}\overline{\mathfrak M}_g$. Here follows
a rough sketch of the argument (cf.~\cite{M2}).

\smallskip
Consider $\on{Hilb}^{p(x)}_r$, the Hilbert scheme
parametrizing all closed subschemes of ${\mathbf P}^r$ with Hilbert polynomial
$p(x)=dx-g+1$ for some $d=2n(g-1)\gg 0$ and $r=d-g$. Let
$\cal H\subset \on{Hilb}^{p(x)}_r$ be the locally closed smooth subscheme
of $n$-canonical stable curves of genus $g$. Then $\overline{\mathfrak M}_g$ is
the GIT-quotient of $\cal H$ by  ${\on{PGL}}_r$. Let 
\[\rho:\cal H\rightarrow \overline{\mathfrak M}_g=\cal H/{\on{PGL}}_r\]
be the natural surjection, and let $(\on{Pic}\cal H)^{{\on{PGL}}_r}$ be
the subgroup of isomorphism classes of line bundles on $\cal H$ invariant 
under the action of ${\on{PGL}}_r$. 

Consider also $\on {Pic}_{\on{fun}}\overline{\mathfrak M}_g$,
the group of line bundles on the {\it moduli functor}. An element $L$ of 
$\on {Pic}_{\on{fun}}\overline{\mathfrak M}_g$ consists of the following
data: for any proper flat family $f:X\rightarrow S$ of stable curves
a line bundle $L_S$ on $S$ natural under base change. Two such elements 
are declared isomorphic under the obvious compatibility conditions.

Naturally, a line bundle on $\overline{\mathfrak M}_g$ gives rise by
pull-back to a line bundle on the moduli functor. In fact,
this map is an inclusion with a torsion cokernel, and 
$\on{Pic}_{\on{fun}}\overline{\mathfrak M}_g$ is torsion free and isomorphic
to $(\on{Pic}\cal H)^{{\on{PGL}}_r}$:
\[\on{Pic}\overline{\mathfrak M}_g\stackrel{\rho^*}{\hookrightarrow}
 \on{Pic}_{\on{fun}}\overline{\mathfrak M}_g\cong (\on{Pic}\cal H)^{{\on{PGL}}_r}.\]
Hence we may regard all these groups as sublattices of
$\on{Pic}_{\mathbb Q}\overline{\mathfrak M}_g$. In particular,
\[\on{Pic}_{\on{fun}}\overline{\mathfrak M}_g\otimes {\mathbb Q}\cong
\on{Pic}_{\mathbb Q}\overline{\mathfrak M}_g,\]
and any line bundle on the moduli functor can be thought of as a rational
class on $\overline{\mathfrak M}_g$. These identifications allow us to make 
the following

\medskip\noindent{\bf Definition 2.1.}
In $\on{Pic}_{\mathbb Q} \overline{\mathfrak M}_g$ we define the line bundles
$\lambda,\,\,\kappa$ and $\delta$ by
\[\lambda=\det\pi_*(\omega_{\cal {C}/\cal H}),\,\, \kappa=\pi_*c_1(\omega
_{\cal {C}/\cal H})^2,\,\,
\delta=\cal{O}_{\cal H}(\Delta\cal H),\]
where ${\cal{C}}\subset \cal H\times {\mathbb P}^r$ is the universal curve 
over $\cal H$, $\pi:\cal C \rightarrow \cal H$ is the projection map,
$\omega_{\cal{C}/\cal H}$ is the relative dualizing sheaf of $\pi$, and 
$\Delta\cal H\subset \cal H$ is the divisor of singular curves on $\cal H$.

\medskip
As defined, $\lambda,\kappa$ and $\delta$ lie in 
$\on{Pic}_{\mathbb Q}\overline{\mathfrak M}_g$, and as such they are only 
{\it rational} Cartier divisors on $\overline{\mathfrak M}_g$. In
\cite{MHE} one can find examples where $\lambda$  does {\it not} 
descend  to a line bundle on $\overline{\mathfrak M}_g$.
On the other hand, it is obvious from which divisor on $\overline{\mathfrak M}_g$
our $\delta$ comes: $\delta=\cal O_{\overline{\mathfrak M}_g}(\Delta)$, where
$\Delta$ denotes the divisor on 
$\overline{\mathfrak M}_g$ comprised of all singular stable curves. 
Again,  due to singularities of the total space of $\overline{\mathfrak
M}_g$, $\Delta$ is only a {\it rational}  Cartier divisor.
In fact, the only locus of $\overline{\mathfrak{M}}_g$ on which
$\lambda,\,\,\delta$ and $\kappa$ are necessarily {\it integer} divisor classes
is $(\overline{\mathfrak{M}}_g)_0$ - the locus of automorphism-free curves.

\smallskip
We can further define the {\it boundary} classes $\delta_0,
\delta_1,...,\delta_{[\frac{g}{2}]}$ in $\on{Pic} _{\mathbb Q}\overline
{\mathfrak M}_g$. Let $\Delta_i$ be the $\mathbb Q$--Cartier divisor on $\overline 
{\mathfrak M}_g$ whose general member is an irreducible nodal curve with
a single node (if $i=0$), or the join of two irreducible 
smooth curves of genera $i$ and $g-i$ intersecting 
transversally in one point (if $i>0$). Setting $\delta_i={\cal O}_{\overline
{\mathfrak M}_g}(\Delta_i)$, we have $\Delta=\sum_i\Delta_i$ and
$\delta=\sum_i\delta_i$.

\smallskip As the following result of Harer \cite{Ha1,Ha2} suggests that,
in order to describe the geometry of the moduli space 
$\overline{\mathfrak{M}}_g$, it will be useful to study
the divisor classes defined above, and to understand the relations 
between them.

\begin{thm}[Harer] The Hodge class $\lambda$ and
the boundary classes $\delta_0,\delta_1,...,\delta_{[\frac{g}{2}]}$ generate
$\on{Pic}_{\mathbb Q} \overline{\mathfrak M}_g$, and for $g\geq 3$ they are
linearly independent.
\end{thm} 

It is easy to recognize the restrictions of $\lambda,\,\,\delta$
and $\kappa$ to a curve $B$ in $\overline{\mathfrak{M}}_g$ as the previously
defined $\lambda|_B,\,\,\delta|_B$ and $\kappa|_B$. For example,
the restriction of $\delta$ to the base curve $B$ counts,
with appropriate multiplicities, the number of intersections of
$B$ with the boundary components $\Delta_i$ of $\overline{\mathfrak{M}}_g$.

As a final remark, applying Grothendieck Riemann-Roch Theorem (GRR) to the map
$\pi:\cal C \rightarrow \cal H$ and the sheaf $\omega_{\cal{C}/\cal H}$, 
implies the basic relation:
\begin{equation}
12\lambda=\kappa+\delta.
\label{GRR}
\end{equation}

\subsection{Slope of non-isotrivial families}
\label{slope-non-isotrivial} Let $f:X\rightarrow B$
be our family of stable curves with a smooth general member. By definition,
$\delta_B\geq 0$. Further, all
locally free quotients of the Hodge bundle $f_*(\omega_f)$  
have non-negative degrees \cite{key13}. If $X$ is a
non-isotrivial family, then $\lambda|_B>0$, and since 
the relative canonical divisor $K_{X/B}$ is nef,  
$\kappa|_B>0$ \cite{key32}. In particular, we can divide by $\lambda|_B$.

\medskip
\noindent{\bf Definition 2.2.} The {\it slope} of a non-isotrivial family
$f:X\rightarrow B$ of stable curves with a smooth general member
 is the ratio 
\[\on{slope}(X/\!_{\displaystyle{B}}):=\frac{\delta|_B}{\lambda|_B}
\cdot\]

Suppose we make a base change $B_1\rightarrow B$ of degree $d$, 
and set $X_1=X\times_{B}B_1$ to be the pull-back of
our family over the new base $B_1$ (cf.~Fig.~\ref{basechange}).
Then the three invariants on $B$ will
pull-back to the corresponding invariants on $B_1$, and
their degrees will be multipied by $d$,
e.g. $\lambda|_{B_1}=d\lambda|_B$, etc. In particular, the
slope of $X/_{\displaystyle{B}}$ will be preserved.

\setlength{\unitlength}{10mm}
\begin{figure}[t]
\begin{picture}(1.8,1.8)(-0.2,3.9)
\put(0,4){$B_1\stackrel{d}{\longrightarrow} B$}
\put(0,5.1){$X_1\,\,{\longrightarrow}\,X$}
\multiput(0.2,5)(1.25,0){2}{\vector(0,-1){0.6}}
\end{picture}
\hspace*{3mm}
\vspace*{-1mm}
\caption{Base change}
\label{basechange}
\end{figure}

\smallskip
 In view of (\ref{GRR}), restricting to the base curve $B$ yields
\begin{equation}
12\lambda|_B=\kappa|_B+\delta|_B.
\end{equation}
From the positivity conditions above, we deduce that
$0\leq \on{slope}(X/\!_{\displaystyle{B}})<12.$

\subsection{Statement of the problem and what is known}
\label{statement}
It is natural to ask whether we can find a better estimate for the slope of
$X$. The first fundamental result in this direction is the following

\begin{thm}[Cornalba-Harris, Xiao] Let $f:X\rightarrow B$ be a nonisotrivial
family with smooth general member. Then the slope of the family satisfies:
\begin{equation}
\frac{\delta|_B}{\lambda|_B}\leq 8+\frac{4}{g}\cdot
\label{8+4/g}
\end{equation}
Equality holds if and only if the general fiber of $f$ is hyperelliptic,
and all singular fibers are irreducible.
\label{CHX}
\end{thm}

Note that the upper bound is achieved only for hyperelliptic families.
Such families are of very special type since 
the hyperelliptic locus $\overline{\mathfrak I}_g$ has codimension $g-2$ in
$\overline{\mathfrak M}_g$.
On the other hand, if the base curve $B$ is
general enough, a much better estimate can be shown (cf.~\cite{MHE}):

\begin{thm}[Mumford-Harris-Eisenbud] If $B$ passes through a general
point $[C]\in \overline{\mathfrak{M}}_g$, then
\vspace*{-4mm}
\begin{equation}
\frac{\delta|_B}{\lambda|_B}\leq 6+\on{o}(\frac{1}{g})\cdot
\label{6+o(1/g)}
\end{equation}
\label{generalbound1}
\end{thm}

For example, when $g$ is odd, we can set $k=(g+1)/2$ and
define the divisor $\overline{\cal{D}}_k$ in ${\overline{\mathfrak M}}_g$ as
the closure of the $k$-sheeted covers of ${\proj}^1$:
\[\overline{\cal{D}}_k=\overline{\{C\in{\mathfrak{M}}_g\,\,|\,\,C\,\,
\on{has}\,\,{g}^1_k\}}.\]
\begin{figure}[h]
$$\psdraw{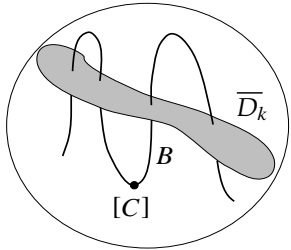}{1.5in}{1.3in}$$
\caption{General curve $B\not \subset \overline{\cal{D}}_k$}
\label{generalcurve}
\end{figure}

If our family is not entirely contained in $\overline{\cal{D}}_k$, or
equivalently, if $B$ passes through a point $[C]\not\in\overline{\cal{D}}_k$
(cf.~Fig.~\ref{generalcurve}),
\begin{equation}
\frac{\delta|_B}{\lambda|_B}\leq 6+\frac{12}{g+1}\cdot
\label{6+12/(g+1)}
\end{equation}

\smallskip\noindent
 Higher than the ``generic'' ratio can be obtained, therefore, only
for a very special type of families: those entirely contained in
$\overline{\cal{D}}_k$, and hence possessing ${g}^1_2$, ${g}^1_3$, etc.

\subsubsection{The rational Picard group of the hyperelliptic locus 
$\overline{\mathfrak{I}}_g$}
\label{rationalhyper}
In proving the maximal bound $8+4/g$, Cornalba-Harris also
describe $\on{Pic}_{\mathbb Q}\overline{\mathfrak I}_g$
by exhibiting generators and relations (cf.~\cite{CH}). 
Here we briefly discuss their result. 

\smallskip
Recall the irreducible divisors $\Delta_i$ on $\overline{\mathfrak M}_g$. 
For $i=1,...,[g/2]$, $\Delta_i$ cuts out an irreducible divisor
on $\overline{\mathfrak I}_g$, while the intersection
$\Delta_0\cap\overline{\mathfrak I}_g$ breaks
up into several components:
\[\Delta_0\cap\overline{\mathfrak I}_g=\Xi_0\cap\Xi_1\cap\cdots
\cap\Xi_{[\frac{g-1}{2}]}.\]
Set $\xi_i=\cal O_{\overline{\mathfrak I}_g}(\Xi_i)$ for the class of $\Xi_i$ in
$\overline{\mathfrak I}_g$, and retain the symbols $\lambda$ and $\delta_i$
 for their corresponding restrictions to 
$\on{Pic}_{\mathbb Q}\overline{\mathfrak I}_g$. Thus,
$\delta_i:=\cal O_{\overline{\mathfrak I}_g}(\Delta_i\cap\overline{\mathfrak I}_g)$
for all $i$. Finally note that the class $\delta_0$ is realised in
$\on{Pic}_{\mathbb Q}\overline{\mathfrak I}_g$ as the sum
\[\delta_0=\xi_0+2\xi_1+\cdots+2\xi_{[\frac{g-1}{2}]}.\]
The coefficient $2$ roughly means that
$\Delta_0$ is {\it double} along $\Xi_i$, for $i>0$.

\begin{thm}[Cornalba-Harris] The classes
$\xi_0,\cdots,\xi_{[\frac{g-1}{2}]}$ and $\delta_1,\cdots,\delta_{[\frac{g}
{2}]}$ freely generate 
$\on{Pic}_{\mathbb Q}{\overline{\mathfrak I}_g}$. The Hodge class
$\lambda\in \on{Pic}_{\mathbb Q}{\overline{\mathfrak I}_g}$ is expressed
in terms of these generators as the following linear combination:
\begin{equation}
(8g+4)\lambda=g\xi_0+\sum_{i=1}^{[(g-1)/2]}2(i+1)(g-i)\xi_i
+\sum_{j=1}^{[g/2]}4j(g-j)\delta_j.
\label{CHPic}
\end{equation}
\label{theoremCHPic}
\end{thm}
For a specific family $f:X\rightarrow B$ of hyperelliptic 
stable curves this relation reads:
\[(8g+4)\lambda|_B=g\xi_{0}|_B+\sum_{i=1}^{[(g-1)/2]}2(i+1)(g-i)\xi_{i}|_B
+\sum_{j=1}^{[g/2]}4j(g-j)\delta_{j}|_B\]
\[\Rightarrow\,\, (8+4/g)\lambda|_B\geq \xi_{0}|_B+\sum_i2\xi_{i}|_B+
\sum_j2\delta_{j}|_B=\delta|_B.\]
This yields the desired $8+4/g$
inequality for the slope of a hyperelliptic family,
and shows that the maximum can be obtained exactly when all 
$\xi_1,\cdots,\xi_{[\frac{g-1}{2}]},\delta_1,\cdots,\delta_{[\frac{g}{2}]}$
vanish on $B$. In other words, the singular fibers of $X$ belong only to
the boundary divisor $\Xi_0$, and hence are irreducible. In Appendix
we review the description of the divisors $\Xi_i$ via admissible covers, and
give an alternative proof of Theorem~\ref{theoremCHPic}.

\subsubsection{Example of a hyperelliptic family with maximal slope}
\label{example} We present here a typical example in which the 
upper bound $8+4/g$ is achieved, and show how to
calculate explicitly the basic invariants $\lambda|_B$ and $\delta|_B$
for this family.

\begin{figure}[t]
\begin{picture}(1,3)(0,3.1)
\put(-0.7,5.5){$X\,\,\hookrightarrow\,\,{\proj}^1\!
\!\times\! {\proj}^1\!\!\times\!{\proj}^1$}
\put(-0.5,5.4){\vector(1,-1){1}}
\put(1.7,5.4){\vector(-1,-1){1}}
\put(0,3.9){${\proj}^1\!\!\times\! {\proj}^1$}
\put(0.4,2.7){${\proj}^1$}
\put(0.6,3.8){\vector(0,-1){0.6}}
\end{picture}
\caption{Ratio $8+4/g$}
\label{ratio8+4/g}
\end{figure}
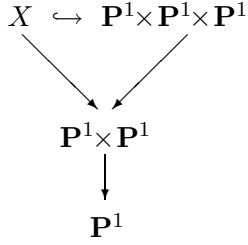

\medskip
\noindent{\bf Example 2.1.} Consider a pencil $\cal{P}$ of hyperelliptic
curves of genus $g$
on ${\proj}^1\!\times \!{\proj}^1$. Because of genus considerations,
its members must be of type $(2,g+1)$.
Our family \newline
$f\!:\!X\!\rightarrow \!{\proj}^1$ will be obtained by blowing-up 
${\proj}^1\!\times\! {\proj}^1$ at the $4(g+1)$ base points of the pencil
in order to separate its members (cf.~Fig.~\ref{ratio8+4/g}). Hence,
$\chi(X)=\chi({\proj}^1\!\times\! {\proj}^1)+4(g+1)$ for the
corresponding topological Euler characteristics. Riemann-Hurwitz formula for
the map $f$ gives a second relation:
$\chi(X)=\chi({\proj}^1)\chi(X_b)+\delta|_B$, where ${\proj}^1$
is the base $B$ and $X_b$ is the general
fiber of $X$. Putting together, $\delta|_B=8g+4.$

The total space
of $X$ is a divisor on ${\proj}^1\!\!\times\! {\proj}^1\!\!\times\!
{\proj}^1$ of type
$(2,g+1,1)$, and the map $f\!:\!X\!\rightarrow \!{\proj}^1$ is the
restriction to $X$ of the third projection $\pi_3\!:\!{\proj}^1\!\!\times\!
{\proj}^1\!\!\times\!{\proj^1}\!\rightarrow \!{\proj}^1$. 
Using standard methods,
we compute $h^0( (f_*(\omega_f))(-2))=0$. From the positivity of all free
quotients of the Hodge bundle on ${\proj}^1$, $f_*(\omega_f)$
splits as a direct sum $\bigoplus_{i=1}^g
{\cal O}_{{\proj}^1}(a_i)$ for some $a_i>0$.
Then, for $f_*(\omega_f)(-2)
=\bigoplus_i\cal{O}_{\proj^1}(a_i-2)$
to have no sections, all $a_i$'s must be at most $1$. Finally,
\begin{equation*}
f_*(\omega_f)=\bigoplus_{i=1}^g{\cal O}_{{\proj}^1}(+1),\,\,\lambda|_B=g,
\,\,\on{and}\,\,\frac{\delta|_B}{\lambda|_B}=8+\frac{4}{g}
\cdot
\end{equation*}

\subsubsection{The Trigonal Locus $\overline{\mathfrak{T}}_g$}
\label{trigonallocus} In a similar vein
as in the above example, we consider pencils of trigonal curves on ruled 
surfaces, and obtain the slope $7+6/g$. It is somewhat reasonable to expect
that this is the maximal ratio. Recall that a bundle $\cal{E}$ on a curve
$B$ is {\it semistable} if for any proper subbundle $\cal{F}$, we have 
$q(\cal{F})\leq q(\cal{E})$, where $q$ is 
the quotient of the degree and the rank
of the corresponding bundle. Following Xiao's approach in the proof of  
Theorem~\ref{CHX}, Konno shows  
that for non-hyperelliptic fibrations of genus $g$ with semistable
Hodge bundle $f_*\omega_{f}$ (cf.~\cite{HN,key5}):
\begin{equation}
\frac{\delta|_B}{\lambda|_B}\leq 7+\frac{6}{g}\cdot
\label{7+6/g}
\end{equation}
As for any trigonal families, he establishes the inequality (cf.~\cite{key6}):
\begin{equation}
\frac{\delta|_B}{\lambda|_B}\leq
\frac{22g+26}{3g+1}\sim 7\frac{1}{3}+\on{o}
(\frac{1}{g})\cdot
\end{equation}
Examples of trigonal families achieving this ratio were not found, which
suggested that this bound might be too big. On the other hand, in trying to
disprove the smaller bound $7+6/g$, we naturally arrived at counterexamples
pointing to a third intermediate ratio (cf.~Theorem~\ref{maximal bound2}): 
\begin{equation}
\frac{36(g+1)}{5g+1}\sim 7\frac{1}{5}+\on{o}
(\frac{1}{g})\cdot
\end{equation}
The difference between the last two estimates may seem negligible, but this
would not be so when both $\lambda|_B$ and $\delta|_B$ become large and
we attempt to bound $\lambda|_B$ from below by $\delta|_B$. What is more
important, the second ratio is in fact {\it exact}, and we give
equivalent conditions for it to be achieved (cf.~Sect.~\ref{whenmaximal},
\ref{Maroni-maximal}). This maximal bound confirms Chen's result
for genus $g=4$ in \cite{Chen}.

\smallskip
The reader may ask why the ``generic'' examples for the maximum in the 
hyperelliptic case fail to provide also the maximum in the trigonal
case. As we noted in the Introduction, the answer is  closely related to 
the so-called {\it Maroni} locus in $\overline{\mathfrak{T}}_g$.  
More precisely, if ${\mathbf F}_k={\proj}({\cal O}_{{\proj}^1}\oplus
{\cal O}_{{\proj}^1}(k))$ denotes the corresponding rational ruled surface,
a general curve $C$ embeds in ${\mathbf F}_0$ is
$g$ is even, and in ${\mathbf F}_1$ if $g$ is odd. The Maroni locus consists
of those curves that embed in ${\mathbf F}_k$ with $k\geq 2$. The number
$k/2$ is referred to as the {\it Maroni invariant} of $C$.  
In these terms, the examples of pencils of trigonal curves on ${\mathbf F}_0$
and ${\mathbf F}_1$ have the lowest possible constant Maroni invariant,
and we shall see that the maximum bound can be obtained only for families 
entirely contained in the Maroni locus, and having very high Maroni 
invariants.

\medskip
The ``semistable'' bound $7+6/g$ appears in Theorem~\ref{7+6/g Bogomolov2},
where we give instead a sufficient {\it Bogomolov-semistability} condition
$4c_2(V)-c_1^2(V)\geq 0$ for a canonically associated to $X$ vector bundle
$V$ of rank 2. The rational Picard group of $\overline{\mathfrak{T}}_g$
is described in terms of generators and relations in
Section~\ref{generators}, providing thus in the trigonal case
an analog of Theorem~\ref{theoremCHPic}. Note the apparent similarity 
of the coefficients $\widetilde{c}_{k,i}$ of the trigonal boundary classes 
and the coefficients of the hyperelliptic boundary classes.
This is not coincidental. In fact, the $\widetilde{c}_{k,i}$'s are
in a sense the ``smallest'' coefficients that could have been
associated to the corresponding classes $\delta_{k,i}$
(cf.~Fig.~\ref{Delta-k,i}): they are symmetric with respect
to the two genera of the components in the general member of $\delta_{k,i}$.
A crucial role in the proof of Theorem~\ref{Pic trigonal} is played
by the interpretation of the above Bogomolov semistability condition
in terms of the Maroni locus it $\overline{\mathfrak{T}}_g$ (cf.~Sect.~\ref
{interpretation}).

\medskip
\subsection{The idea of the proof}
\label{idea} 

Let $f:X\rightarrow B$ be a family of stable curves, whose general member
$X_b$ is a smooth trigonal curve.
By definition, $X_b$ is a triple cover of ${\proj}^1$. We would like
to study how this triple cover varies as $X_b$ moves in the family $X$.
Thus, it would be desirable to represent $X$, by analogy with
$X_b$, as a triple cover of a ruled
surface $Y$, comprised by the image lines ${\proj}^1$. Unfortunately, due
to existence of hyperelliptic and other
special singular fibers, this is not always possible. 

\setlength{\unitlength}{10mm}
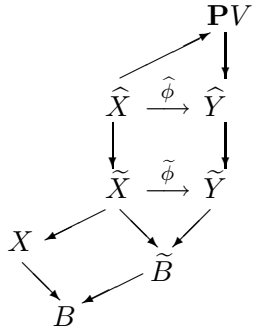
\begin{figure}[h]
\begin{picture}(3,4.9)(-1,2)
\put(0,4){$\widetilde{X}\,\stackrel{\widetilde{\phi}}
{\longrightarrow}\, \widetilde{Y}$}
\put(0,5.1){$\widehat{X}\,\stackrel{\widehat{\phi}}{\longrightarrow}\,
 \widehat{Y}$}
\multiput(0.2,3.85)(-1.3,-0.7){2}{\vector(1,-1){0.5}}
\put(1.4,3.85){\vector(-1,-1){0.5}}
\multiput(0.5,3)(-0.5,0.85){2}{\vector(-2,-1){0.8}}
\put(0.6,2.9){$\widetilde{B}$}
\put(-1.3,3.3){$X$}
\put(-0.7,2.3){$B$}
\multiput(0.1,5)(1.5,0){2}{\vector(0,-1){0.6}}
\put(1.6,6.2){\vector(0,-1){0.6}}
\put(1.35,6.3){${\proj}V$}
\put(0.2,5.6){\vector(2,1){1.2}}
\end{picture}
\vspace*{-3mm}
\caption{Basic construction}
\label{Basic construction idea}
\end{figure}

 \subsubsection{The basic construction} The ``closest''
model of such a triple cover can be obtained after a finite number of 
birational transformations on $X$, and a possible base change over the base 
$B$. This way we construct a {\it
quasi-admissible} cover $\widetilde{\phi}:{\widetilde{X}}
\rightarrow {\widetilde{Y}}$ over a new base
${\widetilde{B}}$ (cf.~Prop.~\ref{propquasi}). Here ${\widetilde{Y}}$ is a
{\it birationally} ruled surface over $\widetilde{B}$ with reduced, but non
necessarily irreducible, special fibers: $\widetilde{Y}$ allows for 
{\it pointed stable} rational fibers, i.e. trees of ${\proj}^1$'s with
points marked in a certain (stable) way.

The map $\widetilde{\phi}$ expresses any fiber
$\widetilde{X}_ b$ 
as a triple quasi-admissible cover of the corresponding {\it pointed 
stable} rational curve $\widetilde{Y}_b$. To calculate
effectively our invariants $\lambda,\delta$ and $\kappa$, we need that
$\widetilde{\phi}$ be {\it flat}, 
which could force a few additional blow-ups on $\widetilde{X}$
and $\widetilde{Y}$. We end up with a flat proper triple
cover $\widehat{\phi}:
\widehat{X}\rightarrow \widehat{Y}$, where certain fibers of 
$\widehat{X}$ and $\widehat{Y}$
are allowed to be {\it non-reduced}: these are the scheme-theoretic
preimages under the blow-ups on $\widetilde{X}$ and $\widetilde{Y}$.
We call such covers $\widehat{\phi}$ {\it effective}.

\smallskip\quad
We observe next that any smooth trigonal curve $C$ can
be naturally embedded in a ruled surface ${\mathbf F}_k$ over $B$. If
$\alpha:C\rightarrow {{\proj}^1}$ is the corresponding triple cover, there
is an exact sequence of locally free sheaves on ${{\proj}^1}$:
\begin{equation*}
0\rightarrow {V}\rightarrow {\alpha}_*{\cal O}_{C}\stackrel
{\on{tr}}{\rightarrow}{\cal O}_{{\proj}^1}\rightarrow 0.
\end{equation*}
The projectivization $\proj V$ of the rank 2 vector bundle $V$ is the
ruled surface ${\mathbf F}_k$.

\smallskip
This construction can be extended as $C$ moves in the effective cover
$\widehat{\phi}:\widehat{X}\rightarrow \widehat{Y}$. The flatness of
$\widehat{\phi}$ forces the pushforward ${\phi}_*{\cal O}_{\widehat{X}}$
to be a locally free sheaf of rank 3 on $\widehat{Y}$, and the finiteness
of $\widehat{\phi}$ ensures the existence of a {\it trace map}
$\on{tr}:{\phi}_*{\cal O}_{\widehat{X}}\rightarrow {\cal O}_{\widehat{Y}}$.
Again, the kernel $V$ of $\on{tr}$ 
is the desired rank 2 vector bundle on $\widehat{Y}$, in
whose projectivization, ${\proj}V$, we embed $\widehat{X}$
(cf.~Fig.~\ref{Basic construction idea}).

\subsubsection{Chow Rings Calculations}
We can now use the relations in the Chow rings of ${\mathbb{A}}({\mathbf P}V)$,
$\mathbb{A}\widehat{Y}$ and $\mathbb{A}\widehat{X}$ to calculate the invariants
$\lambda_{\widehat{X}}$ and $\delta_{\widehat{X}}$, 
appropriately defined for the new family 
$\widehat{X}\rightarrow {\widetilde{B}}$ of semistable and occasionally
non-reduced fibers. Then, of course, we translate 
$\lambda_{\widehat{X}}$ and $\delta_{\widehat{X}}$ into $\lambda_{{X}}$ 
and $\delta_{{X}}$ with the necessary adjustments from the
birational transformations on $X$ and the base change on $B$. 
We compare the resulting expressions to obtain a relation among
$\lambda_X$ and $\delta_X$.  

\subsubsection{Boundary of the Trigonal Locus}
As we vary the base curve $B$ inside $\overline{\mathfrak{T}}_g$, 
we actually obtain a relation among the restrictions of $\lambda$ and $\delta$
in $\on{Pic}_{\mathbb{Q}}\overline{\mathfrak{T}}_g$, rather than just
among $\lambda|_B=\!\lambda_X$ and $\delta|_B=\!\delta_X$ in 
$\on{Pic}B$.

\smallskip
{\it In terms of what} have we thus represented and linked
$\lambda|_{\overline{\mathfrak{T}}_g}$ and $\delta|_{\overline
{\mathfrak{T}}_g}$?
To answer this question, we need first to understand the
boundary divisors of the trigonal locus $\overline{\mathfrak{T}}_g$.
As we shall see, there
are seven types of such divisors, denoted by
$\Delta{\mathfrak{T}}_{0}$ and $\Delta{\mathfrak{T}}_{k,i}$ for $k=1,...,6$. 
Each type is
determined by the specific geometry of its general member. For example,
$\Delta{\mathfrak{T}}_0$ is the closure of all irreducible trigonal curves with
one node, while $\Delta{\mathfrak{T}}_{2,i}$ corresponds to joins in two points
of a trigonal and a hyperelliptic curve with genera $i$ and $g-1-i$,
respectively (cf.~Fig.~\ref{Delta-k,i}). 
Naturally, we derive an expression for the restriction of 
the divisor class $\delta\in\on{Pic}_{\mathbb{Q}}\overline{\mathfrak{M}}_g$ to 
$\overline{\mathfrak{T}}_g$:
\begin{equation*}
\delta|_{\displaystyle{\overline{\mathfrak{T}}_g}}=\delta_0+
\sum_{i=1}^{\scriptscriptstyle{[(g-2)/2]}}3\delta_{1,i}
+\sum_{i=1}^{\scriptscriptstyle{g-2}}2\delta_{2,i}
+\sum_{i=1}^{\scriptscriptstyle{[g/2]}}\delta_{3,i}
+\sum_{i=1}^{\scriptscriptstyle{[(g-1)/2]}}3\delta_{4,i}+
\sum_{i=1}^{\scriptscriptstyle{g-1}}\delta_{5,i}
+\sum_{i=1}^{\scriptscriptstyle{[g/2]}}\delta_{6,i}.
\label{delta}
\end{equation*}
Here $\delta_0$ and $\delta_{k,i}$ are the divisor classes of 
$\Delta{\mathfrak{T}}_0$ and
$\Delta{\mathfrak{T}}_{k,i}$ in $\on{Pic}_{\mathbb{Q}}\overline{\mathfrak{T}}_g$.

\subsubsection{Relations among $\lambda$ and $\delta$}
For a fixed family $X\rightarrow B$ with a smooth trigonal general member,
we establish a relation among the Hodge class $\lambda|_B$, the
boundary classes $\delta_{k,i}|_B$, and the Bogomolov quantity
$4c_2(V)-c_1^2(V)$ for the associated vector bundle $V$:
\begin{equation}
(7g+6)\lambda|_B=g\delta_0|_B+\sum_{k,i}\widetilde{c}_
{k,i}\delta_{k,i}|_B+\frac{g-3}{2}(4c_2(V)-c_1^2(V)).
\label{tobelifted}
\end{equation}
The polynomial coefficients $\widetilde{c}_{k,i}$ are
comparatively larger than the corresponding coefficients of the
boundary divisors in the expression for
$\delta|_{\displaystyle{\overline{\mathfrak{T}}_g}}$.
As a result, we rewrite (\ref{tobelifted}) as
\begin{equation}
(7g+6)\lambda|_B=g\delta|_B+\cal{E}|_B+\frac{g-3}{2}(4c_2(V)-c_1^2(V)),
\label{E-argument}
\end{equation}
where $\cal{E}$ is an effective combination of the boundary classes on 
$\overline{\mathfrak{T}}_g$. In particular, if $V$ is Bogomolov semistable,
the slope satisfies (cf.~Theorem~\ref{7+6/g Bogomolov2}):
\begin{equation}
\on{slope}(X/_{\displaystyle{B}})\leq 7+\frac{6}{g}.
\label{idea7+6/g}
\end{equation}
Further, we describe $\on{Pic}_{\mathbb{Q}}\overline{\mathfrak{T}}_g$
as generated freely by the restriction $\lambda|_{\overline{\mathfrak{T}}_g}$
and the boundary classes of $\overline{\mathfrak{T}}_g$. In the even genus
$g$ case, we can replace $\lambda|_{\overline{\mathfrak{T}}_g}$ by
a geometrically defined class $\mu$, corresponding to the so-called
Maroni divisor in $\overline{\mathfrak{T}}_g$. This, of course, means
that the Hodge class $\lambda|_{\overline{\mathfrak{T}}_g}$ must be
some linear combination of the boundary classes and $\mu$. The Bogomolov 
quantity is interpreted as 
\[4c_2(V)-c_1^2(V)=4\mu|_B+0\cdot \delta_0|_B+\sum_{k,i}\alpha_{k,i}
\delta_{k,i}|_B,\]
which in turn ``lifts'' (\ref{tobelifted}) to the wanted relation in
$\on{Pic}_{\mathbb{Q}}\overline{\mathfrak{T}}_g$:
\begin{equation*}
(7g+6)\lambda|_{\overline{\mathfrak{T}}_g}=g\delta_0+
\sum_{k,i}\widehat{c}_{k,i}\delta_{k,i}+2(g-3){\mu}.
\end{equation*}
We have not yet computed explicitly all coefficients $\widehat{c}_{k,i}$.
In the cases which we have completed ($\Delta_{0}\mathfrak{T}_g$
and $\Delta_{1,i}\mathfrak{T}_g$), these coefficients  
turn out again  sufficiently large so that 
we can repeat the argument in (\ref{E-argument}). 
Thus, if $X$ has
at least one non-Maroni fiber, and its singular fibers belong to 
$\Delta_{0}\mathfrak{T}_g\cup\Delta_{1,i}\mathfrak{T}_g$, 
then $\mu|_B\geq 0$, and hence the stronger bound of (\ref{idea7+6/g}) holds
 (cf.~Prop.~\ref{Maroni inequality} and Conj.~\ref{Maroni-conj}).

\subsubsection{Maximal Bound}
Since the Bogomolov semistability condition 
$4c_1^2(V)-c_2(V)\geq 0$ is not always satisfied, the above discussion 
shows that $7+6/g$ is {\it not} the maximal bound for the
slope of trigonal families, Therefore, we need another, 
more subtle, estimate. The expressions for $\lambda|_B$ and $\delta|_B$ 
suggest that any maximal
bound would be equivalent to an inequality involving $c_1^2(V)$,
$c_2(V)$, and possibly some other invariants. We construct a specific
divisor class $\eta$ on $\widetilde{X}$, for which the {\it Hodge Index}
theorem implies $\eta^2\leq 0$, and we translate this into
$9c_2(V)-2c_1^2(V)\geq 0$  (cf.~Prop.~\ref{genindex}). 
We notice that the only reasonable
way to replace Bogomolov's condition $4c_1^2(V)-c_2(V)\geq 0$
by the newly found inequality is by subtracting the following quantities:
\begin{equation*}
36(g+1)\lambda|_B-(5g+1)\delta|_B= \cal{E}^{\prime}|_B+(g-3)\big(9c_2(V)-
2c_1^2(V)\big),
\label{maximum1}
\end{equation*}
so that the ``left-over'' linear combination of boundary divisors
$\cal{E}^{\prime}$ is again effective (cf.~Theorem~\ref{maximal relation2}). 
Hence, we conclude that for {\it
all} trigonal families:
\[\on{slope}(X/_{\displaystyle{B}})\leq \frac{36(g+1)}{5g+1}
\cdot\]

\subsection{The organization of the paper}
\label{organization}

\begin{figure}[h]
$$\psdraw{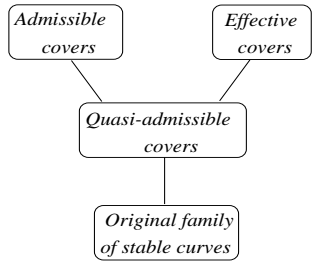}{1.6in}{1.35in}$$
\caption{Types of covers}
\label{stages}
\end{figure}

The presentation of the {\it Basic Construction} is done in several stages.
Fig.~\ref{stages} shows schematically the connection between the
three types of covers, admissible, quasi-admissible and effective, in
relation to the original family $X\rightarrow B$ of stable curves.

\vspace*{4mm}

We start in Section~\ref{hurwitz} by introducing a compactification
 $\overline{\cal{H}}_{d,g}$ of the Hurwitz scheme, parametrizing {\it
admissible} $d$-uple covers of stable pointed rational curves. Using its
coarse moduli properties, we show in Section ~\ref{admissible} the existence
of admissible covers of surfaces $X^a\rightarrow Y^a$ associated to the
original family $f\!:\!X\!\rightarrow \!B$. Next we modify these covers to 
{\it quasi-admissible} covers $\widetilde{\phi}:\widetilde{X}
\rightarrow \widetilde{Y}$ (cf.~Prop.~\ref{propquasi}),
and further to {\it effective} covers $\widehat{\phi}:\widehat{X}\rightarrow
\widehat{Y}$ in order to resolve the technical difficulties arising from the
non-flatness of $\widetilde{\phi}$ (cf.~Sect.~\ref{effectivecovers}).

\bigskip
We devote Section~4
to the study of the boundary components of the trigonal
locus $\overline{\mathfrak{T}}_g$ inside the moduli space
$\overline{\mathfrak{M}}_g$, and express the restriction
$\Delta|_{\overline{\mathfrak{T}}_g}$ as a linear combination of the
boundary divisors (cf.~Prop.~\ref{divisorrel}).
In Section~6
we complete the Basic Construction by embedding the effective
cover $\widehat{X}$ in a rank 1 projective bundle ${\proj}V$ over
$\widehat{Y}$. 

\medskip
For convenience of the reader, the proofs of the maximal $36(g+1)/(5g+1)$
and the semistable $7+6/g$ bounds are presented first in the special, 
but fundamental case when the original family
$f\!:\!X\!\rightarrow \!B$ 
is already an effective triple cover of a ruled surface
$Y$ (cf.~Sect.~7). The discussion results in finding the 
coefficients of $\delta_0$ in two different expressions of
$\lambda|_{\overline{\mathfrak{T}}_g}$, but, as it turns out, the knowledge of
these coefficients is enough to determine the desired two bounds. We refer
to this as the {\it global} calculation.
The Hodge Index Theorem and Nakai-Moishezon criterion
on $X$ complete the global calculation in Sect.~\ref{indextheorem}. 
A discussion of maximal bound examples can be found in
Section~\ref{whenmaximal}.

\medskip
The {\it local} calculations in Sections~8-10 
compute the contributions of the other boundary \vspace*{-1mm}classes
$\delta_{k,i}$, and express $\lambda|_{\overline{\mathfrak{T}}_g}$ in terms
of these contributions and the Chern classes of the rank 2 vector bundle
$V$ on $\widehat{Y}$. For clearer exposition, the proofs of the
two bounds are shown first for a {\it general} base curve $B$ (i.e.
$B$ intersects transversally the boundary components in general points),
and then in Section~11 the results are extended to {\it any}
base curve $B$. We develop the necessary notation and techniques for the
local calculations in Section ~\ref{conventions}.

\medskip
Section~12 discusses the relation between the Bogomolov
semistability condition and the Maroni locus, and describes the structure of
$\on{Pic}_{\mathbb{Q}}\overline{\mathfrak{T}}_g$. In Section~\ref{Maroni-maximal}
we give another interpretation of the conditions for the maximal bound.

\medskip
We present further results and conjectures for $d$-gonal families in
Section~13. In the Appendix, we
give another proof of the $8+4/g$ bound in
the hyperelliptic case and show an application of the
maximal trigonal bound to the study of the discriminant
locus of certain triple covers.

\bigskip 
\section*{\hspace*{1.9mm}3. Quasi-Admissible Covers of Surfaces}

\setcounter{section}{3} 
\setcounter{subsection}{0} 
\setcounter{subsubsection}{0} 
\setcounter{lem}{0}
\setcounter{thm}{0}
\setcounter{prop}{0}
\setcounter{defn}{0} 
\setcounter{cor}{0} 
\setcounter{conj}{0} 
\setcounter{claim}{0} 
\setcounter{remark}{0}
\setcounter{equation}{0}
\label{quasi-admissible}

We first review briefly 
the theory of admissible covers. For more details, we refer
the reader to \cite{MHE,HM}.

\subsection{The Hurwitz scheme $\overline{\cal H}_{d,g}$}
\label{hurwitz}
\begin{figure}[h]
\begin{picture}(5,3.5)(-0.8,2.2)
\put(0,4){${\cal H}_{d,g}\hookrightarrow \overline{\cal H}_{d,g}$}
\put(0.4,3.85){\vector(1,-1){0.9}}
\put(1.9,3.85){\vector(1,-1){0.9}}
\put(0.4,4.2){\vector(1,1){0.9}}
\put(1.9,4.2){\vector(1,1){0.9}}
\put(2.4,4.5){${pr}_1$}
\put(2.4,3.5){${pr}_2$}
\put(1.1,2.5){${\mathfrak P}_{0,b}\hookrightarrow 
\overline{\mathfrak P}_{0,b}$}
\put(1.1,5.2){${\mathfrak M}_g\hookrightarrow \overline{\mathfrak M}_g$}
\end{picture}
\caption{Hurwitz scheme}
\label{Hurwitz figure}
\end{figure}
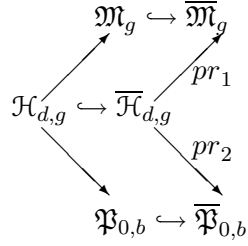

Let ${\cal H}_{d,g}$ be the {\it small Hurwitz scheme}
 parametrizing the pairs
$(C,\phi)$, where $C$ is a smooth curve of genus $g$ and $\phi:C\rightarrow
{\proj}^1$ is a cover of degree $d$, simply branched over $b=2d+2g-2$
distinct points. Since $C\in {\mathfrak M}_g$, there is a
natural map ${\cal H}_{d,g}\rightarrow {\mathfrak M}_g$, whose
image contains an open dense subset of ${\mathfrak M}_g$. The theory of
admissible covers provides the commutative diagram in 
Fig.~\ref{Hurwitz figure}.

There ${\mathfrak P}_{0,b}$ (resp. $\overline{\mathfrak P}_{0,b}$)
is the moduli space of $m$-pointed ${\proj}^1$'s (resp. of stable 
$m$-pointed rational curves), and
$\overline{\cal H}_{d,g}$ is a compactification of the Hurwirz
scheme. The points of $\overline{\cal H}_{d,g}$ correspond to triples 
$(C,(P;p_1,...,p_m),\phi)$,
where $C$ is a connected reduced nodal curve of genus $g$, 
$(P;p_1,...,p_m)$ is a stable
$m$-pointed rational curve, and $\phi:C\rightarrow P$ is a 
so-called {\it admissible cover}.

\medskip
\noindent{\bf Definition 3.1.} Given the curves $C$ and $P$ as above, an
{\it admissible cover} $\phi:C\rightarrow P$ is a regular map
satisfying the following conditions:

\smallskip
(A1) $\phi^{-1}(P_{\on{sm}})=C_{\on{sm}}$ and $\phi:C_{\on{sm}}
\rightarrow P_{\on{sm}}$
is simply branched over the distinct points $p_1,...,p_b\in P_{\on{sm}}$;

(A2) for every $q\in C_{\on{sing}}$ lying over a node $p\in P$, the two
branches through $q$ map with the same ramification index
to the two branches through $p$.

\begin{figure}[h]
$$\psdraw{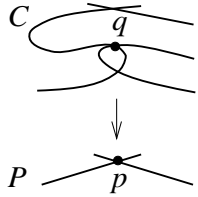}{1in}{1in}$$
\caption{Admissible model}
\label{nonstable}
\end{figure}
\smallskip
Note that $C$ is not necessarily a stable curve, but contracting its
destabilizing rational chains yields the corresponding stable
curve $pr_1(C)\in \overline{\mathfrak M}_g$. In such a case, we say that
$C$ is the ``admissible model'' for $pr_1(C)$ (cf.~Fig.~\ref{nonstable}). 
Harris-Mumford have shown that
the compactification $\overline{\cal H}_{d,g}$  is in fact
 a {\it coarse moduli space} for the admissible covers $\phi:C\rightarrow P$. 

\begin{figure}[h]
$$\psdraw{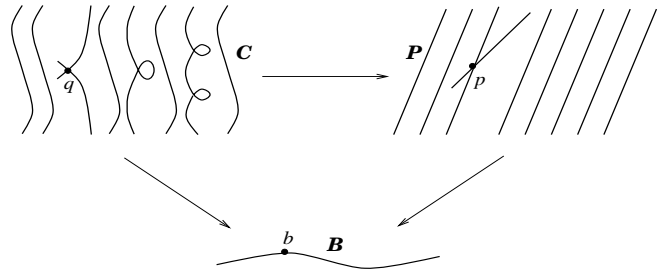}{3.4in}{1.4in}$$
\hspace*{4.5mm}\vspace*{-5mm}
\caption{Admissible family}
\label{admfamilies}
\end{figure}

\subsection{Local properties of admissible covers} 
\label{localproperties}
When we vary the admissible covers of curves in families, the
local structure of the corresponding total spaces becomes apparent.
Let $\phi:\cal C\rightarrow \cal P$ be a proper flat family (over a scheme
$\cal B$) of admissible covers of curves (cf.~Fig.~\ref{admfamilies}). 
Assume that $\phi$ is \'{e}tale
everywhere except over the nodes of the fibers of $\cal P/_
{\textstyle{\cal B}}$,
and except over some sections $\sigma_i:\cal B\rightarrow \cal C$ and their
images $\omega_i:\cal B\rightarrow \cal P$: there $\phi$ is simply branched
along $\sigma_i$ over $\omega_i$ for all $i$. If $q\in {\cal C}_b$ is a point
lying above a node $p\in {\cal P}_b$ for some $b\in \cal B$, then
$\cal C_b$ has a node at $q$, and locally analytically  we can describe
$\cal C,\cal P$ and $\phi$ near $q$ and $p$ by:

\[\left\{\begin{array}{lll}
\cal C: & xy=a, &x,y\,\,\on{generate}\,\,\widehat{\mathfrak m}_{q,\cal
C_b},\,\, a\in \widehat{\cal O}_{b,\cal B},\\
\cal P: & uv=a^n, &u,v\,\,\on{generate}\,\,\widehat{\mathfrak m}_{p,\cal
P_b},\\
\phi:   & u=x^n,v=u^n.
\end{array}\right.\]

\smallskip
One can see that $n$ is the index of ramification of $\phi$ at $q$,
and that fiberwise $\cal C_b\rightarrow \cal P_b$ is
an admissible cover (of curves). From now on, by 
{\it admissible covers} we mean, more generally,
families $\cal C\rightarrow \cal P$ over $\cal B$ with the above description. 

\smallskip
The local properties of the admissible cover $\phi:\cal C\rightarrow \cal P$
over the nodes in $\cal P_b$ forces singularities on the total spaces
of $\cal C$ and $\cal P$. Since we will be interested only in the cases when
the base $\cal B$ is a smooth projective curve $B$  and the general
fiber of $\cal C$ 
is smooth, we can always pick a generator $t$ for $\widehat{\cal O}_{b,B}$,
and express $a=t^l$ for some $l\in {\mathbb N}$.

\begin{figure}[h]
$$\psdraw{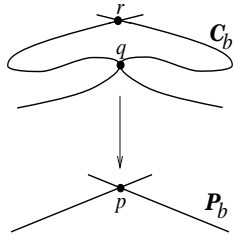}{1.2in}{1.2in}$$
\caption{Singularity of $\cal C$}
\label{singular}
\end{figure}

\noindent{\bf Example 3.1.} Let the triple
admissible cover $\phi\!:\!\cal C\!\rightarrow\!\cal P$ 
contain the fiber $\cal C_b$ as in Fig.~\ref{singular}. At $q$,  
${\cal C}$ is given by 
$\,\,xy=t^{l}$, and at $p$, $\cal P$ is given by 
$uv=t^{2l}$, where $u=x^2,\,\,v=y^2$. This
forces  at $r$ the local equation 
$xy=t^{2l}$ ($u=x,\,\,v=y$). 
Even if $\cal C$ is smooth at $q$ ($l=1$), $\cal C$
and $\cal P$ will be singular at $r$ and $p$, respectively
($xy=t^2,\,\,uv=t^2$). Compare this with the non-flat cover of ramification
index 1 in Fig.~\ref{mult4}.

\smallskip
Recall that a rational double point $s$ on a surface $S$ is of type
$A_{l-1}$ if locally analytically $S$ is given at $s$ by the
equation $xy=t^l$. Thus, $r$ and $p$ above are rational double points on
$\cal C$ and $\cal P$, respectively, of type $A_{l-1}$.

\medskip
\noindent{\bf Remark 3.1.}
In the sequel, we use the fact that the projection 
$pr_1\!:\!\overline{\cal H}_{d,g}\!\rightarrow
\!\overline{\mathfrak P}_{0,b}$ is a {\it finite} map. From the weak valuative
criterion for properness, this means that given a family of admissible
covers $\phi:{\cal C}^*\rightarrow {\cal P}^*$ over the punctured disc
${\on {Spec}}\,{\mathbb C}((t))$, there is some $n\in {\mathbb N}$ for which
$\phi$ extends to a family $\phi_n:{\cal C}_n\rightarrow {\cal P}_n$
of admissible covers
over ${\on {Spec}}\,{\mathbb C}[[t^{1/n}]]$. In particular, if the base for
the admissible cover $\phi:X^*\rightarrow Y^*$ is an open set $B^*$ of a smooth
projective curve $B$, modulo a finite base change, we can extend this
to a family of admissible covers $X^a\rightarrow Y^a$ over the whole curve $B$.

\subsection{Admissible covers of surfaces}
\label{admissible}
Consider a family $f:X\rightarrow B$ of stable curves of genus $g$, whose
general member is smooth and $d$-gonal. Let
$\psi:B\rightarrow {\overline{\mathfrak M}_g}$ be
the canonical map, and let $\overline B$ denote
the fiber product $B\times_{\overline{\mathfrak M}_g}
\overline{\cal H}_{d,g}$.
 
\begin{figure}
\begin{picture}(5,4)(-0.7,1.9)
\put(-2,4){$\overline{B}_0\subset
\overline B\stackrel{\eta}{\longrightarrow}
\overline{\cal H}_{d,g}$}
\multiput(-0.8,3.85)(1.4,0){2}{\vector(0,-1){0.9}}
\put(-1,2.5){$B\stackrel{\psi}{\longrightarrow}
\,\overline{\mathfrak M}_{g}$}
\put(-2.5,3.3){$X$}
\put(0.7,3.35){$\scriptstyle{pr_1}$}
\put(-3.5,4.8){$\overline{X}$}
\multiput(-3.2,4.7)(1.5,-0.8){2}{\vector(2,-3){0.7}}
\multiput(-2.1,3.4)(-1,1.5){2}{\vector(3,-2){1.1}}

\put(3,3.7){$B^*\subset B\stackrel{\eta}{\longrightarrow}
\overline{\cal H}_{d,g}$}
\multiput(4.2,5.05)(1.4,-1.5){2}{\vector(0,-1){0.9}}
\put(3.2,5.05){\vector(0,-1){0.9}}
\put(5.4,2.2){$\overline{\mathfrak M}_{g}$}
\put(5.7,3.1){$\scriptstyle{pr_1}$}
\put(3,5.2){$X^*\subset X$}
\put(4.3,3.6){\vector(1,-1){1}}
\put(4.8,3.1){$\scriptstyle{\psi}$}
\end{picture}
\caption{$\eta:\overline{B}\rightarrow\overline{\cal {H}}_{d,g}$
\hspace*{10mm}{\sc Figure 12.} Simply branched $C$\hspace*{-20mm}}
\label{map eta}
\end{figure}

\addtocounter{figure}{1}
If the general member of $X$ has infinitely many ${g}_d^1$'s,
the variety $\overline B$ will have dimension $\geq 2$.
We can resolve this by considering an intersection of
the appropriate number of hyperplane sections of 
$\overline B$, and picking a one-dimensional
component $\overline B_0$ dominating $B$. The curve
$\overline B_0$ might be
singular, but by normalizing it and pulling $X$
 over it, we get another family of stable curves (cf.~Fig.~\ref{map eta}):
 \[\overline X=X\times_B(\overline B_0)^{\on{norm}}
\rightarrow (\overline B_0)^{\on{norm}}.\]
 Since the two
families have the same basic invariants, we can replace
the original with the new one, and assume
the existence of a map $\eta:B\rightarrow
\overline{\cal H}_{d,g}$ compatible with $\psi:
B\rightarrow \overline{\mathfrak M}_g$. In other
words, $\eta$ associates to every fiber $C$ of $X$ a
specific $g^1_d$ on $C$ or, possibly, a $g^1_d$ 
on an admissible model $C^a$ of $C$.

\smallskip
 Let $B^*$ be the open subset of
$B$ over which {\it all} fibers are smooth and $d$-gonal. For simplicity, 
assume for now that all the fibers over $B^*$ can be 
represented as admissible covers of ${\proj}^1$ via the chosen
$g^1_d$'s, i.e. they are {\it simply branched} covers of ${\proj}^1$
over $m$ distinct points of ${\proj}^1$.
Denote by $X^*$ the restriction of $X$ over $B^*$ (cf.~Fig.12).

\smallskip
The map $\eta:B^*\rightarrow {\cal H}_{d,g}$ induces
a section 
\[\sigma:B^*\rightarrow {\on{Pic}}^d(X^*/B^*),\] where
${\on{Pic}}^d(X^*/B^*)$ is the {\it relative degree $d$ Picard variety} of
$X^*$ over $B^*$. ${\on{Pic}}^d(X^*/B^*)$ parametrizes
the line bundles on $X^*$ of relative degree $d$.
The image $\sigma(B^*)\subset {\on{Pic}}^d(X^*/B^*)$ is a class of line bundles
on $X^*$ whose fiberwise restrictions are the chosen $g^1_d$'s.
Let $\cal L$ be a representative of this class, and let
$Y^*$ be the ruled surface ${\proj}((f_*{\cal L})^{\widehat
{\phantom{n}}})$ over $B^*$.
The map $\phi:X^*\rightarrow Y^*$ induced by $\cal L$
defines an admissible cover over $B^*$, as shown in Fig.~\ref{construction
Y*}.

\setlength{\unitlength}{10mm}
\begin{figure}[h]
\begin{picture}(5,2.2)(-0.3,2.6)
\put(0,4){$X^*\stackrel{\phi}{\longrightarrow} Y^*={\proj}((f_*{\cal L})^
{\widehat{\phantom{n}}})$}
\put(0.2,3.85){\vector(1,-1){0.5}}
\put(1.4,3.85){\vector(-1,-1){0.5}}
\put(0.6,2.9){$B^*$}
\put(0.05,3.4){$f$}
\put(1.3,3.4){$h$}
\end{picture}
\caption{ Construction of $Y^*$}
\label{construction Y*}
\end{figure}

\bigskip
From Remark 3.1, $\phi$ extends to a family of admissible covers
${\phi}^a:{X}^a\rightarrow {Y}^a$ over the whole base $B$.
 Since ${X}^a$ and $X$ are 
isomorphic over $B^*$, they are birational to each other.
In other words, the fibers $C$ of $X$, over which $\cal L$ does not
extend to the base-point free linear series $g^1_d=\sigma_1(b)$, are 
modified by blow-ups and blow-downs so as to arrive at their admissible 
models in ${X}^a$. We have thus proved the following

\begin{lem} Let $f:X\rightarrow B$ be a family of stable curves, 
whose general member over an open subset $B^*\subset B$
is a smooth $d$-uple admissible cover of ${\proj}^1$.
Then, modulo a finite base change, there exists an admissible cover
of surfaces ${X}^a\rightarrow {Y}^a$
over $B$ such that ${X}^a$ is obtained from $X$ by a finite number of
birational transformations performed on the fibers over $B-B^*$. 
\label{quasicov}
\end{lem}

\subsection{Quasi-admissible covers}
\label{quasi-covers}
In case the general member of $X$ is {\it not} an admissible cover
of ${\proj}^1$, e.g. it is trigonal with a total point of ramification,
we have to modify the above construction. To start with, we cannot
expect to obtain an {\it admissible} cover $X^*\rightarrow Y^*$, even
modulo a finite base change. This leads us to consider a different kind
of covers, which we call {\it quasi-admissible}.

\medskip
\noindent{\bf Definition 3.2.} A {\it quasi-admissible cover}
$\widetilde{\phi}:
C\rightarrow P$ of a nodal curve $C$ over a semistable pointed rational  
curve $P$
is a regular map which behaves like an admissible cover over the singular 
locus of $P$, i.e. for any $q\in C$ lying over a node $p\in P$
the two branches through $q$ map with the same  ramification index
to the two branches through $p$.

\smallskip
\begin{figure}[h]
$$\vspace*{5mm}\psdraw{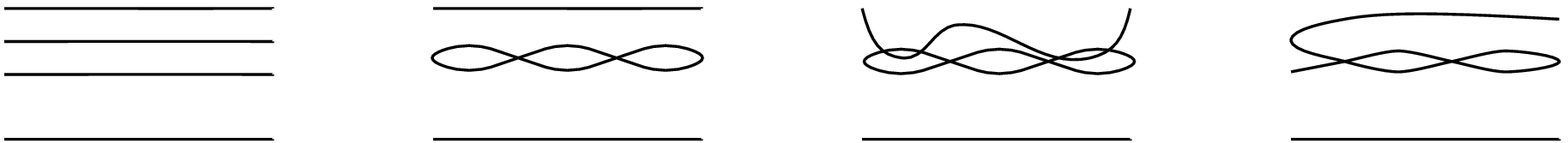}{4.5in}{0.5in}$$
\vspace*{-6mm}
\caption{Quasi-admissible covers over $\proj^1$}
\label{quasicovers}
\end{figure}

Quasi-admissible covers differ from admissible covers in allowing
more diverse behavior of $C$ over $P_{\on{sm}}$, e.g. having singularities,
higher ramification points and multiple simple ramification points.
Fig.~\ref{quasicovers} 
displays several degree 3 quasi-admissible covers over ${\proj}^1$:

However, any quasi-admissible cover can be obtained 
from an admissible cover $\phi^a\!:\!C^a\!\rightarrow\! P^a$
by simultaneous contractions of components in $P^a$
and their (rational) inverses on $C^a$.

\medskip
\noindent{\bf Definition 3.3.} A {\it minimal} quasi-admissible cover 
$\widetilde{\phi}:C\rightarrow P$ is minimal with respect to the number of 
components of $P$. In other words, one cannot apply more simultaneous
contractions on $C\rightarrow P$ and end with another quasi-admissible cover.

\begin{figure}[h]
$$\psdraw{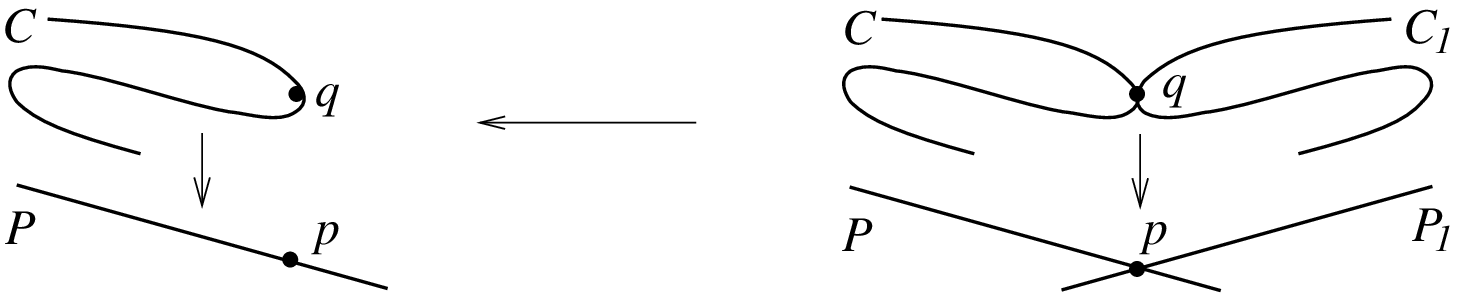}{2.4in}{0.8in}$$
\caption{Quasi virsus admissible covers}
\label{quasi/adm}
\end{figure}

\smallskip
\noindent{\bf Example 3.2.}
A smooth trigonal curve $C$ with a total point of ramification $q$ is a 
minimal quasi-admissible cover of $P={\proj}^1$. Blowing up $q$ on $C$
and $p=\widetilde{\phi}(q)\in P$, gives an admissible cover $C^a=C\cup
C_1\rightarrow P\cup P_1$, where $C_1\cong {\proj}^1$  maps
three-to-one onto $P_1\cong {\proj}^1$ with a total
point of ramification $q=C_1\cap C$ (cf.~Fig.~\ref{quasi/adm}).

\bigskip

\medskip
The motivation for using {\it minimal} quasi-admissible covers, 
instead of just admissible or quasi-admissible covers, is that the former are
the closest covers to the original families $X\rightarrow B$ of stable
curves, and calculations
on them will yield the best possible estimate for the ratio 
$\delta_X/\lambda_X $ (cf.~Fig.~\ref{stages}).

\subsubsection{Quasi-admissible covers for families with
higher ramification sections}
Now let us consider the remaining case of 
a family $X\rightarrow B$, whose general member over
$B^*$ is smooth and $d$-gonal, but {\it not} an admissible cover of 
${\proj}^1$. After a possible base change, we still have the map
(cf.~Fig.~{12})
\[\eta:B\longrightarrow {\overline {\cal H}}_{d,g}.\]
It associates to every fiber $C$ a $g^1_d$ on its admissible 
model $C^a$. Let $C^a\rightarrow P^a$ be the 
corresponding admissible cover. Since $C$ itself is $d$-gonal, and 
by assumption it does not possess
a $g^1_e$ with $e<d$, $C$ must be a $d$-uple cover of some component of
$P^a$. In particular, the $g^1_d$ on 
$C^a$ restricts to a $g^1_d$ on $C$. Thus, in effect, $\eta$
gives again a section $\sigma:B^*\rightarrow {\on{Pic}}^d(X^*/B^*)$.
As before, we obtain a degree $d$ finite map $\phi:X^*\rightarrow Y^*$ to 
the ruled surface ${\mathbf P}((f_*{\cal L})^{\widehat{\phantom{n}}})$ 
over $B^*$. Note that
this is a family of {\it minimal quasi-admissible} covers.
 
\medskip
We extend $\phi$ over the curve $B$ as follows.
For simplicity, assume that $d=3$. Let $R$ be the
ramification divisor of $\phi$ in $X^*$. By hypothesis,
there is a component $R_0$ of $R$ which passes through
total  ramification points and dominates $B^*$. 
Letting $\overline{R}_0$
be the closure of $R_0$ in $X$, we can normalize it and pull the family $X$
over it. So we may assume that $\overline{R}_0$ is a section of $X\rightarrow
B$. If there are some other components $R_1,R_2,...,R_l$ of the
ramification divisor $R$ passing through higher ramification points,
we repeat the same procedure for them, until we ``straighten out'' all
$\overline{R}_i$'s into sections of $X\rightarrow B$.
Let $E_i=\phi(R_i)$ be the corresponding sections of $Y^*$ over $B^*$.
We can shrink $B^*$ in order to exclude any fibers with isolated
higher ramification points. 

\smallskip
Consider a fiber $C$ in $X^*$. Let $\{r_i=C\cap R_i\}$ be its
total ramification points, and let $\{p_i=\phi(r_i)\}$ be their images
on $P=\phi(C)$ in $Y^*$. It is clear that blowing-up all $r_i$'s 
and $p_i$'s will give an admissible triple cover $C^a=
\on {Bl}_{\{r_i\}}(C)\rightarrow P^a=\on {Bl}_{\{p_i\}}(P)$. 
The $g^1_d$, giving
this cover, is the original one assigned by $\eta:B^*\rightarrow
{\overline{\cal H}}_{d,g}$. We globalize this construction by
blowing-up the sections ${R}_i$ on $X^*$ and $E_i$ on $Y^*$.
Similarly as above, we obtain a triple admissible cover of surfaces
$\phi^*:\on{Bl}_{\cup R_i}(X^*)\rightarrow \on{Bl}_{\cup E_i}
(Y^*)$ over $B^*$. The properness of $pr_1:
\overline{\cal H}_{d,g}\rightarrow \overline{\mathfrak M}_g$ allows us
to extend this to an admissible cover ${\phi}^a:
\overline{\on{Bl}_{\cup R_i}(X^*)}
\rightarrow\overline{\on{Bl}_{\cup E_i}(Y^*)}$ over $B$ (cf.~
Fig.~\ref{blowing up}).

\begin{figure}[h]\hspace*{-30mm}
\begin{picture}(3,5)(4.3,-0.5)
\put(2,3.7){$\overline{{\cal R}_i}\subset 
           \overline{\on{Bl}_{\cup R_i}(X^*)}\stackrel{{\phi}^a}
           {\longrightarrow} \overline{\on{Bl}_{\cup E_i}(Y^*)}\supset 
             \overline{{\cal E}_i}$}
\put(2,2.2){${{\cal R}_i}\subset {\on{Bl}_{\cup R_i}(X^*)}
              \stackrel{\phi^*}{\longrightarrow}
              {\on{Bl}_{\cup E_i}(Y^*)}\supset{{\cal E}_i}$}
\put(2,0.7){$R_i\hspace{0.5mm}\subset \hspace{6.6mm}
X^*\hspace{6mm}\stackrel{\phi}{\longrightarrow}\hspace{7.3mm}Y^*
\hspace*{5.6mm}\supset E_i$}
\multiput(2.2,3.5)(6.1,0){2}{\vector(0,-1){0.9}}
\multiput(3.9,3.5)(2.8,0){2}{\vector(0,-1){0.9}}
\multiput(2.2,2)(6.1,0){2}{\vector(0,-1){0.9}}
\multiput(3.9,2)(2.8,0){2}{\vector(0,-1){0.9}}
\end{picture}
\vspace*{-10mm}
\caption{Blowing up ${R}_i$ and ${E}_i$ \hspace*{15mm}
{\sc Figure 17.} Over $B^*$\hspace*{-10mm}}
\label{blowing up}
\hspace*{80mm}\begin{picture}(3,0)(6.7,-2.3)
\put(7.8,2.2){${X}^q\stackrel{\phi^q}{\longrightarrow} {Y}^q$}
\multiput(8,2)(1.4,0){2}{\vector(0,-1){0.9}}
\multiput(8.1,1.5)(1.4,0){2}{$\wr$}
\put(7.8,0.7){$X^*\longrightarrow Y^*$}
\put(8,0.5){\vector(1,-1){0.5}}
\put(9.4,0.5){\vector(-1,-1){0.5}}
\put(8.5,-0.4){$B^*$}
\end{picture}
\label{over B*}
\end{figure}

\addtocounter{figure}{1}
\smallskip
Denote by ${\cal R}_i$ 
the component of $\on {Bl}_{\cup R_i}(X^*)$, obtained by blowing up 
${R}_i\subset X^*$, 
and let $\overline{{\cal R}_i}$ be its closure in 
$\overline{\on{Bl}_{\cup R_i}(X^*)}$.
Define similarly ${\cal E}_i\subset \on{Bl}_{\cup E_i}(Y^*)$ and 
$\overline{\cal E}_i
\subset \overline{\on{Bl}_{\cup E_i}(Y^*)}$. The admissible cover 
$\phi^a$ maps
$\overline{\cal R}_i$ to $\overline{\cal E}_i$, so that after removing 
all the $\overline{\cal R}_i$'s and $\overline{\cal E}_i$'s we still have
a triple cover 
\[{\phi}^q:{X}^q=\overline{\on{Bl}_{\cup R _i}
(X^*)}-\cup\overline{\cal R_i}\longrightarrow 
{Y}^q=\overline{\on{Bl}_{\{{ E_i}\}}
(Y^*)}-\cup\overline{\cal E_i}.\]

Note that ${X}^q\cong X$ and ${Y}^q\cong Y$ over the open set
$B^*$, and that ${Y}^q$ is a birationally ruled surface over $B$
(cf.~Fig.~17).
Finally, note that from the quasi-admissible cover
${\phi}^q:{X}^q\rightarrow {Y}^q$ we obtain a family 
$\widetilde{\phi}:\widetilde{X}\rightarrow
\widetilde{Y}$ of {\it minimal} 
quasi-admissible covers: simply contract the unnecessary
rational components in the fibers of ${X}^q$ and
${Y}^q$, and observe that the triple map ${\phi}^q$ restricts to the
corresponding triple map $\widetilde{\phi}$. 

\smallskip This completes the construction of minimal quasi-admissible
covers for any family $X\rightarrow B$ with general smooth trigonal member.
The cases $d>3$ are only notationally more difficult. One has to keep
track of the possibly different higher multiplicities in $C$ and multiple 
double points in $C$ over the same $p\in P$. The construction 
of an admissible cover ${X}^a\rightarrow {Y}^a$ 
goes through with minimal modifications. We combine the results of this
section in the following

\begin{prop} 
Let $f:X\rightarrow B$ be a family of stable curves, 
whose general member over an open subset $B^*\subset B$
is  smooth and $d$-gonal.
Then, modulo a finite base change, there exists a minimal quasi-admissible 
cover of surfaces $\widetilde{X}\rightarrow \widetilde{Y}$
over $B$ such that $\widetilde{X}$ is obtained from $X$ by a finite number of
birational transformations performed on the fibers over
$B-B^*$.
\label{propquasi}
\end{prop}

\medskip
\section*{\hspace*{1.9mm}4. The Boundary $\Delta{\mathfrak{T}}_g$
of the Trigonal Locus $\overline{\mathfrak{T}}_g$}

\setcounter{section}{4} 
\setcounter{subsection}{0} 
\setcounter{subsubsection}{0} 
\setcounter{lem}{0}
\setcounter{thm}{0}
\setcounter{prop}{0}
\setcounter{defn}{0} 
\setcounter{cor}{0} 
\setcounter{conj}{0} 
\setcounter{claim}{0} 
\setcounter{remark}{0}
\setcounter{equation}{0}

\label{boundarycomponents}
\subsection{Description and notation for the boundary of
$\overline{\mathfrak{T}}_g$}
\label{description} In this section we shall see that
there are {\it seven types} of boundary divisors of
$\overline{\mathfrak{T}}_g$, each denoted by $\Delta{\mathfrak{T}}_{k,i}$
for $k=0,1,...,6$. The second index $i$ is determined in the following
way. Let $C=C_1\cup C_2$ be the 
general member of $\Delta{\mathfrak{T}}_{k,i}$, where $C_1$ and $C_2$ are smooth curves.
If $C_1$ and $C_2$ are both trigonal or both hyperelliptic, then we set 
$i$ to be the smaller of the two genera $p(C_1)$ or $p(C_2)$. If, say, $C_1$ is
a trigonal, but $C_2$ is hyperelliptic, then we set $i$ to be
genus of the trigonal component $C_1$.
The only exception to this rule occurs 
when $C$ is irreducible (and hence of genus $g$ with exactly one
node). We denote this boundary component by $\Delta{\mathfrak{T}}_0$.

\smallskip
When we view a general member $C$ roughly as a triple
cover of ${\proj}^1$'s in the Hurwitz scheme (consider the pull-back
$pr_1[C]\in\overline{\cal{H}}_{3,g}$), then it may or may not be
ramified. If there is no ramification, then $C$ lies in one of the first
four types of trigonal boundary divisors $\Delta{\mathfrak{T}}_{k,i}$, 
$k=0,1,2,3$. Ramification index 1 characterizes the general members of
$\Delta{\mathfrak{T}}_{4,i}$ and $\Delta{\mathfrak{T}}_{5,i}$, and in case of 
$\Delta{\mathfrak{T}}_{6,i}$ the ramification index is 2 
(cf.~Fig.~\ref{Delta-k,i}).

\smallskip
There is an alternative description of the boundary components
$\Delta{\mathfrak{T}}_{k,i}$'s of $\overline{\mathfrak{T}}_g$. 
\vspace*{-1mm}If one such $\Delta{\mathfrak{T}}_{k,i}$ lies in the
restriction $\Delta_0\big|_{\displaystyle{\overline{\mathfrak{T}}_g}}$ 
of the divisor $\Delta_0$ in $\overline{\mathfrak{M}}_g$, 
\vspace*{-1mm}then $\Delta{\mathfrak{T}}_{k,i}$ is one of $\Delta{\mathfrak{T}}_0,\,\,
\Delta{\mathfrak{T}}_{1,i},\,\,\Delta{\mathfrak{T}}_{2,i},$
or $\Delta{\mathfrak{T}}_{4,i}$. The partial normalization of their general 
members $C$ is still connected, i.e. $C$ is either irreducible, 
or the join of two smooth curves meeting in at least two
points. Correspondingly, for the general member $C$ of the remaining
three types of boundary components,
$\Delta{\mathfrak{T}}_{3,i},\,\,\Delta{\mathfrak{T}}_{5,i}$ and $\Delta{\mathfrak{T}}_{6,i}$, the irreducible
components of $C$ intersect transversely in exactly one point, so that the
normalization of $C$ is  disconnected.

\bigskip
\begin{figure}[h]
$$\psdraw{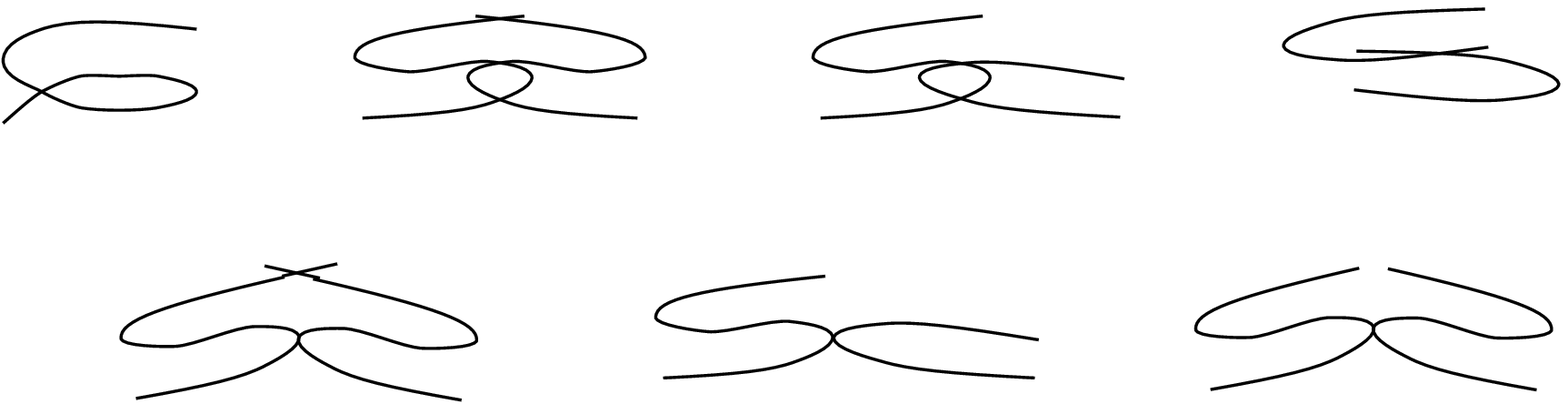}{4.5in}{1in}$$
\begin{picture}(6,1)(2.7,-1.4)
\put(-0.7,2.45){$\Delta{\mathfrak{T}}_{0} \hspace{19mm}\Delta{\mathfrak{T}}_{1,i}
\hspace{25mm}\Delta{\mathfrak{T}}_{2,i}
\hspace{23mm}\Delta{\mathfrak{T}}_{3,i}$}
\put(2.8,1.2){$\scriptstyle{i=1,2,...,}\left[\frac{g-2}{2}\right]\hspace{17mm}
\scriptstyle{i=1,2,..., g-2}\hspace{17mm}\scriptstyle{i=1,2,
...,} \left[\frac{g}{2}\right]$}
\put(-0.2,0.5){$\Delta{\mathfrak{T}}_{4,i}\hspace{31mm}\Delta{\mathfrak{T}}_{5,i}
\hspace{31mm}\Delta{\mathfrak{T}}_{6,i}$}
\put(01.1,-0.6){$\scriptstyle{i=1,2,...,\left[\frac{g-1}{2}\right]\hspace{24mm}
i=1,2,...,g-1 \hspace{22mm}i=1,2,...,\left[\frac{g}{2}\right]}$}
\end{picture}
\vspace*{-10mm}
\caption{Boundary Components $\Delta\mathfrak{T}_{k,i}$
of $\overline{\mathfrak{T}}_g$}
\label{Delta-k,i}
\end{figure}

\begin{prop} The boundary divisors of $\overline{\mathfrak{T}}_g$ can
be grouped in seven types: $\Delta{\mathfrak{T}}_{0}$ and $\Delta{\mathfrak{T}}_{k,i}$
for $k=1,...,6$. Their general members and range of index $i$ are shown in
Fig.~\ref{Delta-k,i}. The boundary of $\overline{\mathfrak{T}}_g$ consists of
$\Delta{\mathfrak{T}}_{0}$, $\Delta{\mathfrak{T}}_{k,i}$, and the codimension 2
component $\overline{\mathfrak{I}}_g$ of hyperelliptic curves.
\label{boundary}
\end{prop}

\medskip
Consider the projection map
$pr_1:\overline{\cal{H}}_{3,g}\rightarrow \overline{\mathfrak{M}}_g$, whose
image is the trigonal locus $\overline{\mathfrak{T}}_g$. Thus, the
 inverse image of 
each boundary divisor $\Delta{\mathfrak{T}}_{k,i}$ will be a boundary
divisor $\Delta{\cal{H}}_{k,i}$ in $\overline{\cal{H}}_{3,g}$.
The converse, however, is not always true, i.e. certain boundary divisors
of $\overline{\cal{H}}_{3,g}$ contract under $pr_1$ to smaller subschemes of
$\overline{\mathfrak{T}}_g$, e.g. the hyperelliptic locus $\overline{\mathfrak{I}}_g$. 
With the description of the Hurwitz scheme
$\overline{\cal{H}}_{3,g}$, given in Section~3,
it is easier to determine first $\overline{\cal{H}}_{3,g}$'s boundary
divisors. Thus, we postpone the proof of Proposition~\ref{boundary} until
the end of the next subsection.

\subsubsection{The Boundary of $\overline{\cal{H}}_{3,g}$.}
\label{admissibleaboundary}
\begin{prop} The boundary divisors of $\overline{\cal{H}}_{3,g}$ can
be grouped in six types: 
$\Delta{\cal{H}}_{k,i}$
for $k=1,...,6$. Their general members and range of index $i$ are shown 
in Fig.~\ref{admissible-k,i}.
\begin{figure}[h]
\bigskip
$$\psdraw{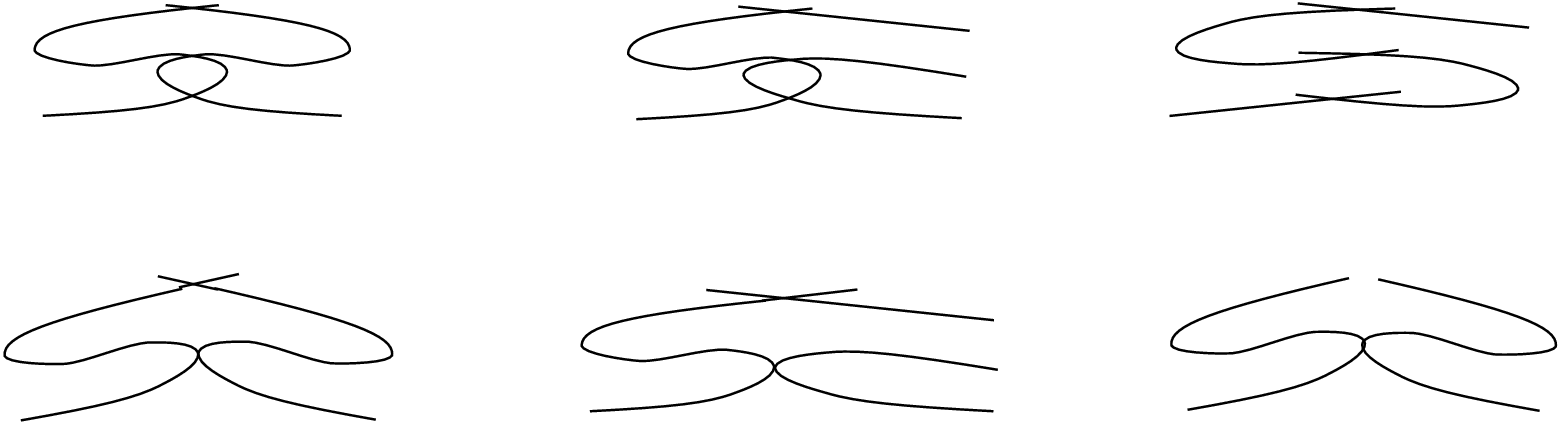}{4.5in}{1in}$$
\begin{picture}(5,1)(3,-2)
\put(-0.9,1.85){$\Delta{\cal{H}}_{1,i}
\hspace{33mm}\Delta{\cal{H}}_{2,i}
\hspace{32mm}\Delta{\cal{H}}_{3,i}$}
\put(0.3,0.7){$\scriptstyle{i=1,2,...,}\left[\frac{g-2}{2}\right]\hspace{29mm}
\scriptstyle{i=1,..., g-1}\hspace{24mm}\scriptstyle{i=0,1,
...,} \left[\frac{g}{2}\right]$}
\put(-0.9,0.15){$\Delta{\cal{H}}_{4,i}\hspace{33mm}\Delta{\cal{H}}_{5,i}
\hspace{32mm}\Delta{\cal{H}}_{6,i}$}
\put(0.3,-1.2){$\scriptstyle{i=1,2,...,\left[\frac{g-1}{2}\right]\hspace{29mm}
i=1,2,...,g-1 \hspace{23mm}i=1,2,...,\left[\frac{g}{2}\right]}$}
\end{picture}
\vspace*{-5mm}
\caption{Boundary Components of $\overline{\cal{H}}_{3,g}$}
\label{admissible-k,i}
\end{figure}
\label{boundary2}
\end{prop}

\begin{proof}
A general member $A$ of the boundary $\Delta{\cal{H}}$ is a
triple admissible cover of a chain of {\it two} ${\proj}^1$.
(From the dimension calculations that follow it will become clear that
an admissible cover of a chain of three or more $\proj^1$'s will
generate a subscheme in $\overline{\cal{H}}_{3,g}$ of codimension $\geq 2$.)
Note that {\it
three} connected components of $A$ over one ${\proj}^1$ means that they are
all smooth ${\proj}^1$'s themselves, and hence they can all be
contracted simultaneously, leaving us with a smooth trigonal curve, or
with a hyperelliptic curve with an attached ${{\proj}}^1$, neither of which
cases by dimension count corresponds to a {\it general} member of 
a boundary component $\Delta{\cal{H}}_{k,i}$. Considering all
combinations of one or two connected components of $A$ over each
${\proj}^1$, we generate a list of the possible general members of
the boundary divisors $\Delta{\cal{H}}_{k,i}$.

To see which of these are indeed of codimension 1 in
$\overline{\cal{H}}_{3,g}$, we do the following calculation. First we
note that, for
a fixed set of $2i+4$ ramification points in ${\proj}^1$, there are finitely
many covers of degree $3$ and genus $i$, that is,
\[\on{dim}\overline{\mathfrak{T}}_i=2i+4-3=2i+1.\]
Substracting $3$ takes
into account the projectively equivalent triples of points on ${\proj}^1$.
In particular, $\on{dim}\overline{\mathfrak{T}}_g=2g+1$.
A similar agrument (with $2i+2$ ramification points) shows that for
the hyperelliptic locus:
\[\on{dim}\overline{\mathfrak{I}}_i=2i+2-3=2i-1.\]
These computations are valid for $i>0$, whereas $0=
\on{dim}\overline{\mathfrak{T}}_i=\on{dim}\overline{\mathfrak{I}}_i$.

\smallskip
Thus, to compute the dimensions of the six types of subschemes of
$\overline{\cal{H}}_{3,g}$, one adds the corresponding dimensions of
$\overline{\mathfrak{T}}_i$ and $\overline{\mathfrak{I}}_j$, making 
the necessary adjustments for the choice of intersection points on
the components of each curve $A$. For example, when $i>0$
the dimension of the subscheme with general member $A$,
shown in  Fig.~\ref{admissible-k,i}, is 
\[\on{dim}\overline{\mathfrak{T}}_i+\on{dim}\overline{\mathfrak{T}}_{g-i-2}+1+1=2g.\]
The final 1's account for the choice of triples of points in the
$g^1_3$'s on each component. We conclude that for
$i=1,2,...,[(g-2)/2]$ the join at three points of two trigonal curves,
one of genus $i$ and the other of genus $g-i-2$, is the general member of
a boundary component of $\overline{\cal{H}}_{3,g}$. We denote it by
$\Delta{\cal{H}}_{1,i}$. The range of $i$ stops at 
$[(g-2)/2]$ for symmetry considerations. When $i=0$, the corresponding
subscheme has a smaller dimension of $2g-2$ and hence no boundary divisor
is generated by such curves.

\smallskip
As another example, consider the fifth sketch in Fig.~\ref{admissible-k,i}.
It corresponds to the join at one point 
of a trigonal curve $C_1$ of genus $i$, a hyperelliptic curve $C_2$
of genus $g-i$, and an attached ${\proj}^1$ to $C_2$ to make the whole curve a 
triple cover. Note that $C_1$ and $C_2$ intersect transversally at a point
$q$, but when presented as covers of ${\proj}^1$ they both have 
ramifications at $q$ of index 1. On all such curves $C_1$ and $C_2$
the total number of ramification points over ${\proj}^1$ is finite, and
hence their choice does not affect the dimension of our subscheme. Thus,
\[\on{dim}\overline{\mathfrak{T}}_i+\on{dim}\overline{\mathfrak{I}}_{g-i}=2i+1+2(g-i)=
2g.\]
Therefore, this subscheme is in fact a divisor in
$\overline{\cal{H}}_{3,g}$, which we denote by $\Delta
{\cal{H}}_{5,i}$. The cases of $i=0$ or $i=g$ lead to
contractions of unstable rational components ($C_1$ or $C_2$), and do not
yield the necessary dimension of $2g$. Hence, $i=1,2,...,g-1$.

\smallskip
In the case of $\Delta{\cal{H}}_{6,i}$, the two components $C_1$ and
$C_2$ meet transversally in one point $q$, but both have ramification of
index $2$ at $q$ as triple covers of ${\proj}^1$. Smooth trigonal curves of
genus $i$ with such high ramification form a codimension 1 subscheme of the 
trigonal locus ${\mathfrak{T}}_i$, hence the dimension of
$\Delta{\cal{H}}_{6,i}$ is
\[\on{\dim}\overline{\mathfrak{T}}_i-1+\on{dim}\overline{\mathfrak{T}}_{g-i}-1=
(2i+1)-1+(2(g-i)+1)-1=2g.\]
Thus, $\Delta{\cal{H}}_{6,i}$ is a boundary divisor in 
$\overline{\cal{H}}_{3,g}$ for $i=1,2,...,[g/2]$. The case of $i=0$ yields
dimension $2g-1$, and hence we disregard it.

\smallskip
The remaining cases are treated similarly. We conclude 
that $\overline{\cal{H}}_{3,g}$ has
six types of boundary divisors, $\Delta{\cal{H}}_{k,i}$,
 whose general members and range of indices
are indicated in Fig.~\ref{admissible-k,i}. \end{proof}

\subsubsection{Boundary of $\overline{\mathfrak{T}}_g$. 
Proof of Proposition~\ref{boundary}} 
\label{trigonalboundary}

Having described the boundary of $\overline{\cal{H}}_{3,g}$, it remains to
check which of the divisors $\Delta\cal{H}_{k,i}$ preserve their
dimension under the map $pr_1$ and hence map into divisors of $\overline
{\mathfrak{T}}_g$. The only ``surprises'' can be expected where $pr_1$
contracts unstable ${\proj}^1$, such as in $\Delta{\cal{H}}_{2,i}$,
$\Delta{\cal{H}}_{3,i}$, and $\Delta{\cal{H}}_{5,i}$. In fact, only
$\Delta{\cal{H}}_{2,g-1}$ and $\Delta{\cal{H}}_{3,0}$ diverge from
the common pattern; in all other cases, we set $\Delta{\mathfrak{T}}_{k,i}:=
pr_1\left(\Delta{\cal{H}}_{k,i}\right)$ 
to be the corresponding boundary divisor in $\overline{\mathfrak{T}}_g$.

\smallskip
The map $pr_1$ contracts the three rational components of
the general member of $\Delta{\cal{H}}_{3,0}$, leaving only
a smooth hyperelliptic curve of genus $g$. Thus, the image
$pr_1\left(\Delta{\cal{H}}_{3,0}\right)$ is the hyperelliptic locus
$\overline{\mathfrak{I}}_g$, which is of dimension $2g-1$. Hence 
$\Delta{\cal{H}}_{3,0}$ does not yield a divisor in
$\overline{\mathfrak{T}}_g$, but a boundary component of codimension 2.

\smallskip
 Finally we consider $\Delta{\cal{H}}_{2,g-1}$. After we contract
its two rational components, we arrive at an {\it irreducible nodal}
trigonal curve with exactly one node. The dimension of the subscheme
of such curves is  
\[\on{dim}\overline{\mathfrak{T}}_{g-1}+1=2(g-1)+1+1=2g,\]
where the final $1$ indicates the choice of a triple of points on
a smooth trigonal curve (belonging to the $g^1_3$), two of which
will be identified as a node. Correspondingly, we obtain another divisor
in $\overline{\mathfrak{T}}_g$, which we denote by $\Delta{\mathfrak{T}}_0.$
\qed

\subsection{Multiplicities of the boundary divisors $\Delta{\mathfrak{T}}_{k,i}$
in the restriction $\delta|_{\overline{\mathfrak{T}}_g}$} 
\label{multiplicities}
By abuse of
notation, we will denote by  $\delta_0$ and $\delta_{k,i}$
the classes in $\on{Pic}_{\mathbb{Q}}\overline{\mathfrak{T}}_g$
of $\Delta{\mathfrak{T}}_0$ and $\Delta{\mathfrak{T}}_{k,i}$, respectively.

\begin{prop} The divisor class $\delta\in\on{Pic}_{\mathbb{Q}}
\overline{\mathfrak{M}}_g$
restricts to $\overline{\mathfrak{T}}_g$ as the following linear combination
of the boundary classes in $\overline{\mathfrak{T}}_g$:
\begin{equation}
\delta|_{\displaystyle{\overline{\mathfrak{T}}_g}}=\delta_0+
\sum_{i=1}^{\scriptscriptstyle{[(g-2)/2]}}3\delta_{1,i}
+\sum_{i=1}^{\scriptscriptstyle{g-2}}2\delta_{2,i}
+\sum_{i=1}^{\scriptscriptstyle{[g/2]}}\delta_{3,i}
+\sum_{i=1}^{\scriptscriptstyle{[(g-1)/2]}}3\delta_{4,i}+
\sum_{i=1}^{\scriptscriptstyle{g-1}}\delta_{5,i}
+\sum_{i=1}^{\scriptscriptstyle{[g/2]}}\delta_{6,i}.
\label{divisorrel}
\end{equation}
\end{prop}

\noindent{\it Proof.} Let us rewrite equation (\ref{divisorrel}) in the form
\[\delta|_{\displaystyle{\overline{\mathfrak{T}}_g}}=(\on{mult}_{\delta}
\delta_0)\delta_0+\sum_{k,i}(\on{mult}_{\delta}\delta_{k,i})\delta_{k,i},\]
and call $\on{mult}_{\delta}\delta_{k,i}$ the {\it multiplicity} of
$\delta_{k,i}$ in $\delta|_{\overline{\mathfrak{T}}_g}$.
This linear relation simply counts the contribution
of each singular curve of a specific boundary type in $\Delta\mathfrak{T}_g$
to the degree of
$\delta$. Recall that for any trigonal family $f:X\rightarrow B$:
\[\on{deg}\delta|_B=\sum_{q\in X}m_q.\]
Here $m_q$ denotes the local analytic multiplicity of the total space of
$X$ nearby $q$ measured by the equation $xy=t^{m_q}$, where
$x$ and $y$ are local parameters on the singular fiber $X_b$, and
$t$ is a local parameter on $B$ near $b=f(q)$.

\smallskip
For each boundary class $\Delta{\mathfrak{T}}_{k,i}$
of $\overline{\mathfrak{T}}_g$, we consider its general member
$C\!=\!C_1\cup C_2$,
and a base curve $B$ in $\overline{\mathfrak{T}}_g$ which intersects
transversally $\Delta{\mathfrak{T}}_{k,i}$ in $[C]$. In the
corresponding one-parameter trigonal family $f:X\rightarrow B$, 
we must find the sum of the multiplicities $m_q$ of the singularities
of $C$. Thus, \[\on{mult}_{\delta}\delta_{k,i}=\sum_{\,\,q\in
C_{\on{sing}}}\!\!m_q.\]

For most of the divisors classes, this sum is actually quite straight forward.
For example, the general member $[C]\in \Delta{\mathfrak{T}}_{3,i}$
is the join of two smooth hyperelliptic curves $C_1$ and $C_2$, which
intersect transversally in one point $q$. The family $X$,
constructed as above, will be given locally analytically nearby $q$ by
$xy=t$, and hence $\on{mult}_{\delta}\delta_{k,i}=m_q=1$.
A similar situation occurs in the cases of $\Delta\mathfrak{T}_0,
\Delta\mathfrak{T}_{5,i}$ and $\Delta\mathfrak{T}_{6,i}$: there is one point
of transversal intersection (or one node) forcing 
\[\on{mult}_{\delta}\delta_0=
\on{mult}_{\delta}\delta_{k,i}=1\,\,\on{for}\,\, k=3,5,6.\] 

\smallskip
In the cases of $\Delta\mathfrak{T}_{2,i}$ and $\Delta\mathfrak{T}_{1,i}$ there
are correspondingly two or three points of transversal intersection,
forcing \[\on{mult}_{\delta}\delta_{2,i}=2\,\,\on{and}\,\,
\on{mult}_{\delta}\delta_{1,i}=3.\] This can be also interpreted
by the fact that $\Delta\mathfrak{T}_{2,i}$ and
$\Delta\mathfrak{T}_{1,i}$ lie entirely in the divisor
 $\Delta_0$ in $\overline{\mathfrak{M}}_g$ with,
 $\Delta_0$ being {\it double} along $\Delta\mathfrak{T}_{2,i}$ and
{\it triple} along $\Delta\mathfrak{T}_{1,i}$.

\begin{figure}
$$\psdraw{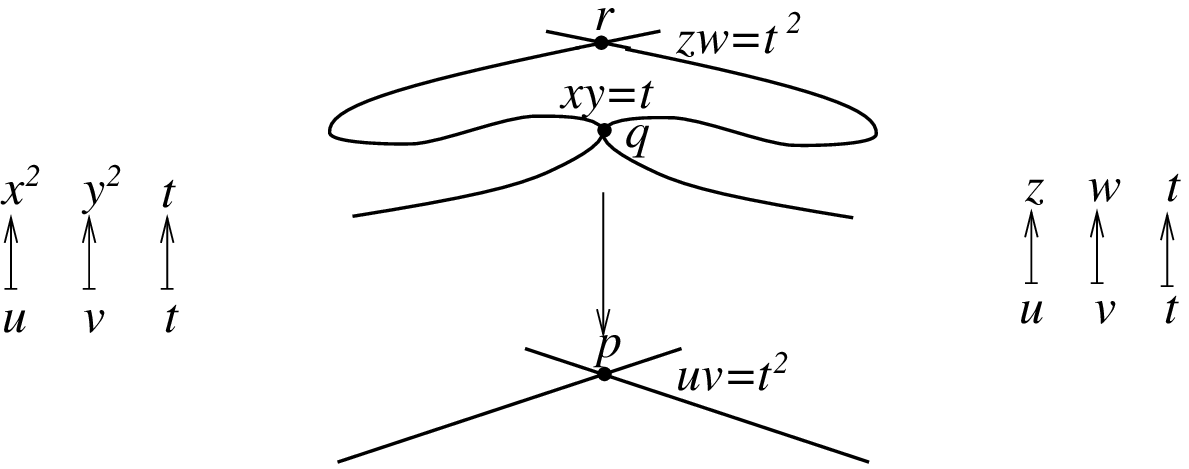}{2.8in}{1.3in}$$
\caption{The multiplicity mult$_{\delta_0}\delta_{4,i}$}
\label{mult4}
\end{figure}

\medskip
A slightly more complex situation occurs in the case of
$\Delta\mathfrak{T}_{4,i}$. The general member $C$ consists of two curves 
$C_1$ and $C_2$, meeting transversally in two points $q$ and $r$
(see Fig.~\ref{mult4}). But,
as in an admissible triple cover of two ${\proj}^1$'s, the points $q$ and
$r$ behave differently: at one of them, say $r$, the triple cover
is {\it not} ramified, while at $q$ there is ramification of index $1$.
In the local analytic rings of $p,q$ and $r$
the generators of $\widehat{\cal{O}}_{Y,p}$ map into the squares of 
the generators of $\widehat{\cal{O}}_{X,q}$: $u\mapsto x^2, v\mapsto y^2$, and
of course, $t\mapsto t$, so that the local equation of $Y$ near $p$ is
$uv=t^2$, and that of $X$ near $q$ is $xy=t$. But since the triple cover
is a local isomorphism of $\widehat{\cal{O}}_{Y,p}$ into $\widehat{\cal{O}}
_{X,r}$, the total space of $X$ near $r$ is given locally analytically
by $zw=t^2$ ($u\mapsto z, v\mapsto w, t\mapsto t$). 
Therefore, $m_q=1$, but $m_r=2$, and
\[\on{mult}_{\delta_0}\delta_{4,i}=m_q+m_r=3.\,\,\,\qed\]

\subsection{The hyperelliptic locus $\overline{\mathfrak{I}}_g$ inside
$\overline{\mathfrak{T}}_g$}
\label{hyperelliptic locus} 
Although the relations proved in this paper
will be valid on the Picard group $\on{Pic}_{\mathbb
Q}\overline{\mathfrak{T}}_g$, it will be interesting to check what happens
with the hyperelliptic curves inside the trigonal locus
$\overline{\mathfrak{T}}_g$. 

We noted that $\overline{\mathfrak{I}}_g$ is the
only boundary component of $\overline{\mathfrak{T}}_g$ of
codimension 2. It is obtained as the image $pr_1(\Delta{\cal H}_{3,0})$. 
By blowing up a point on a smooth hyperelliptic
curve $C$, we add a $\proj^1$--component to $C$ to make it a triple cover
$C^{\prime}$ (cf.~Fig.~\ref{smoothhyper}). 
It terms of the quasi-admissible covers, such $C^{\prime}$
behaves exactly as an irreducible singular trigonal curve in
$\Delta{\mathfrak{T}}_0$. However, $C$ 
does not contribute to the invariant $\delta|_B$, as defined in
Section~\ref{definition}. In fact, in a certain sense,
it even decreases $\delta|_B$. 

To simplify the exposition, we shall postpone the discussion
of families with hyperelliptic fibers until Section~11, where
we will explain the behavior of trigonal families with finitely many
hyperelliptic fibers in terms of the exceptional divisor $\Delta{\cal H}_{3,0}$
of the projection $pr_1$.

A similar phenomenon occurs with the boundary component
$\Delta\mathfrak{T}_{1,0}=pr_1(\Delta\cal{H}_{1,0})$, but it does not
make sense to exclude its members from our discussion, since they behave
exactly as members of the boundary divisor $\Delta\mathfrak{T}_{1,i}$ for
$i\geq 1$. 

\subsection{The invariants $\mu(C)$}
\label{The invariants} In the transition from the original
family $X\rightarrow B$ to the minimal quasi-admissible family
$\widetilde{X}\rightarrow \widetilde{Y}$ over $\widetilde{B}$, certain
changes occur in the calculation of the basic invariants. To start with,
it is easy to redefine
$\lambda_{\widetilde{X}},\kappa_{\widetilde{X}}$ and
$\delta_{\widetilde{X}}$ for $\widetilde{X}\rightarrow \widetilde{B}$:
simply use the corresponding definitions from Section~\ref{definition}.
Since we are interested in the slope of the family, which
is invariant under base change, we may assume that $\widetilde{B}:=B$
and that $X$ is the pull-back over the new base $\widetilde{B}$.
Now the difference between $X$ and $\widetilde{X}$ 
is reduced to several ``quasi-admissible'' blow-ups on $X$.

\smallskip
Blowing up smooth or rational double points
on a surface does not affect its structure sheaf. Therefore, 
the degrees of the Hodge bundles on the two surfaces $X$ and $\widetilde{X}$
will be the same: $\lambda_{\widetilde{X}}=\lambda_X$. 
On the other hand, blowing up a smooth point on a surface 
decreases the square of its dualizing sheaf by 1, while there is no effect when
blowing up a rational double point. Each type of singular fibers $C$ in $X$
requires apriori 
different quasi-admissible modifications (or no modifications at all),
and thus decreases $\kappa_X$ by some nonnegative integer, denoted by
$\mu(C)$: 
\begin{equation}
\kappa_X=\kappa_{\widetilde{X}}+\sum_{C}\mu(C).
\end{equation}
Thus, $\mu(C)$ counts the number
of ``smooth blow-ups'' on $C$, which are needed to obtain
the minimal 
quasi-admissible cover $\widetilde{C}\rightarrow C$
within the surface quasi-admissible cover $\widetilde{\phi}:
\widetilde{X}\rightarrow \widetilde{Y}$.

\smallskip In the following Lemma, we compute the invariants $\mu(C)$
only for the general members of the boundary $\Delta{\mathfrak{T}}_g$
(cf.~Fig.~\ref{Delta-k,i}). The remaining, more special, singular curves
in $\Delta{\mathfrak{T}}_g$ will be linear combinations of these $\mu(C)$'s
(cf.~Sect.~11).

\begin{lem} If $\mu_{k,i}$ denotes the invariant $\mu(C)$ for a general
curve $C\in\Delta\mathfrak{T}_{k,i}$, then
\begin{eqnarray*}
\on{(a)}&&\mu_0=\mu_{1,i}=\mu_{4,i}=\mu_{6,i}=0;\\
\on{(b)}&&\mu_{2,i}=1;\\
\on{(c)}&&\mu_{3,i}=\mu_{5,i}=2.
\end{eqnarray*}
\label{mu(C)}\vspace*{-5mm}
\end{lem}
\noindent{\it Proof.} The general members of the boundary 
$\Delta{\cal H}$ are in fact the minimal quasi-admissible covers associated to 
the general members of the boundary $\Delta{\mathfrak{T}}$, except for
$\Delta_0$ which has $\mu_0=0$. Thus, we trace the blow-ups necessary to
transform the curves in Fig.~\ref{Delta-k,i} to the curves in
Fig.~\ref{admissible-k,i}. For example, no blow-ups are needed in the case
of $\Delta_{1,i}$, so that $\mu_{1,i}=0$, while we need 2 blow-ups in
the case of $\Delta_{3,i}$, and hence $\mu_{3,i}=2$.

\smallskip
The only interesting situation occurs for $\Delta_{5,i}$. Apparently,
there is only {\it one} added component $\proj ^1$ to the original
$C\in\Delta_{3,i}$, but the lemma states that $\mu_{3,i}=2$. The difference
comes from the fact that near the intersection $r=C\cup\proj^1$
the surface $\widetilde{X}$ has equation $xy=t^2$, i.e. $r$ is a rational
double point on $\widetilde{X}$ of type $A_1$ (a similar situation occurred
in Fig.~\ref{mult4}). To obtain such a point $r$ in place of a smooth
point $r_1$ on $X$, we first blow up $r_1$, and then on the
obtained exceptional divisor we blow up another point $r_2$, so as to end
with a {\it chain of two} $\proj^1$'s (cf.~Fig.~\ref{mu-5,i}). Finally,
we blow down the first $\proj^1$, and develop the required rational
double point $r$. As a result, we have two ``smooth'' and one
``singular'' blow-ups, which implies $\mu_{5,i}=2.$ \qed

\bigskip
\begin{figure}[h]
$$\psdraw{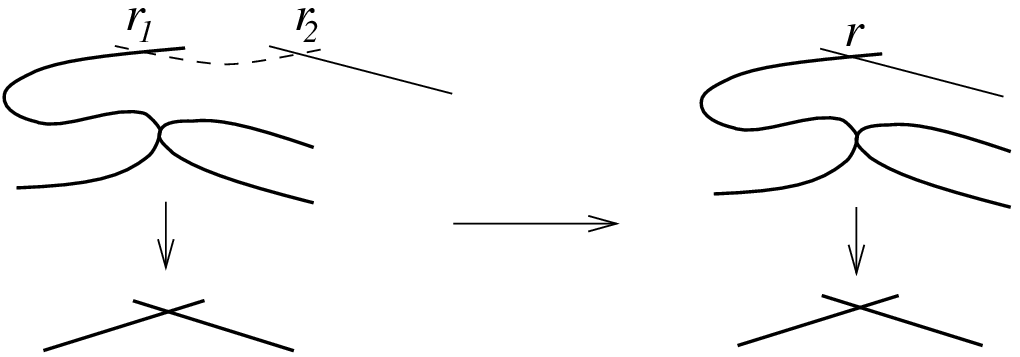}{3in}{1in}$$
\caption{Quasi-admissible blow-ups on $\Delta_{5,i}$}
\label{mu-5,i}
\end{figure}

\section*{\hspace*{1.9mm}5. Effective Covers}

\setcounter{section}{5} 
\setcounter{subsection}{0} 
\setcounter{subsubsection}{0} 
\setcounter{lem}{0}
\setcounter{thm}{0}
\setcounter{prop}{0}
\setcounter{defn}{0} 
\setcounter{cor}{0} 
\setcounter{conj}{0} 
\setcounter{claim}{0} 
\setcounter{remark}{0}
\setcounter{equation}{0}
\label{effectivecovers}

In this section we construct the final type of triple covers in
the Basic Construction. These will not be necessary for the global
calculation in Section~7, so the reader may wish to skip
this more technical part on a first reading, and assume in
Section~6 that all covers are flat.

\subsection{Construction of
effective covers $\widehat{X}\rightarrow\widehat{Y}$}
\label{constructioneffective} Consider a quasi-admissible cover
$\widetilde{\phi}:\widetilde{X}\rightarrow \widetilde{Y}$, as given in
Prop.~\ref{propquasi}. In order to use the fact that the pushforward
$\widetilde{\phi}_*{\cal{O}_{\widetilde{X}}}$
is locally free on $\widetilde{Y}$, we need to assure that the map
$\widetilde{\phi}$ is {\it flat}. Unfortunately, there are certain points
on $\widetilde{X}$ where this fails to be true: exactly where the fibers
of $\widetilde{X}$ are ramified as triple covers of the corresponding
fibers of $\widetilde{Y}$. The situation can be resolved by several further
blow-ups. 

\smallskip
Namely, we work locally analytically near the points in $\widehat{X}$
of ramifications index 1 or 2, and consider correspondingly two cases.

\subsubsection{Case of ramification index 1} 
\label{caseram1}
\begin{figure}[t]
$$\psdraw{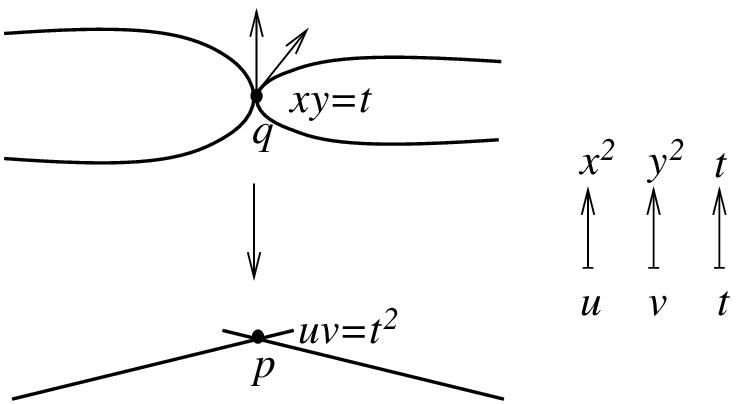}{2.7in}{1.2in}$$
\caption{Non-flat cover of ramification  index 1}
\label{ram}
\end{figure}

This case involves
members of the boundary divisors $\Delta{\mathfrak{T}}_{4,i}$ and 
$\Delta{\mathfrak{T}}_{5,i}$. Let $q$ be the point of ramification in the
fiber of $\widetilde{X}$ over the point $p$ in the fiber of $\widetilde{Y}$
(cf.~Fig.~\ref{ram}).

We use the pull-back
of the map $\widetilde{\phi}$ to study the embedding of the completion of
the local ring of $p$ into that of $q$:
\[\widehat{\cal{O}}_{\widetilde{Y}\!,p}=\mathbb{C}[[u,v,t]]
\big/_{\displaystyle{(uv-t^2)}}
\stackrel{\widetilde{\phi}^{\#}}
{\hookrightarrow}\widehat{\cal{O}}_{\widetilde{X}\!,q}=
\mathbb{C}[[x,y,t]]\big/_{\displaystyle{(xy-t)}}.\]

\noindent Therefore, as an $\widehat{\cal{O}}_{\widetilde{Y}\!,p}$-module,
\[\widehat{\cal{O}}_{\widetilde{X}\!,q}= \widehat{\cal{O}}_{\widetilde{Y}\!,p}+
\widehat{\cal{O}}_{\widetilde{Y}\!,p}x+
\widehat{\cal{O}}_{\widetilde{Y}\!,p}y.\]
However, this is not a locally-free 
$\widehat{\cal{O}}_{\widetilde{Y},p}$-module: for instance, one
relation among the generators is $(v-t)x+(u-t)y=0$.

\smallskip
Alternatively, the fiber of $\phi$ over $p$ is supported at $q$, but 
it is of degree 3 rather than 2,
which would have been necessary for the flatness of a degree $2$ map.
 Indeed, as $\mathbb{C}-$vector spaces:
\[\widehat{\cal{O}}_{\widetilde{X}\!,q}\otimes_{\widehat{\cal{O}}_
{\widetilde{Y}\!,p}}
\on{Spec}k(p)\cong \widehat{\cal{O}}_{\widetilde{X}\!,q}\big/_{\displaystyle
{\widehat{\mathfrak{m}}_{Y\!,p}\widehat{\cal{O}}_{\widetilde{X}\!,q}}}\cong
\mathbb{C}[[x,y]]\big/_{\displaystyle{(x^2,y^2,xy)}}=\mathbb{C}\oplus
\mathbb{C}x\oplus\mathbb{C}y.\]
In Fig.~\ref{ram} one can visually observe the two distinct tangent
directions at $q$ making it a {\it fat} point of degree $3$.

\medskip
We conclude that $\widetilde{\phi}$ is indeed
non-flat at $q$. To resolve this, we blow-up $\widetilde{Y}$ at $p$
and $\widetilde{X}$ at $q$, denoting the new surfaces by $\widehat{Y}$ and
$\widehat{X}$. It is easy to see that they fit into the following coming
diagram:

\smallskip
\begin{figure}[h]
$$\psdraw{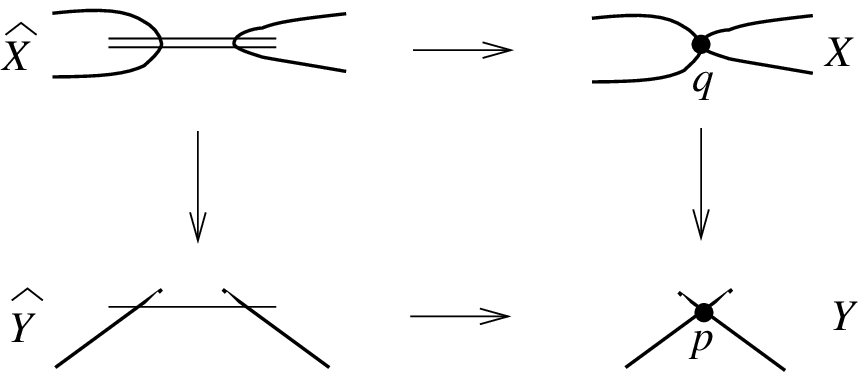}{3.5in}{1.4in}$$
\caption{Resolving the case of ram. index 1}
\label{resolve1}
\end{figure}

\smallskip
In order to keep the map $\widehat{\phi}:\widehat{X}\rightarrow \widehat{Y}$
of degree 3, we need to blow-up one further point on $\widetilde{X}$:
if the inverse image of $p$ is $\{q,r\}$ we blow-up $r$, and thus we add
the necessary component to $\widehat{X}$ to make it a triple cover of
$\widehat{Y}$ (cf.~Fig~\ref{coef2.fig}).

\medskip
\subsubsection{Case of ramification index 2} 
\label{caseram2}
\begin{figure}[h]
$$\psdraw{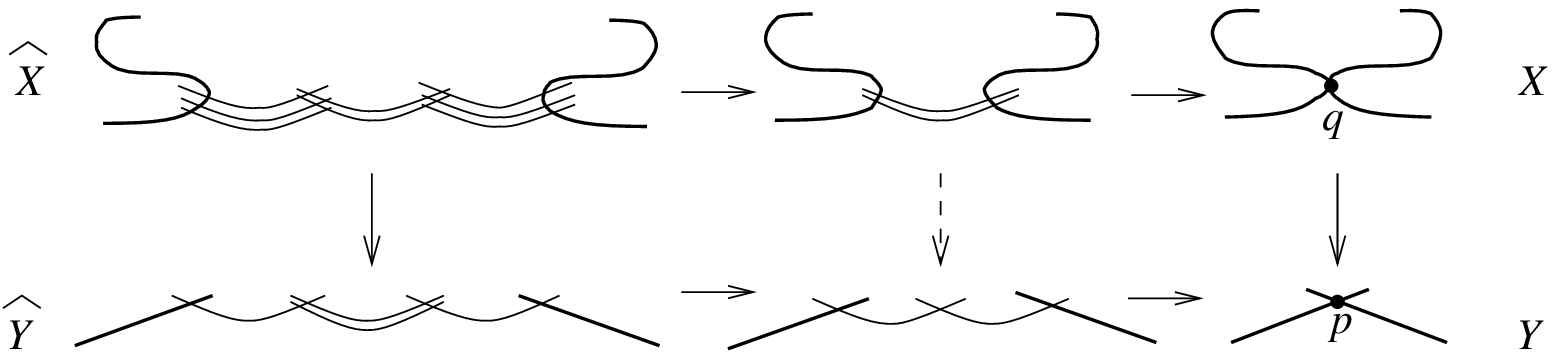}{5.3in}{1.4in}$$
\caption{Resolving the case of ram. index 2}
\label{resolve2}
\end{figure}

The only boundary component,
where ramification index 2 occurs, is $\Delta{\mathfrak{T}}_{6,i}$.
Similarly as above, $\widetilde{\phi}:
\widetilde{X}\rightarrow \widetilde{Y}$ is non-flat at $q$. Indeed,
$\widehat{\cal{O}}_{\widetilde{X}\!,q}$ is generated as an 
$\widehat{\cal{O}}_{\widetilde{Y}\!,p}$-module 
by $1,x,y,x^2,y^2$, but not-freely
due to the relation $u\cdot x+v\cdot y-t\cdot x^2-t\cdot y^2=0$. To 
resolve the apparent non-flatness of $\widetilde{\phi}$, we can blow-up
once $\widetilde{X}$ and $\widetilde{Y}$ at 
$q$ and $p$, but this would not be sufficient. In fact, we must make
further blows-ups on each surface, as Fig.~\ref{resolve2} suggests:
two more on $\widetilde{X}$ and one more on $\widetilde{Y}$.

\medskip

\medskip
In both cases of ramification index 1 or 2,
the new map $\widehat{\phi}:\widehat{X}\rightarrow
\widehat{Y}$ is obtained from $\widetilde{\phi}$ by a base change, and hence 
$\widehat{\phi}$
is {\it proper} and {\it finite}, and by construction, also a {\it flat}
morphism. We call such covers {\it effective}.

\medskip
The above considerations combined with Prop.~\ref{quasicov}
imply the existence of effective covers for our families of trigonal curves:

\begin{prop}  Let $X/\!_{\displaystyle{B}}$ be a family of
trigonal curves with smooth general member. After several blow-ups (and
possibly modulo a base change) we can associate to it an effective cover 
$\widehat{\phi}:\widehat{X}\rightarrow \widehat{Y}$.
\label{effexist}
\end{prop}

Here $\widehat{Y}$ is a birationally ruled surface over $B$. 
If the base curve $B$ is {\it not} tangent to the boundary divisors $\Delta
{\mathfrak{T}}_{k,i}$, then the resulting surfaces $\widehat{X}$ and
$\widehat{Y}$ will have smooth total spaces. If, moreover, $B$ intersects
the $\Delta_{k,i}$'s only in their general points (as given in
Fig.~\ref{Delta-k,i}), then the special fibers of $\widehat{Y}$ and
$\widehat{X}$ are easy to describe
(cf.~Fig.~\ref{coef1.fig}-\ref{coef3.fig}). For example, $\widehat{Y}$'s
special fibers are either
chains of two or three reduced projective lines, or chains of five
smooth rational curves with non-reduced middle component of
multiplicity two. The special fibers of $\widehat{X}$ can also contain 
nonreduced components (of multiplicity 2 or 3), and this occurs only in the
ramification cases discussed above ($\Delta{\mathfrak{T}}_{k,i}$ for $k=4,5,6$).

\subsection{Change of $\lambda_X,\kappa_X$ and $\delta_X$ in
the effective covers}
\label{change}
This is an analog to the discussion in Section~\ref{The invariants}.
After the necessary base changes  
we again identify, without loss of generality,
the new base curve $\widetilde{B}$ with $B$, 
and the pull-back of $X$ over $\widetilde{B}$ with $X$, and we redefine
the basic invariants $\lambda_{\widehat{X}}$ and $\kappa_{\widehat{X}}$
for the effective family $\widehat{X}$ over $\widetilde{B}$.
(It doesn't
make sense to define directly $\delta_{\widehat{X}}$, because of the
nonreduced fiber components in $\widehat{X}$. We could, of course, set
$\delta_{\widehat{X}}=12\lambda_{\widehat{X}}-\kappa_{\widehat{X}}$,
but we will not need this in the sequel.) Now the original $X$ and
the effective $\widehat{X}$ differ by ``quasi-admissible'' and ``effective''
blow-ups. The connections between the invariants of $X$, $\widetilde{X}$ and
$\widehat{X}$ are expresssed by the following 

\begin{lem} With the above notation,
\begin{eqnarray*}
\on{(a)}&\!\!\!\!&\displaystyle{\lambda_X}=\lambda_{\widetilde{X}}=
\lambda_{\widehat{X}};\\
\on{(b)}&\!\!\!\! &\kappa_X=\kappa_{\widetilde{X}}+\sum_{C}\mu(C);\\
\on{(c)}&\!\!\!\! &\displaystyle{\kappa_{\widetilde{X}}=
\kappa_{\widehat{X}}+\sum_{\on{ram}1}1+\sum_{\on{ram}2}3}.
\end{eqnarray*}
\label{changeinv}\vspace*{-10mm}
\end{lem}

\begin{proof} In view of Lemma~\ref{mu(C)}, the first and the
second statements are obvious. Obtaining a flat cover
$\widehat{X}\rightarrow \widehat{Y}$ requires blowing up on $\widetilde{X}$
one smooth point for each 
ramification index 1, and three smooth points for each ramification 
index 2. Hence the relation between $\kappa_{\widehat{X}}$ and
$\kappa_{\widetilde{X}}.$ \end{proof}

\section*{\hspace*{1.9mm}6. 
Embedding $\widehat{X}$ in a Projective Bundle over $\widehat{Y}$}

\setcounter{section}{6} 
\setcounter{subsection}{0} 
\setcounter{subsubsection}{0} 
\setcounter{lem}{0}
\setcounter{thm}{0}
\setcounter{prop}{0}
\setcounter{defn}{0} 
\setcounter{cor}{0} 
\setcounter{conj}{0} 
\setcounter{claim}{0} 
\setcounter{remark}{0}
\setcounter{equation}{0}
\label{embedding}

Given the effective degree $3$ map $\widehat{\phi}:\widehat{X}\rightarrow
\widehat{Y}$, our next step is to embed 
$\widehat{X}$ into a projective bundle $\proj V$ of rank $1$ over
the birationally ruled surface $\widehat{Y}$. We shall consider a degree 
$3$ morphism $\widehat{\phi}$, but the same
discussion is valid for any degree $d$.

\subsection{Trace map}
\label{tracemap} Since $\widehat{\phi}$ is flat and finite,
the pushforward $\widehat{\phi}_*(\cal O_{\widehat{X}})$
is a locally free sheaf on $\widehat{Y}$ of rank $3$.
Define the {\it trace} map
\[\on{tr}:\widehat{\phi}
_*(\cal O_{\widehat{X}})\rightarrow \cal O_{\widehat{Y}}\]
as follows. The finite field extension $K(\widehat{X})$ of
$K(\widehat{Y})$ induces the {\it algebraic} trace map
$\on{tr^\#}:K(\widehat{X})\rightarrow K(\widehat{Y})$, defined by 
$\on{tr^\#}(a)=\textstyle{\frac{1}{3}}(a_1+a_2+a_3)$. Here
the $a_i$'s are the conjugates of $a$ over $K(\widehat{Y})$ 
in an algebraic closure of $K(\widehat{X})$. The restriction
 $\on{tr^\#}{|}_{K(\widehat{Y})}=\on{id}_{K(\widehat{Y})}$.
Over an affine open $U={\on {Spec}}\,A\subset
\widehat{Y}$ and its affine inverse $\widehat{\phi}^{-1}(U)=
{\on {Spec}}\,B\subset
\widehat{X}$, $B$ is the integral closure of $A$ in its field of fractions
$K(\widehat{X})$. Therefore, 
the trace map restricts to the $A$-module homomorphism
$\on{tr^\#}:B\rightarrow A$. We have a commutative diagram:

\medskip
\begin{figure}[h]
\vspace*{2mm}
\begin{picture}(4,2.5)(4,0)
\put(7.8,2.2){$B\hspace{1mm}\hookrightarrow \hspace{1mm}K(\widehat{X})$}
\multiput(8,1.2)(1.6,0){2}{\vector(0,1){0.9}}
\put(7.8,0.7){$A\hspace{1mm}\hookrightarrow \hspace{1mm}K(\widehat{Y})$}
\put(7.8,1.65){\oval(0.35,0.9)[l]}
\put(7.8,1.2){\vector(1,0){0.1}}
\put(7.8,2.1){\line(1,0){0.1}}
\put(9.8,1.65){\oval(0.35,0.9)[r]}
\put(9.8,1.2){\vector(-1,0){0.1}}
\put(9.8,2.1){\line(-1,0){0.1}}
\multiput(7,1.5)(3.1,0){2}{$\on{tr}^{\#}$}
\put(3.2,0.7){$U=\on{Spec}\,A$}
\put(2.3,2.2){$\widehat{\phi}^{-1}(U)=\on{Spec}\,B$}
\put(3.7,1.5){$\big\downarrow$}
\put(3.7,1.4){$\big\downarrow$}
\put(3.3,1.5){$\widehat{\phi}$}
\end{picture}
\vspace*{-6mm}
\caption{The trace map}
\end{figure}

The so-defined local maps $\on{tr}:\widehat{\phi}_*\cal O_{\on {Spec}\,B}
\twoheadrightarrow 
\cal O_{\on {Spec}\,A}$ patch up to give a global trace map
 $\on{tr}:\widehat{\phi}_*\cal O_{\widehat{X}}
\twoheadrightarrow \cal O_{\widehat{Y}}.$ 
Let $\cal V$ be the kernel of $\on{tr}$:
\begin{equation}
0\rightarrow {\cal V}\rightarrow {\widehat{\phi}}
_*{\cal O}_{\widehat X}\stackrel
{\on{tr}}{\rightarrow}{\cal O}_{\widehat Y}\rightarrow 0.
\label{splitting}
\end{equation}
Note that $\cal V$ is locally free of rank $2$. 
The natural inclusion
$\cal O_{\widehat{Y}}\hookrightarrow \widehat{\phi}
_*\cal O_{\widehat{X}}$, composed
with $\on{tr}$, is the identity on $\cal O_{\widehat{Y}}$, hence the
exact sequence splits:
\begin{equation}
{\widehat{\phi}}_*{\cal O}_{\widehat X}={\cal O}_{\widehat Y}\oplus {\cal V}.
\label{directsum}
\end{equation}

\smallskip
\subsubsection{Geometric interpretation of the trace map}
It is useful to interpret the trace map geometrically in terms of the 
corresponding vector bundles $\widehat{\phi}_*{O_{\widehat X}}$, $O_{\widehat Y}$ and
$V$ associated to the sheaves
$\widehat{\phi}_*{\cal O}_{\widehat X}$, ${\cal O}_{\widehat Y}$ and $\cal V$.
We again localize over affine opens,
and if necessary, we shrink $U=\on{Spec}\,A$ so that $\widehat{\phi}_*
\cal O_{\widehat X}$ becomes a {\it free} $\cal O_{\widehat Y}$-module.

\smallskip
Let $p$ be a closed point in $\on{Spec}\,A$ with maximal ideal
$\mathfrak{p} \subset A$, having three distinct preimages $q,r,s\in\on{Spec}\,B$
with maximal ideals $\mathfrak{q,r,s}\subset B$. 
Since $B$ is a free $A$-module, the
quotient $B/{\mathfrak{p}}B$ is a 3-dim'l algebra over the ground
field $\on{k}(p)=A/{\mathfrak{p}}$, i.e. a 3-dim'l vector space over ${\mathbb C}$. 
On the other hand, from $\mathfrak{qrs}=\mathfrak{q}\cap
\mathfrak{r}\cap \mathfrak{s}$ and the Chinese Remainder Theorem, it is clear that 
$B/{\mathfrak{p}}B\cong B/{\mathfrak{q}} \oplus B/{\mathfrak{r}} \oplus
B/{\mathfrak{s}}\cong {\mathbb C}\overline{f}_q\oplus{\mathbb
C}\overline{f}_r\oplus{\mathbb C}\overline{f}_s.$
The generators $\overline{f}_q,\overline{f}_r,\overline{f}_s$ are
chosen as usual: $\overline{f}_q$, for instance,
 is the residue in $\on{k}(q)$
of a function $f_q\in B$ such that $f_q\equiv 1(\on{mod}\mathfrak{q}),\,\,
f_q\equiv 0(\on{mod}\mathfrak{r,s})$. 

\smallskip
In the Groethendieck style, the vector bundle
over $\widehat Y$ associated to $\widehat{\phi}_*\cal O_{\widehat X}$ is 
$\on{ Spec}\on{S}(B_A)$, where $\on{S}(B_A)$ is the symmetric algebra
of $B$ over $A$. Its fiber over $p$ is the pull-back 
$\on{ Spec}\on{S}(B_A)\times _{\on{Spec\,A}}\on{Spec\,k}(p) =
\on{ Spec}(\on{S}(B_A)\times_A A/\mathfrak {p})=\on{Spec}\on{S(B}/{\mathfrak{p}}B).$
We prefer to work with the dual $\widehat{\phi}_*{O_{\widehat X}}$
of this bundle, and the same goes for
projectivizations: we projectivize the 1-dim'l subspaces of
$\widehat{\phi}_*{O_{\widehat X}}$ rather than its 1-dim'l quotients. 

In view of this convention, the fiber of the bundle
$\widehat{\phi}_*{O_{\widehat X}}$ is canonically identified as
\[(\widehat{\phi}_*{O_{\widehat X}})_p=B/{\mathfrak p}B\cong {\mathbb C}\overline{f}_q
\oplus{\mathbb
C}\overline{f}_r\oplus{\mathbb C}\overline{f}_s \cong {\bf A}^3_{\mathbb C}.\]
The generators $\overline{f}_q,\overline{f}_r,\overline{f}_s$
span three lines in ${\bf A}^3_{\mathbb C}$, which can be
canonically described: the line $\Lambda_q=
{\mathbb C}\overline{f}_q$, for example, corresponds
to all functions regular at $q,r$ and $s$, and vanishing at $r$ and $s$.

\smallskip
Similarly, the vector bundle $O_{\widehat Y}$ associated to 
$\cal O_{\widehat Y}$ has fiber $(O_{\widehat Y})_p=A/{\mathfrak {p}}
\cong {\mathbb C}\overline{h}_p$, where $h_p$ is a
function near $p$ having residue $h_p(p)=1$ in $\on{k}(p)$. The trace
map $\on{tr}:\widehat{\phi}_*\cal O_{\widehat{X}}
\twoheadrightarrow \cal O_{\widehat{Y}}$ translates fiberwise in terms
of the vector bundles $\widehat{\phi}_*{O_{\widehat X}}$ and 
$O_{\widehat Y}$ as:
\[\on{tr}_p:{\mathbb C}\overline{f}_q\oplus{\mathbb
C}\overline{f}_r\oplus{\mathbb C}\overline{f}_s \rightarrow 
{\mathbb C}\overline{h}_p,\,\,\overline{f}\mapsto \frac{1}{3}
(f(q)+f(r)+f(s)).\] 

Finally, the locally free subsheaf ${\cal V}=\on{Ker(tr)}\subset
\widehat{\phi}_*{\cal
O}_{\widehat X}$ is associated to a vector bundle $V$ with fiber
$V_p=\{\overline{f}\,\,|\,\,f(q)+f(r)+f(s)=0\}
\subset (\widehat{\phi}_*{O_{\widehat X}})_q.$
Equivalently, from the direct sum (\ref{directsum}), 
$V_p=(\widehat{\phi}_*{O_{\widehat X}})_p \big{/}_{\textstyle{\Lambda}}$, where the line
$\Lambda=\{\overline{f}\,\,|\,\,f(q)=f(r)=f(s)\}$
corresponds to pull-backs of functions regular at $p$.

\begin{figure}[h]
$$\psdraw{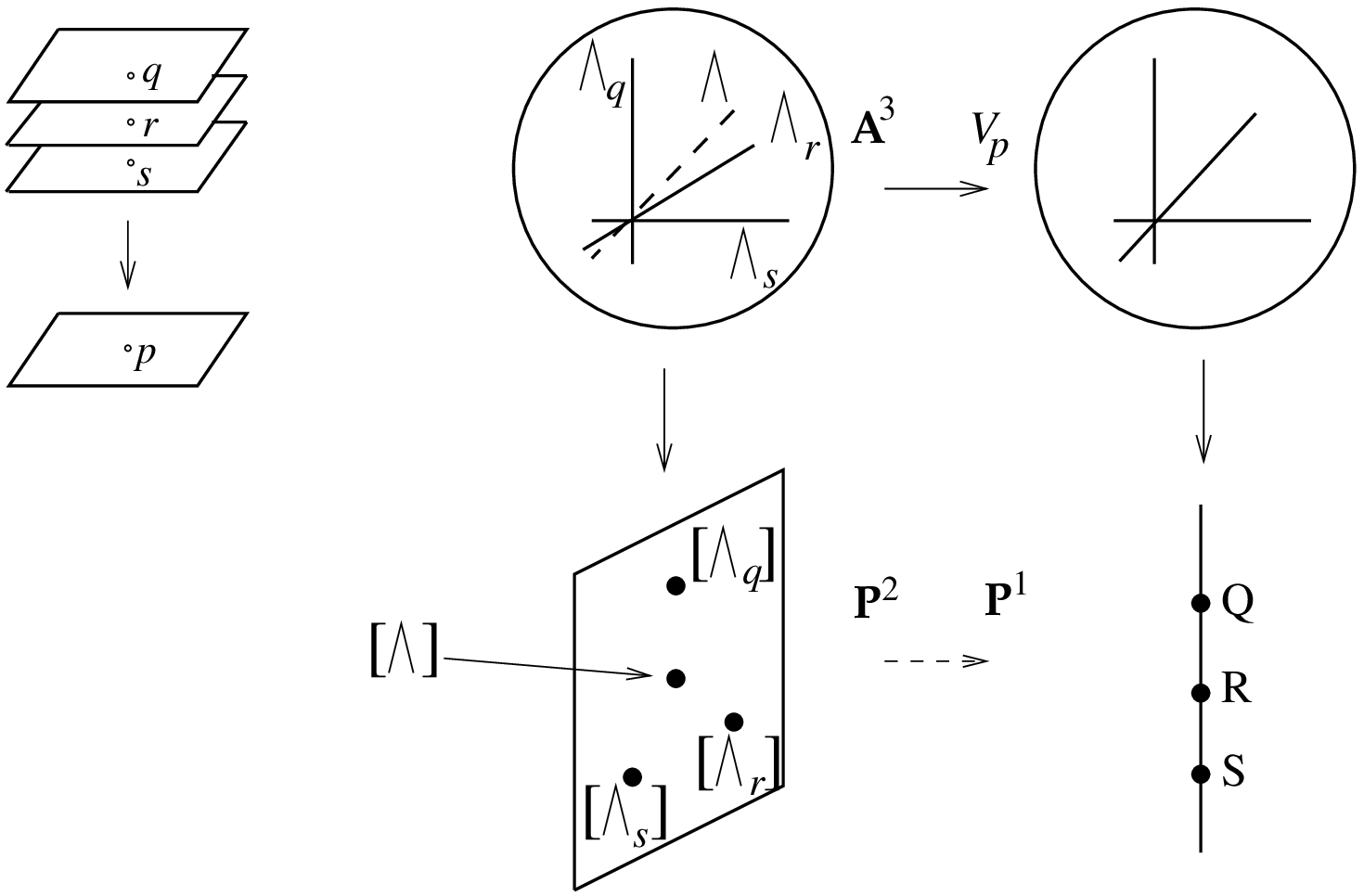}{3in}{1.88in}$$
\caption{Geometric interpretation of $tr$}
\label{geometric}
\end{figure}

 Since the four lines $\Lambda_q,\Lambda_r,\Lambda_s$
and $\Lambda$ are in general position in the fiber $(\widehat{\phi}
_*{O_{\widehat X}})_p$, modding out by $\Lambda$ yields three
distinct lines in the fiber $V_p$ (cf.~Fig.~\ref{geometric}). 
Therefore, projectivizing $V_p$ naturally induces three distinct points
$Q,R,S$ in the fiber ${\proj}^1$ of ${\proj}V$.
Going the other way around the diagram, we first projectivize
$(\widehat{\phi}_*{O_{\widehat X}})_q\cong{\bf A}^3$, 
and then we project from the point 
$[\Lambda]$ onto the fiber of ${\proj}V$. In other words, $\pi_{[\Lambda]}
:{\proj}^2\dashrightarrow {\proj}^1$ is well-defined at
$[\Lambda_q],[\Lambda_r]$ and $[\Lambda_s]$.

\medskip This completes the interpretation of the trace map in the case
of three distinct preimages  $q,r,s$ in $\widehat{X}$.
In case of only two or one preimage of $p$ in $\widehat{X}$, one modifies
 correspondingly the above interpretation.

\subsection{$\widehat{X}$ embeds naturally in 
$\proj V$ over $\widehat{Y}$}
\label{naturally}
We construct the map $i:\widehat{X}\hookrightarrow \proj V$ via the use of
 an invertible relative dualizing sheaf $\omega_{\widehat X/ \widehat Y}$. 
Its existence imposes a mild technical condition on the schemes $\widehat X$ 
and $\widehat Y$: we assume that they are Gorenstein, i.e. Cohen-Macaulay
with invertible dualizing sheaves $\omega_{\widehat X/{\mathbb C}}$ and 
$\omega_{\widehat Y/{\mathbb C}}$.
In our situation this will be sufficient. As we noted in
Section~\ref{constructioneffective}, when the base curve $B$ is {\it not}
tangent to the boundary divisors $\Delta{\mathfrak{T}}_{k,i}$, then 
$\widehat{X}$ and $\widehat{Y}$ are smooth surfaces. The remaining
cases are ``local'' base changes of these, and the construction
carries over.

\begin{prop} 
Let $\widehat{\phi}:\widehat X\rightarrow \widehat Y$ be a flat and finite
degree $d$ morphism of Gorenstein schemes, with $\widehat Y$ integral. Then
$\widehat{\phi}$ factors through a natural embedding of $\widehat  X$ into
 the projective bundle $\proj V$, followed by the projection $\pi:
\proj V\rightarrow \widehat Y$ (cf.~Fig.~\ref{basicconstruction}).
\label{propembedding}
\end{prop}

For easier referencing in the sequel, 
we have kept the notation $\widehat X$ and $\widehat Y$,
but the statement is true for {\it any} schemes
satisfying the Gorenstein condition. For another proof of
Prop.~\ref{propembedding}, see \cite{embedding}.

\medskip
\noindent{\it Proof of Prop.~\ref{propembedding}}.
Here we construct the map $i:\widehat X\rightarrow \proj V$, 
give the proof of its regularity, and point out how to show its injectivity.

\begin{figure}[h]
\begin{picture}(6,2.5)(-2.3,2.5)
\put(-1.7,4){$\proj({\cal O}_{\widehat Y})\stackrel 
{\on{tr}^{\#}}{\hookrightarrow}\proj({\widehat{\phi}}_*{\cal O}_{\widehat X})
\stackrel{\rho}{\dashrightarrow} \proj V$}
\put(0.95,3.05){\vector(0,1){0.8}}
\put(0.75,2.6){$\widehat X$}
\put(0.6,3.35){$\psi$}
\put(1.1,2.95){\vector(2,1){1.85}}
\put(1.7,3.35){$i$}
\end{picture}
\caption{Embedding $i:\widehat X\hookrightarrow \proj V$}
\label{construction of i}
\end{figure}

\subsubsection{Construction of the embedding map}
 According to Prop. II.7.12 
\cite{Hartshorne}, to give
a morphism $\psi:\widehat X\rightarrow {\proj}(\widehat{\phi}
_*(\cal O_{\widehat X}))$ over $\widehat Y$ is equivalent to give an 
invertible sheaf $\cal L$ on $\widehat Y$ and a surjective map of sheaves
$\widehat{\phi}^*(\widehat{\phi}_*(\cal O_{\widehat X}){\textstyle
{\widehat{\phantom{n}}}})\twoheadrightarrow  \cal L$. Recall from {\it
relative Serre-duality} that $(\widehat{\phi}_*\cal O_{\widehat X}){\textstyle
{\widehat{\phantom{n}}}}\cong \widehat{\phi}
_*\omega_{\widehat X/\widehat Y}$, and let
$\cal L=\omega_{\widehat X/\widehat Y}$. The natural morphism 
\[\sigma:\widehat{\phi}^*\widehat{\phi}
_*\omega_{\widehat X/\widehat Y}\rightarrow\omega_{\widehat X/\widehat Y}\]
is {\it surjective}. This is in fact true for any quasicoherent sheaf
$\cal F$ on $\widehat X$. Indeed, restricting to the
affine open sets $\widehat{\phi}:\on{Spec}\,B\rightarrow \on{Spec}\,A$, we have
$\cal F=M^{\sim}$ for some finitely generated $B$-module $M$, and
$\widehat{\phi}^*\widehat{\phi}
_*{\cal F}=\widehat{\phi}^*(M_A)^{\sim}=(M_A\otimes_A B)^{\sim}.$
The surjective $B$-module homomorphism $M_A\otimes_A B \twoheadrightarrow M$,
given by $m\otimes b \mapsto b\circ m$, induces
$\widehat{\phi}^*\widehat{\phi}_*{\cal F}\twoheadrightarrow \cal F$.

Thus, the  above map $\sigma$ naturally defines a morphism
$\psi:\widehat X\rightarrow {\proj}(\widehat{\phi}_*(\cal O_
{\widehat X}))$ over $\widehat Y$. Projectivizing
$0\rightarrow {\cal O}_{\widehat Y} \rightarrow {\widehat{\phi}}
_*{\cal O}_{\widehat X}\rightarrow {\cal V} \rightarrow 0$
gives a sequence of projective bundles, as in Fig.~\ref{construction of i}.
The map $\rho$ is undefined exactly on the image of $\on{tr}
^{\#}$. Composing $\rho$ with the map $\psi$ yields the
map $i:\widehat{X}\dasharrow \proj V$, which we claim is
a regular map.

\subsubsection{Regularity and injectivity of $i$.}
To see regularity, we show that the restriction
of $\sigma|_{\widehat{\phi}^*(\cal V{{\widehat{\phantom{n}}}})}
:{\widehat{\phi}^*(\cal V{{\widehat{\phantom{n}}}})}
\rightarrow \omega_{\widehat X/\widehat Y}$ is also
surjective. Indeed, we again work locally, and let $B=A\oplus C$
be the decomposition of $B$ via the trace map as a free $A$-module, where
$C=A\cdot b_1\oplus A\cdot b_2$ with $\on{tr}b_1=\on{tr}b_2=0$. 
Let $\omega_{\widehat X/\widehat Y}=M^{\sim}$ for some finitely
generated $B$-module $M$. Recall that
$\widehat{\phi}_*\omega_{\widehat X/\widehat Y}\cong 
(\widehat{\phi}_*\cal O_{\widehat X}){\textstyle{\widehat{\phantom{n}}}}$,
so that as  $A$-modules: $M\cong (B_A){{\widehat{\phantom{n}}}}=\on{Hom}_A
(B,A)$, and $B$ acts on $M$ by 
\[(b\circ f)(x)=f(bx)\,\,\,\on{for}\,\,\,f\in \on{Hom}_A(B,A),\,x\in B.\]
 Naturally, the sheaf $\cal V=C^{\sim}$, and $\widehat{\phi}^*
(\cal V{{\widehat{\phantom{n}}}})=(\on{Hom}_A(C,A)\otimes_A B)^{\sim}$, 
where we think of $f\in \on{Hom}_A(C,A)$ as an element
of $\on{Hom}_A(B,A)$ by extending it via $f(1)=0$. Our statement is
equivalent to showing that the $B$-module homomorphism
\[\sigma:\on{Hom}_A(C,A)\otimes_A B \rightarrow \big{(}\on{Hom}_A(B,A)
\big{)_B},\,\,f\otimes b \mapsto b\circ f,\]
is surjective. In fact, it suffices to show that the trace map is in the
image of $\sigma$, i.e. to find $c_1,c_2\in B$ such that 
\begin{equation}
\on{tr}\equiv c_1\circ {\pi_1}+c_2\circ {\pi_2}.
\label{trace equation}
\end{equation}
Here $\pi_j:B\rightarrow A$ 
gives the $b_j$-coordinate of $b\in B$, $j=1,2$. Set $c_1=b_1-\pi_1(b_1^2)$ and
$c_2=-\pi_1(b_1b_2)$. Evaluating both sides of (\ref{trace equation})
at $1,b_1$ and $b_2$ yields the same result, hence the identity is 
established, and $\sigma|_{\widehat{\phi}^*(\cal V{{\widehat{\phantom{n}}}})}
:{\widehat{\phi}^*(\cal V{{\widehat{\phantom{n}}}})}
\rightarrow \omega_{\widehat X/\widehat Y}$ is surjective.

\smallskip We have shown that the composition $\rho\circ \psi=
i:\widehat{X}\dasharrow \proj V$ is a regular map on $\widehat{X}$.
Alternatively, one could employ the geometric interpretation
of the trace map. A {\it general} point $p\in {\widehat Y}$ has three preimages
$q,r,s$ in $\widehat X$, each of which defines canonically a distinct point
$[\Lambda_q],[\Lambda_r]$ or $[\Lambda_s]$ in the fiber 
of ${\proj}(\widehat{\phi}_*{\cal O}_{\widehat X})$. As we pointed above, the
projection  $\pi_{[\Lambda]}
:{\proj}^2\dashrightarrow {\proj}^1$ is well-defined at
$[\Lambda_q],[\Lambda_r]$ and $[\Lambda_s]$. But $\pi_{[\Lambda]}$ is 
precisely the fiberwise restriction of $\proj({\widehat{\phi}}_*
{\cal O}_{\widehat X})
\stackrel{\rho}{\dashrightarrow} \proj V$, which shows again that the
composition $i=\rho\circ\psi:\widehat{X}\dasharrow \proj V$ is regular
on an open set of $\widehat X$. One makes the necessary modifications in
the cases of fewer preimages of $p$ in $\widehat{X}$. 
Finally, one can show, using similar methods (either algebraically or
geometrically), that the map $i$ is also an embedding. \qed

\medskip
\noindent{\bf Remark 6.1.} Since the
general fiber $C$ of $\widehat{X}$ is a smooth trigonal curve, 
the restriction $i|_{\displaystyle{C}}$ embeds $C$ in a ruled surface ${\mathbf
F}_k$ over the corresponding fiber $F_{\widehat{Y}}=
\proj^1$ of $\widehat{Y}$, where ${\mathbf F}_k=\proj(V|_{F_{\widehat{Y}}})$.
In Section~\ref{Maroniinvariant}
we will describe how the ruled surface ${\mathbf F}_k$ varies
as the fiber $C$ varies in $\widehat{X}$.

 \bigskip
\section*{\hspace*{1.9mm}7. 
Global Calculation on a Triple Cover $X\rightarrow Y$}

\setcounter{section}{7} 
\setcounter{subsection}{0} 
\setcounter{subsubsection}{0} 
\setcounter{lem}{0}
\setcounter{thm}{0}
\setcounter{prop}{0}
\setcounter{defn}{0} 
\setcounter{cor}{0} 
\setcounter{conj}{0} 
\setcounter{claim}{0} 
\setcounter{remark}{0}
\setcounter{equation}{0}
\label{triplecover}

In this section we consider the simplest case of effective covers, namely, 
when the original family $X$ is itself a triple cover of
a {\it ruled surface} $Y$ over the base curve $B$. This happens exactly when
all fibers of $X$ are irreducible, and the 
linear system of $g^1_3$'s on the
smooth fibers extends over the singular fibers to base-point free
line bundles  of degree 3 with two linearly independent sections. As we saw in 
Section~\ref{quasi-admissible}, we can patch together
all these $g^1_3$'s in a line bundle ${\cal L}$ on the total
space of $X$. Thus, $X$ will map to ${\proj}(H^0(X,{\cal
L})^{\widehat{\phantom{n}}})$ via 
$\phi_{\cal L}$, and this map will factor through our
ruled surface $Y$ over $B$:

\begin{figure}[h]
$$\psdraw{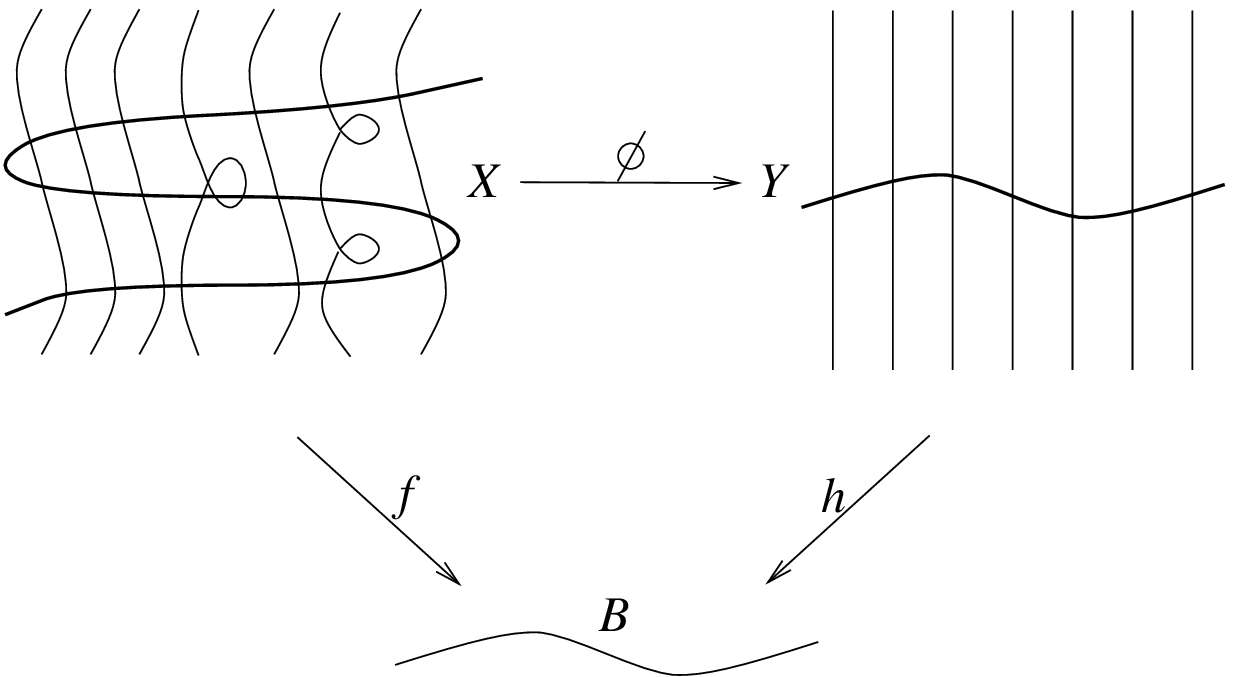}{2.7in}{1.31in}$$
\caption{Basic triple cover case.}
\label{fig.triplecover}
\end{figure}

Equivalently, we can describe such a 
family $X\rightarrow B$ via the classification of the boundary components
of the trigonal locus in Section~4: 
in $\overline{\mathfrak T}_g$ the base
curve $B$ meets only the boundary component 
$\Delta{\mathfrak{T}}_0$ of irreducible curves ($\delta_0|_B>0$), and 
there are no hyperelliptic fibers in $X$
($B\cap \overline{\mathfrak I}_g=\emptyset$).

\subsection{Global versus local calculation}
\label{global} As it will turn out, the
calculation of the slope $\delta_X /\lambda_X$ in this basic case yields
the actual maximal bound $\frac{36(g+1)}{5g+1}$: any addition of
singular fibers belonging to
other boundary components of $\overline{\mathfrak T}_g$ will
only decrease the ratio. Henceforth, we distinguish among two types of
calculation: {\it global} and {\it local}. The {\it global}
calculation refers to finding the
coefficients of $\delta_0$ and 
the Hodge bundle $\lambda|_{\overline{\mathfrak T}_g}$
in a relation in $\on {Pic}_{\mathbb Q}\overline{\mathfrak T}_g$ involving {\it all}
boundary classes.  The {\it local}
calculation provides the remaining coefficients by considering {\it local 
invariants} of each individual boundary class (cf.~Sect.~8).

\subsection{The basic construction}
\label{basic} For the remainder of this section, we consider a family
$X\rightarrow B$ such that, as above, $X$ is a triple cover of the 
corresponding ruled surface $Y$, and we carry out the proposed global 
calculation.

\unitlength 0.11in
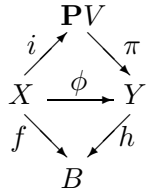
\begin{figure}[h]
\begin{picture}(10,9)(25,21)
\put(26.5,24){$X$}
\put(28.4,24.2){\vector(1,0){3.1}}
\put(32,24){$Y$}
\put(29,20){$B$}
\put(27.3,23.5){\vector(1,-1){2}}
\put(32.3,23.5){\vector(-1,-1){2}}
\put(26.6,22){$f$}
\put(31.75,22){$h$}
\put(29.45,24.7){$\phi$}
\put(27.3,25.4){\vector(1,1){2}}
\put(30.3,27.4){\vector(1,-1){2}}
\put(29,27.8){${\proj}V$}
\put(32,26.4){$\pi$}
\put(27.4,26.4){$i$}
\end{picture}
\caption{Triple Case}
\label{basicconstruction}
\end{figure}

Recall that
the pushforward of the structure sheaf ${\cal O}_X$ to $Y$ is a locally
free sheaf of rank $3$. In the exact sequence:
\[0\rightarrow E\rightarrow 
{\phi}_*{\cal O}_X\stackrel {\on{tr}}{\rightarrow}{\cal O}_Y\rightarrow 0,\]
the kernel of the trace map $\on{tr}$ is a vector
bundle $E$ on $Y$ of rank $2$,
 and $X$ naturally
embeds in the rank 1 projective bundle
${\proj}V$ over $Y$, where $V=E\,\widehat{\phantom{n}}$.
Any rank 2 vector bundle $E$ has the same
projectivizations as its dual bundle $V$ since $E\cong \bigwedge^2E\otimes V$,
where $\bigwedge^2E$ is a line bundle. For easier notation, in the trigonal
case we use the dual $V$ instead of $E$ from Section~\ref{embedding}.

A basis for $\on{Pic}Y$ can be chosen by letting $F_Y$ be the 
fiber of $Y$, and $B^{\prime}$ be any section of
$Y\rightarrow B$. Hence $\on{Pic}Y={\mathbb Z}B^{\prime}\oplus {\mathbb Z}F_Y$.
We normalize by replacing $B^{\prime}$ with
the ${\mathbb Q}$-linear combination
$B_0=B^{\prime}-\displaystyle{\frac{(B^{\prime})^2}{2}}{F_Y}$, and
provide a $\mathbb{Q}$-basis for $\on{Pic}_{\mathbb Q}Y$:
\begin{equation}
\on{Pic}_{\mathbb Q}Y={\mathbb Q}B_0\oplus {\mathbb Q}F_Y\,\,\,\on{with}\,\,\,
B_0^2=F_Y^2=0\,\,\,\on{and}\,\,\,B_0\cdot F_Y=1.
\label{normalize}
\end{equation} 
Let $\zeta$ denote the class of the
hyperplane line bundle ${\cal O}_{{\proj}V}(+1)$ 
on ${\proj}V$, and let $c_1(V)$ and $c_2(V)$ be the Chern classes 
of $V$ on $Y$. The Chow ring $A({\proj}V)$ is
generated as a $\pi^*(A(Y))$-module by $\zeta$
 with the only relation:
\begin{equation}
\zeta^2+\pi^*c_1(V)\zeta+\pi^*c_2(V)=0.
\label{zeta-relation}
\end{equation}
In particular, for the Picard groups:
\begin{equation}
\on {Pic}_{\mathbb Q}({\proj}V)=\pi^*(\on{Pic}
_{\mathbb Q}Y)\oplus {\mathbb Q}\zeta.
\label{Q-basis}
\end{equation}

\smallskip
\noindent As a divisor on ${\proj}V$, the surface
$X$ meets the fiber $F_{\pi}$ of $\pi$ 
generically in three points ($X$ maps three-to-one onto $Y$).
Thus in the Chow ring $A({\proj}V)$ we have $[X]\cdot [F_{\pi}]=3$, which
simply means that $X$ can be expressed as 
\[X\sim 3\zeta + \pi^*D\] 
for some divisor $D$ on $Y$ (see (\ref{Q-basis})).
With respect to the chosen basis for $\on{Pic}_{\mathbb Q}Y$:
\begin{equation}
D\sim aB_0+bF_Y\,\,\on{and}\,\, 
c_1(V)\sim cB_0+dF_Y
\label{define D,c1(V)}
\end{equation}
for some $a,b,c,d\in {\mathbb Z}$.
Note that $\deg(D|_{B_0})=b$ and $\deg(c_1(V)|_{B_0})=d$.

\subsection{Relation among the divisor classes $D$ and $c_1(V)$}
\label{relation}
It is evident that the divisors $D$ and $c_1(V)$ cannot be
independent on ruled surface 
$Y$ since both are canonically determined by the surface $X$.
The relation is in fact quite straightforward.

\begin{lem} With the above notation for the triple cover 
$\phi:X\rightarrow Y$, we have $D=2c_1(V)$ in $\on{Pic}Y$.
\label{D=2c_1(V)}
\end{lem}

\begin{proof} We start with the standard exact sequence for the divisor
$X$ on ${\proj}V$:
\begin{equation}
0\rightarrow \cal{O}_{{\proj}V}(-X)\rightarrow \cal{O}_{{\proj}V}
\rightarrow \cal{O}_X\rightarrow 0.
\label{X-divisorsequence}
\end{equation}
Pushing to $Y$ yields:
\begin{equation}
0 \!\rightarrow \!\pi_*\cal{O}_{{\proj}V}(-X) \!\rightarrow\! 
\pi_*\cal{O}_{{\proj}V}\!\rightarrow\! \pi_*\cal{O}_X \!\rightarrow \!
 R^1\pi_*\cal{O}_{{\proj}V}(-X) \!\rightarrow \!R^1
\pi_*\cal{O}_{{\proj}V} \!\rightarrow \cdots
\label{pushforward}
\end{equation}

\noindent It is easy to show that $R^1\pi_*\cal{O}_{{\proj}V}=0$ and
$\pi_*\cal{O}_{{\proj}V}(-X)=0$. This follows from
Grauert's theorem \cite{Hartshorne}:
$h^1(\cal{O}_{{\proj}V}|_{F_{\pi}})=h^1(\cal
O_{{\proj}^1})=0,$ and
\[h^0(\cal{O}_{{\proj}V}(-X)|_{F_{\pi}})=h^0(\cal{O}_{{\proj}V}
(-3\zeta-\pi^*D)|_{F_{\pi}})=h^0(\cal{O}_{{\proj}^1}(-3))=0.\]
Furthermore, $\pi_*\cal{O}_{{\proj}V}=
\cal{O}_Y$ and $\pi_*\cal{O}_X=\phi_*\cal{O}_X$, so that
(\ref{pushforward}) becomes
\begin{equation}
0\rightarrow \cal{O}_Y\rightarrow \phi_*\cal{O}_X\rightarrow 
R^1\pi_*\cal{O}_{{\proj}V}(-X)\rightarrow 0.
\label{remainingnonzero}
\end{equation}
From relative Serre-duality, and using the first Chern class of the 
relative dualizing sheaf, $c_1(\omega_{\pi})$ (cf.~(\ref{omega-pi})):
\[R^1\pi_*\cal{O}_{{\proj}V}(-X) \cong \big(\pi_*(\omega_{\pi}\otimes
\cal {O}_{{\proj}V}(X))\big)\widehat{\phantom{t}}=
\big(\pi_*\cal{O}_{{\proj}V}(\zeta+\pi^*D-\pi^*c_1(V))\big)
\widehat{\phantom{t}}.\]
\noindent
Since $\pi_*\cal{O}_{{\proj}V}(\zeta)=V\widehat{\phantom{t}}$ (cf.~
~\cite{BPV}),  
\begin{equation}
R^1\pi_*\cal{O}_{{\proj}V}(-X) \cong \big[V\widehat{\phantom{t}}
\otimes \cal{O}_Y(D-c_1(V))\big]\widehat{\phantom{t}}.
\label{tranformedsequence}
\end{equation}
Finally, combining (\ref{tranformedsequence}) with (\ref{remainingnonzero})
and $\phi_*\cal{O}_X/\cal{O}_Y=V\widehat{\phantom{t}}$, we arrive at
\[V\cong V\widehat{\phantom{t}}\otimes \cal{O}_Y(D-c_1(V))\,\,\Rightarrow\,\,
\cal{O}_Y(D-c_1(V))\cong \bigwedge ^2 V\cong\cal{O}_Y(c_1(V)),\]
and hence $D=2c_1(V)$ in $\on{Pic}Y. \,$ \end{proof}

\subsection{Global calculation of $\lambda_X$ and $\kappa_X$.}
\label{globalcalculation}
In the following proposition we express $\lambda_X$ and
$\kappa_X$ in terms of 
$\deg(c_1(V)|_{B_0})=d$ and the Chern class polynomial
$c_1^2(V)-4c_2(V)$, both of which are
independent of the choice of the vector bundle $V$. 
Indeed, if we twist
$V$ by a line bundle $M$ on $Y$ and set $V^{\prime}=V\otimes M$, then
\[c_1(V^{\prime})=c_1(V)+2c_1(M),\,\,
c_2(V^{\prime})=c_2(V)+c_1(V)c_1(M)+c_1(M)^2,\]
\[\Rightarrow\,\,c_1(V^{\prime})^2-4c_2(V^{\prime})=c_1(V)^2-4c_2(V).\]

 Recall the notation of (\ref{define D,c1(V)}).
In order to make $d$ also invariant, we use $b=2d$ from 
Lemma~\ref{D=2c_1(V)} and write $d=2b-3d$.
Now, if we replace ${\proj}V$ with its isomorphic
${\proj}V^{\prime}$ (cf.~Fig.~\ref{invariance}), and set $\zeta^{\prime}=i^*
\zeta\otimes(\pi^{\prime})^*M^{-1}$ to be the new hyperplane bundle,
then in $\on{Pic}({\proj}V)$: 
$X\sim 3\xi^{\prime}+(\pi^{\prime})^*D^{\prime}$ with $D^{\prime}\sim
D+3c_1(M)$. Hence
\[2D^{\prime}-3c_1(V^{\prime})=2D+6c_1(M)-3c_1(V)-6c_1(M)=2D-3c_1(V),\]
and equating their degrees on $B_0$, we obtain $2b^{\prime}-3d^{\prime}=2b-3d$.
\unitlength 0.11in
\begin{figure}[h]
\begin{picture}(10,7)(25,20.5)
\put(26,24){${\proj}V^{\prime}$}
\put(28.4,24.2){\vector(1,0){3.1}}
\put(32,24){${\proj}V$}
\put(29,20){$Y$}
\put(27.3,23.5){\vector(1,-1){2}}
\put(32.3,23.5){\vector(-1,-1){2}}
\put(26.5,22){$\pi^{\prime}$}
\put(31.75,22){$\pi$}
\put(29.45,24.7){$i$}
\end{picture}
\caption{$V^{\prime}=V\otimes M$}
\label{invariance}
\end{figure}
 
In other words, the following formulas for $\lambda_X$ and $\kappa_X$
would be valid for any vector bundle $V^{\prime}$ in place of the
canonically defined $V$, as long as the diagram of the basic construction
(cf.~Fig.~\ref{basicconstruction}) is satisfied, and as long as we adjust the
degree $d=\deg(c_1(V)|_{B_0})$ by its invariant form $2b-3d=
2\deg(D|_{B_0})-3\deg(c_1(V)|_{B_0})$.

\begin{prop}
Let $\phi:X\rightarrow Y$ be a degree 3 map from the original family
$X$ of trigonal curves to the ruled surface $Y$ over $B$. The
invariants $\lambda_X$ and $\kappa_X$ are given by the formulas:
\begin{eqnarray*}
\lambda_X&=&\displaystyle{\frac{g}{2}\deg\big(c_1(V)|_{B_0}\big)+\frac{1}{4}
\big(c_1(V)^2-4c_2(V)\big)},\\
\kappa_X&=&\displaystyle{\frac{5g-6}{2}\deg\big(c_1(V)|_{B_0}\big)+\frac{3}{4}
\big(c_1(V)^2-4c_2(V)\big)}.
\end{eqnarray*}
\label{lambda_X,kappa_X}\vspace*{-6mm}
\end{prop}
\noindent We defer the proof of Prop.~\ref{lambda_X,kappa_X}
to Subsections 7.4.2-3.

\subsubsection{Notation and Basic Tools.}
\label{notation} The proof of Prop.~\ref{lambda_X,kappa_X} 
consists of two calculations in the Chow ring
of ${\proj}V$; one uses versions of Riemann-Roch theorem on $X$ and ${\proj}V$,
and the other uses the adjunction formula on ${\proj}V$ for the
divisor $X$. Here we discuss these statements and set up the necessary
notation. 

\medskip
In order to work in ${\mathbb A}({{\proj}}V)$, we express the Chern classes of
${\proj}V$ in terms of the hyperplane class $\zeta$ and the Chern classes
of $Y$. We first recall that
$\pi_*\cal{O}_{{\proj}V}(+1)\cong V\widehat{\phantom{t}}$.
In the Euler sequence on ${{\proj}V}$:
\begin{equation}
0\rightarrow \cal{O}_{{\proj}V} \rightarrow \cal{O}_{{\proj}V}(+1)\otimes
\pi^*V\rightarrow \cal{T}_{\pi} \rightarrow 0,
\label{Eulersequence}
\end{equation}
we compare the Chern polynomials $c_t(\cal{O}_{{\proj}V}(+1)\otimes
\pi^*V)=c_t(\cal{T}_{\pi})$, and obtain:
\begin{eqnarray}
K_{{\proj}V}&\!\!\!=&\!\!\!\!-2\zeta-\pi^*c_1(V)+\pi^*K_Y,
\label{K_PV}\\
c_1(\omega_{\pi})&\!\!\!=&\!\!\!\!-2\zeta-\pi^*c_1(V),
\label{omega-pi}\label{omega_pi}\\
c_2({\proj}V)&\!\!\!=&\!\!\!\!-2\zeta\pi^*K_Y+\pi^*\big(c_1(V)K_Y+c_2(Y)\big).
\label{c_2(PV)}
\end{eqnarray}
Here $\cal{T}_{\pi}$ and $\omega_{\pi}$ are correspondingly the relative
tangent and the 
relative dualizing sheaves of $\pi$, while $K_{{\proj}V}$ is the 
class of the canonical sheaf on ${\proj}V$. On the ruled surface $Y$ over
the curve $B$ of genus $g_{\scriptscriptstyle B}$ we similarly have
\begin{eqnarray}
\hspace{9mm}K_Y&\!\!\!=&\!\!\!\!-2B_0+h^*(K_B)
\equiv -2B_0+(2g_{\scriptscriptstyle B}-2)F_Y
\label{K_Yglobal},\\
\hspace{9mm}c_2(Y)&\!\!\!=&\!\!4(1-g_{\scriptscriptstyle B}).
\label{c_2(Y)global}
\end{eqnarray}

\medskip
Now let $C$ be the general fiber of $X$, i.e. a smooth trigonal curve of
genus $g$. Assuming the Basic construction for the triple cover
$X\rightarrow Y$ (cf.~Fig.~\ref{basicconstruction}), 
we have the following lemmas.

\begin{lem}If $\chi(\cal{E})$ denotes the holomorphic Euler characteristic
of any sheaf $\cal{E}$, then the invariant $\lambda_X$ is expressible as
$\lambda_X=\chi({\cal O}_X)-\chi({\cal O}_B)\cdot 
\chi({\cal O}_C).$
\label{Euler}
\end{lem}

\noindent{\it Proof.} From
Grothendieck-Riemann-Roch theorem for the map $f:X\rightarrow B$,
\[\on{ch}(f_{!}\cal{O}_X).\on{td}\cal{T}_B=
f_*(\on{ch}\cal{O}_X.\on{td}\cal {T}_X),\]
where $\cal{T}_X$ and $\cal{T}_B$ are the corresponding tangent sheaves.
Since the fibers of $f$ are one-dimensional, $f_{!}\cal{O}_X=
f_*\cal{O}_X-R^1\!f_*\cal{O}_X=\cal{O}_B-(f_*\omega_f)\widehat{\phantom{t}}$. 
Substituting:
\[\big(1-g+c_1(f_*\omega_f)\big).\big(1-\frac{1}{2}K_B\big)=
f_*\big(1-\frac{1}{2}K_X+\frac{1}{12}(K^2_X+c_2(X))\big),\]
\[\Rightarrow\,\,
c_1(f_*\omega_f)=\frac{1}{12}f_*(K^2_X+c_2(X))-\frac{g-1}{2}K_B,\]
\[\,\,\,\,\,\,\,\,\,\,\Rightarrow\,\,
\lambda_X=\on{deg}(f_*\omega_f)=\chi(\cal{O}_X)-\chi(\cal{O}_B)\cdot
\chi(\cal{O}_C).\,\,\,\qed\]

\medskip
Note the similarity between this formula and the formula for $\delta_B$
in Example 2.1. Both quantities are expressed as differences of the
Euler characteristic (holomorphic or topological)
on the total space of $X$ and the product of the
corresponding Euler characteristics on the base $B$ and the general fiber $C$.
Lemma~\ref{Euler} suggests that in order to calculate $\lambda_X$,
we must have control over $\chi(\cal{O}_X)$. 

\begin{lem} In the Chow ring of ${\proj}V$:
\[\chi({\cal O}_X)=\frac{1}{12}X\big[\big(X+K_{{\proj}V}\big)
\big(2X+K_{{\proj}V}\big)+
c_2({\proj}V)\big]\]
\label{holomorphicEuler}\vspace*{-7mm}
\end{lem}
\begin{proof} From the standard exact sequence 
(\ref{X-divisorsequence}) for the divisor $X$ on ${\proj}V$ we have
$\chi({\cal O}_X)=\chi({\cal O}_{{\proj}V})-\chi({\cal O}_{\mathbf
{P}V}(-X))$. On the other hand, Hirzebruch-Riemann-Roch claims that
for any sheaf $\cal E$ on ${\proj}V$: 
$\chi(\cal{E})=\on{deg}\big(\on{ch}(\cal{E}).\on{td}\cal{T}_{{\proj}V}
\big)_3$. Applying this to
the line bundles $\cal{O}_{{\proj}V}$ and ${\cal O}_{{\proj}V}(-X)$, 
and subtracting the results completes the proof of the lemma. \end{proof}

\medskip
The reader may have noticed that all quantities discussed in the above
lemmas are elements 
of the third graded piece ${\mathbb A}^3({\proj}V)\otimes \mathbb{Q}$ of the Chow ring 
${\mathbb A}({\proj}V)\otimes \mathbb{Q}$. Hence they are cubic polynomials in the 
class $\zeta$, whose coefficients are appropriate products of pull-backs 
from ${\mathbb A}(Y)\otimes \mathbb{Q}$.
The higher degrees $\zeta^3$ and $\zeta^2$ can be decreased using the 
basic relation (\ref{zeta-relation}), while $\zeta$ itself can be
altogether eliminated by noticing that for any $\vartheta\in {\mathbb A}^2(Y)$:
\begin{equation}
\zeta.\pi^*(\on{point})=\zeta.F_{\pi}=1\,\,
\Rightarrow\,\,\zeta . \pi^*\vartheta=\on{deg}\vartheta.
\label{trivial}
\end{equation}
It is also useful to remember the trivial fact that for any divisors $D_i$
on $Y$, $\on{dim}Y=2$ implies
$D_1.D_2.D_3=0=D_i.c_2(V)$.  

\smallskip
\begin{figure}[h]
$$\psdraw{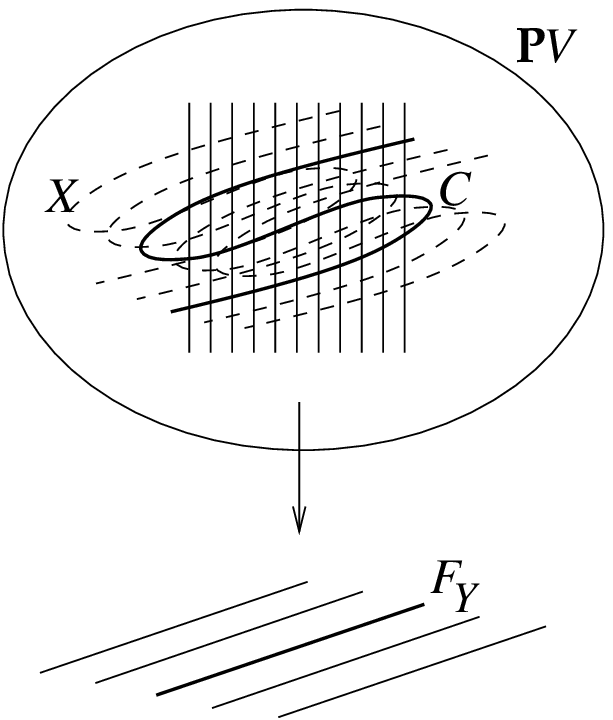}{1.8in}{1.8in}$$
\caption{Intersection of $X$ and $\pi^*F_Y$}
\label{intersection}
\end{figure}

\begin{lem}[Adjunction formula] The canonical bundle $\omega_Z$
of a smooth divisor $Z$ on the smooth variety $T$ can be expressed as
$\omega_Z\cong \omega_T\otimes \cal{O}_T(Z)\otimes
\cal{O}_Z$. Consequently, 
\[K_X^2=\big(K_{{\proj}V}+X\big)^2X\,\,\,{and}\,\,\,
g+2=\on{deg}c_1(V)|_{F_Y}.\]
\label{adjunction}\vspace*{-4mm}
\end{lem}

\noindent{\it Proof.} For the general statement of the adjunction formula
see \cite{Hartshorne}. The expression for $K_X^2$ is a straightforward
application to the divisor $X$ on ${\proj}V$:
$K_X=\big(K_{{\proj}V}+X\big)\big|_X$ is being squared in ${\mathbb A}({\proj}V)$.
As for the genus $g$ of the
general member $C$ of our family, we consider a general fiber
$F_Y$ of $Y$ (cf.~Fig.~\ref{intersection}).
Its pullback $\pi^*F_Y$ is a rational ruled surface $\mathbf F$
over $F_Y$, embedded in the 3-fold ${\proj}V$. The intersection of $\mathbf F$
with the surface $X$ is the trigonal fiber $C=X\cdot\pi^*F_Y=
(3\zeta+2\pi^*c_1(V))|_{\pi^*F_Y}$.

From the adjunction formulas
for $C\subset \pi^*F_Y$ and $\pi^*F_Y\subset {\proj}V$:
\begin{eqnarray*}
2g-2&\!\!=&\!\!(K_{\pi^*F_Y}+C)\cdot C=\big((K_{{\proj}V}+\pi^*F_Y)\big|_{\pi^*F_Y}
+X\big|_{\pi^*F_Y}\big)\cdot X\big|_{\pi^*F_Y}\\
    &\!\!=&\!\!\big(\zeta+\pi^*c_1(V)+\pi^*K_Y+\pi^*F_Y\big)\big(
    3\zeta+2\pi^*c_1(V)\big)\cdot \pi^*F_Y\\
    &\!\!=&\!\!(2c_1(V)+3K_Y)\cdot
F_Y=2\on{deg}c_1(V)\big|_{F_Y}-6. \,\,\,\qed
\end{eqnarray*}

\subsubsection{Global Calculation of $\lambda_X$.}
\label{globallambda} 
We substitute in Lemma~\ref{holomorphicEuler}
the expressions (\ref{K_PV}--\ref{c_2(PV)}) for $X,\,\,
K_{{\proj}V}$ and $c_2(V)$, as well as the identity $D=2c_1(V)$:
\[\chi(\cal{O}_X)=\frac{3\xi+2\pi^*c_1(V)}{12}\big[
\big(\xi+\pi^*c_1(V)+\pi^*K_Y\big)
\big(4\xi+3\pi^*c_1(V)+\pi^*K_Y\big)\]
\[-2\xi\pi^*K_Y+\pi^* c_1(V)\pi^*K_Y+\pi^*c_2(Y)\big].\]
Applying the necessary reductions, we arrive at:
\[\chi({\cal O}_X)=\frac{1}{2}\big(c_1^2(V)-2c_2(V)\big)+\frac{1}{2}c_1(V)K_Y+
\frac{1}{4}\big(K_Y^2+c_2(Y)\big).\]
We expect our formula for $\lambda_X$ to be independent of the base curve
$B$. The contribution of $g_{\scriptscriptstyle B}$ in $\chi(\cal{O}_X)$ can be
written as: $(g_{\scriptscriptstyle B}-1)\on{deg}c_1(V)\big|_{F_Y}+\chi(\cal{O}_Y)=
(g_{\scriptscriptstyle B}-1)(g-1)$,
but this is precisely the adjustment $\chi(\cal{O}_B)\chi(\cal{O}_C)$
given in Lemma~\ref{Euler}. Thus,
\[\lambda_X=\frac{1}{2}\big(c_1^2(V)-2c_2(V)\big)-\on{deg}c_1(V)\big|_{B_0}.\]
It remains to notice that $c_1^2(V)=2\on{deg}c_1(V)|_{F_Y}
\on{deg}c_1(V)|_{B_0}=2(g+2)\on{deg}c_1(V)|_{B_0}$
and rewrite $\lambda_X$ in the form
\[\lambda_X=\frac{1}{4}\big
(c_1^2(V)-4c_2(V)\big)+\frac{g}{2}\on{deg}c_1(V)\big|_{B_0}.\,\,\,\qed\] 

\subsubsection{Global Calculation of $\kappa_X$.}
\label{globalkappa} Since $\omega_f=
\omega_X\otimes \omega_B^{-1}$,
\begin{equation}
\kappa_X=(K_X-\pi^*K_B)^2=K_X^2-8(g_{\scriptscriptstyle B}-1)(g-1).
\label{kappa}
\end{equation}
\noindent From Lemma~\ref{adjunction} we calculate
\begin{eqnarray*}
K_X^2&\!\!=\!\!\!&(K_{{\proj}V}+X)^2X=\big(\xi+\pi^*c_1(V)+\pi^*K_Y\big)^2
(3\xi+2\pi^*c_1(V))\\
     &\!\!=\!\!\!\!&2c_1^2(V)-3c_2(V)+4c_1(V)K_Y+3K_Y^2.
\end{eqnarray*}

\noindent
We calculate the contribution of $g_{\scriptscriptstyle B}$
in $K_X^2$: $8(g_{\scriptscriptstyle B}-1)\on{deg}c_1(V)|_{F_Y}+
24(1-g_{\scriptscriptstyle B})=8(g_{\scriptscriptstyle B}-1)(g-1)$, 
which is exactly the necessary adjustment for $\kappa_X$ 
in (\ref{kappa}). Therefore,
\begin{eqnarray*}
\kappa_X&\!\!\!=\!\!\!&2c_1^2(V)-3c_2(V)-8\on{deg}c_1(V)\big|_{B_0}\\
        &\!\!\!=\!\!\!&\frac{3}{4}\big(c_1^2(V)-4c_2(V)\big)+\frac{5}{2}
\on{deg}c_1(V)\big|_{B_0}\on{deg}c_1(V)\big|_{F_Y}-8\on{deg}c_1(V)\big|_{B_0}\\
        &\!\!\!=\!\!\!&\frac{3}{4}(c_1(V)^2-4c_2(V))+\frac{5g-6}{2}
\on{deg}c_1(V)\big|_{B_0}.\,\,\,\qed
\end{eqnarray*}

\subsection{Index theorem on the surface $X$.}
\label{indextheorem} Now that we have completed
the proof of Prop.~\ref{lambda_X,kappa_X}, we notice that any
bound on the ratio $\delta_X/\lambda_X$ would be equivalent to some
inequality involving the genus $g$ and the two invariants
discussed earlier: $\on{deg}c_1(V)|_{B_0}$ and the quantity
$c_1(V)^2\!-\!4c_2(V)$. This inequality
should be a fairly general one, since the only relevant information in our
situation is that $X$ is a triple cover of a ruled surface $Y$. 
One way of obtaining such general inequalities in ${\mathbb A}^2(X)\otimes 
\mathbb{Q}$
is via 

\begin{thm}[Hodge Index] Let $H$ be an ample
divisor on the smooth surface $X$, and let $\eta$ be a divisor on $X$,
numerically not equivalent to 0. If $\eta \cdot H=0$, then $\eta^2<0$.
\label{Index}
\end{thm}

The question here, of course, is how to find suitable divisors $H$ and
$\eta$ that would yield our result for the maximal slope bound. 
For that, we make use of the triple cover
$\phi:X\rightarrow Y$. If $H$ is any {\it ample} divisor on $Y$, then
its pullback $\phi^*H$ to $X$ is also ample. This follows from 

\begin{thm}[Nakai-Moishezon Criterion] A divisor $A$ on the smooth surface
$X$ is ample if and only if $A^2>0$ and $A\cdot C>0$ for all irreducible
curves $C$ in $X$.
\label{Nakai}
\end{thm}

\noindent
Since $H$ is ample itself, $(\phi^*H)^2=3H^2>0$ and $(\phi^*H)\!\cdot\!
C=H\!\cdot\!\phi_*(C)>0$
for any curve $C$ on $X$, so that $\phi^*H$ is also ample on $X$.
Now, if we find a divisor $\eta$ on $X$ such that $\eta\cdot\phi^*\on{Pic}
Y=0$, we will have assured that $\eta\cdot\phi^*H=0$, and then 
the Index theorem
will assert $\eta^2\leq 0$. As $X$ is a divisor itself on ${\proj}V$,
its Picard group naturally contains the restriction of $\on{Pic}{\proj}V$ to
$X$. We look for $\eta$ inside this subgroup, and for our purposes we may write
it in the form $\eta=\big(\zeta+\pi^*C_1\big)\big|_X$
for some $C_1\in \on{Pic}_{\mathbb Q}Y$. Let $C$ be any divisor class 
$\on{Pic}_{\mathbb {Q}}Y$. We compute
\[\eta\cdot \phi^*C=\big(\zeta+\pi^*C_1\big)\big(3\zeta+2\pi^*c_1(V)\big)
\pi^*C=C(3C_1-c_1(V)).\]
We want this to be zero for all $C$, so we naturally take
$C_1={\displaystyle\frac{1}{3}}c_1(V)\in\on{Pic}_{\mathbb {Q}}Y$. 
We summarize the above discussion in

\begin{lem}[Index Theorem on $X$] The divisor class 
$\eta=\big(\zeta+\frac{1}{3}\pi^*c_1(V)\big)\big|_X$ on $X$
satisfies $\eta\cdot \phi^*\on{Pic}(Y)=0$. In particular, for
an ample divisor $H$ on $Y$, the pullback $\phi^*H$ is also ample on $X$
and $\eta \cdot \phi^*H=0$. Consequently, $\eta^2\leq 0$ with equality
if and only if $\eta$ is numerically equivalent to $0$ on $X$.
\label{Eta}
\end{lem}
We have shown that 
\begin{equation}
0\geq3\eta^2=3\big(\zeta+\frac{1}{3}\pi^*c_1(V)\big)^2
\big(3\zeta+2\pi^*c_1(V)\big)=2c_1^2(V)-9c_2(V),
\label{eta}
\end{equation}
or equivalently,
\begin{equation}
2(g+2)\on{deg}c_1(V)\big|_{B_0}-9\big(c_1^2(V)-4c_2(V)\big)\geq 0.
\label{indexinequality}
\end{equation}

We are now ready to find a maximal bound for the slope of $X$. Recall the
formulas for $\lambda_X$ and $\kappa_X$ (cf.~Prop.~\ref{lambda_X,kappa_X}),
and write
\[\delta_X=12\lambda_X-\kappa_X=\displaystyle{
\frac{7g+6}{2}\on{deg}c_1(V)\big|_{B_0}
+\frac{9}{4}\big(c_1^2(V)-4c_2(V)\big)}.\]

\medskip
In view of the type of bound for the ratio $\delta_X/\lambda_X$, which
 we aim to achieve, we have to eliminate any extra terms and use inequality
(\ref{eta}). Our only choice is to subtract
\begin{eqnarray*}
36(g+1)\lambda_X-(5g+1)\delta_X&\!\!\!=\!\!\!&
\frac{1}{2}\big(36(g+1)g-(5g+1)(7g+6)\big)
\on{deg}c_1(V)\big|_{B_0}+\\
&& +\frac{1}{4}\big(9(g+1)-9(5g+1)\big)\big(c_1^2(V)-4c_2(V)\big)\\
&\!\!\!=\!\!\!&
\frac{1}{2}(g^2-g-6)\on{deg}c_1(V)\big|_{B_0}-\frac{9}{4}(g-3)
\big(c_1^2(V)-4c_2(V)\big)\\
&\!\!\!=\!\!\!&\frac{g-3}{4}\big[2(g+2)\on{deg}c_1(V)\big|_{B_0}-
9\big(c_1^2(V)-4c_2(V)\big)\big]\\
&\!\!\!=\!\!\!&(g-3)(9c_2(V)-2c_1^2(V))\geq 0
\end{eqnarray*}

\smallskip
As a result, we establish an exact maximal bound for the slopes of our
triple covers:

\begin{thm}[Main Theorem in Triple Cover Case] 
Given a triple cover \newline$\phi:X\!\rightarrow \!Y$ satisfying in
the Basic construction, the slope of $X$ satisfies
\[\frac{\delta_X}{\lambda_X}\leq \frac{36(g+1)}{5g+1}\cdot\]
Equality is achieved if and only if $g=3$, or $g>3$ and $\eta\equiv 0$ on $X$.
\label{maintheorem}
\end{thm}

\subsection{When is the maximal bound achieved?}
\label{whenmaximal} 
\subsubsection{The branch divisor of $\phi$}
From GRR, applied to $\phi:X\rightarrow Y$ and the sheaf $\cal{O}_X$, we
obtain a description of $c_1(V)$:
\[\on{ch}(\phi_{!}\cal{O}_X).\on{td}\cal{T}_Y=
\phi_*(\on{ch}\cal{O}_X.\on{td}\cal {T}_X),\]
\[\on{ch}(\phi_*\cal{O}_X)\big(1-\frac{1}{2}K_Y+\frac{1}{12}(K_Y^2+c_2(Y))\
\big)=\phi_*\big(1-\frac{1}{2}K_X+\frac{1}{12}(K_X^2+c_2(X)\big)\]
\[\Rightarrow c_1(\phi_*\cal{O}_X)=-\frac{1}{2}\big(\phi_*K_X-3K_Y).\]
For the {\it ramification} divisor $R$ on $X$ we know 
 $K_X=\phi^*K_Y+R$, so that
$\phi_*K_X=3K_Y+\phi_*R$. Hence $c_1(V)=-c_1(\phi_*\cal{O}_X)=
\frac{1}{2}\phi_* R$. In other words, from Lemma~\ref{D=2c_1(V)}
 we conclude that $c_1(V)$ is half of the {\it branch}
divisor $D$ on $Y$.
On the other hand, we can rewrite the condition $\eta\equiv 0$ 
in the following way:
\[0 \equiv 3\eta=\big(3\zeta+\pi^*c_1(V)\big)\big|_X=
\big(X-\pi^*c_1(V)\big)\big|_X=
c_1\big(\cal{O}_{{\proj}V}(X)\big|_X\big)-\pi^*c_1(V)\big|_X\]
\[\Leftrightarrow\,\,c_1\big(\cal{O}_{{\proj}V}(X)\big|_X\big)\equiv
\frac{1}{2}\phi^*D.\]
The self-intersection of $X$ on ${\proj}V$ satisfies (cf.~
\cite{self-intersection})
\[i^*i_*(1_X)=c_1(\cal{N}_{X/{\proj}V})\,\,\Rightarrow\,\,
X\cdot X=i_*c_1(\cal{N}_{X/{\proj}V}).\]
In particular, our condition $\eta\equiv 0$ can be expressed as
$\displaystyle
{c_1(\cal{N}_{X/{\proj}V})\equiv \displaystyle{\frac{1}{2}}\phi^*D}$.

\subsubsection{Examples of the maximal bound}
Constructing examples of families achieving the maximal bound is not
so easy, considering that the condition $\eta\equiv 0$ is not 
useful in practice. Instead, we start from the Basic construction and
attempt to find a ruled surface $Y$ and a rank 2 vector bundle $V$ on it
satisfying the equality in (\ref{indexinequality}), as well as
the ``genus condition'' given in Lemma~\ref{adjunction}. The former will
ensure the maximal ratio $\delta/\lambda = 36(g+1)/(5g+1)$, while the
latter will imply that the fibers of our family are indeed of genus $g$.
The remaining question is what linear series $3\zeta+\phi^*D$
has an irreducible member with at most rational double points
as singularities, which would serve as the total space of our family $X$.

\medskip
It is hard to work with the canonically defined bundle $V=
\phi_*(\cal{O}_X)/\cal{O}_Y$, since not every vector bundle $W$ of rank $2$
on $Y$ is of this form for some surface $X$. But any $W$ is
congruent to some $V$ after a twist by an appropriate line bundle $M$:
$V=W\otimes M$, and ${\proj}V\cong {\proj}W$.
 So, it seems reasonable to start with $W$ rather than $V$,
and use the invariant forms of our required equalities
(cf.~Sect.~\ref{globalcalculation}). This means 
replacing the degrees of $c_1(V)$ on $B_0$ and $F_Y$ by the corresponding
invariant 
degrees of $2D-3c_1(V)$. Thus, we need for some divisor $\widehat{D}$ on $Y$:
\begin{equation}
2(g+2)(2\on{deg}\widehat{D}\big|_{B_0}-3\on{deg}c_1(W)\big|_{B_0})
=9\big(c_1^2(W)-4c_2(W)\big),
\label{condition1}
\end{equation}
\begin{equation}
g+2=2\on{deg}\widehat{D}\big|_{F_Y}-3\on{deg}c_1(W)\big|_{F_Y}.
\label{condition2}
\end{equation}
For a general fiber $F_Y$ of $Y$ consider the rational ruled surface 
(cf.~Fig.~\ref{intersection}):
\[{\mathbf F}_e=\pi^*F_Y={\proj}(W|_{F_Y})=
{\proj}(\cal{O}_{{\proj}^1}\oplus
\cal{O}_{{\proj}^1}(e)),\,\,\,\on{with}\,\,
e\geq 0,\]
Let $S^{\prime}$ be the section in ${\mathbf F}_e$ with 
self-intersection $(S^{\prime})^2=-e$, and let $F_{\pi}$ be the fiber of
${\mathbf F}_e$ (in terms of the map $\pi:{\mathbf P}V\rightarrow Y$,
$F_{\pi}=\pi^*(\on{pt})$). Since a general
fiber $C$ of our family is embedded in ${\mathbf F}_e$,
the linear system \[|C|=\big|3S^{\prime}+\frac{g+2+3e}{2}F_{\pi}\big|\] 
has an irreducible nonsingular member. Equivalently, 
$C\cdot S^{\prime}\geq 0$, i.e. 
\begin{equation}e\leq (g+2)/3\,\,\,\on{and}\,\,\, 
e\equiv g(\on{mod}2)
\label{e-conditions}
\end{equation}
(compare with Lemma~\ref{gentrig}).
This forces three types of extremal examples 
according to the remainders $g(\on{mod}3)$.

\bigskip
\noindent{\bf Example 7.1 ($g\equiv 0(\on{mod}3)$).} Let $g=3e$ for 
some $e\in \mathbb{N}$. Set the base curve $B={\proj}^1$, and the ruled surface
\[Y={\proj}(\cal{O}_B\oplus \cal{O}_B(6))={\mathbf F}_6.\] 
Let $B^{\prime}$ be the section in $Y$ with smallest
self-intersection: $(B^{\prime})^2=-6$, thus 
$B_0=B^{\prime}+3F_Y$ with $B_0^2=0$. Let $Q=B^{\prime}+6F_Y$, and  
we choose two divisors $\widehat{D}$ and $E$ on $Y$ as follows:
\[\widehat{D}=(g+1)Q\,\,\,\on{and}\,\,\,E=eB^{\prime}+2(g+1)F_Y.\]
For the vector bundle $W$ on $Y$ we set $W=\cal{O}_Y\oplus
\cal{O}_Y(E)$ so that $c_1(W)=E$ and $c_2(W)=0$. 
We claim that the linear system $L=|3\zeta+\pi^*\widehat{D}|$ on
the 3-fold ${\proj}W$ contains an irreducible smooth member,
which we set to be our surface $X$ with maximal ratio $\delta/\lambda$.
\begin{figure}[h]
$$\psdraw{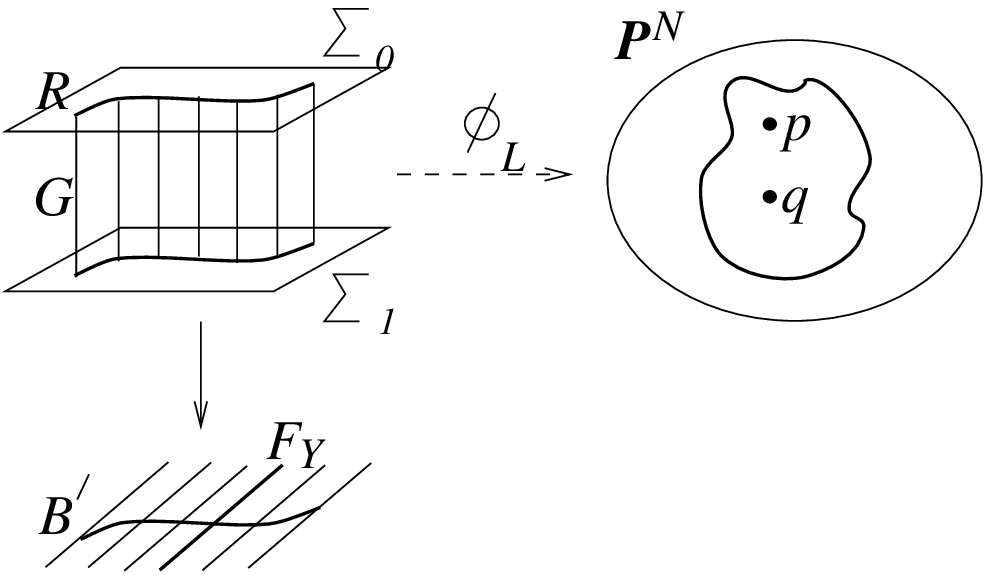}{2.6in}{1.2in}$$
\caption{Example with $g\equiv 0$(mod$3$)}
\label{maximal1.fig}
\end{figure}
Indeed, it is trivial to check conditions
(\ref{condition1}--20). Further, 
for {\it any} fiber $F_Y$ of $Y$: 
\[\pi^*F_Y={\proj}\big(\cal{O}_{{\proj}^1}\oplus \cal{O}_{{\proj}^1}
(E\cdot F_Y)\big)={\proj}\big(\cal{O}_{{\proj}^1}\oplus \cal{O}_{{\proj}^1}
(e)\big)={\mathbf F}_e,\]
so that $e=g/3$ satisfies the required conditions (\ref{e-conditions}).

\medskip
The only nontrivial fact is the existence of the wanted member $X$ in the 
linear system $L$ on ${\proj}W$. Consider two sections $\Sigma_0$
and $\Sigma_1$ of ${\proj}W$ corresponding to the subbundles
$\cal{O}_Y$ and $\cal{O}_Y(E)$ of $W$, respectively:
$\Sigma_0\in |\zeta|,\,\,\,\Sigma_1\in |\zeta+\pi^*(E)|,$ so that
$\Sigma_1\sim \Sigma_0+E$ (cf.~Fig.~\ref{maximal1.fig}).

Note that $\Sigma_0\cdot \Sigma_1=0$ and
$\Sigma_0\cdot L=\Sigma_0\cdot \pi^*B^{\prime}$. 
In other words, if $G=\pi^*B^{\prime}$ is the ruled surface
over $B^{\prime}$, then $\Sigma_0$ intersects every irreducible
member of $L$ in the curve $R=\Sigma_0\cap G$. On the other hand,
if a member of $L$ meets $\Sigma_0$ in a point outside $R$, then this
member contains entirely $\Sigma_0$. Thus, $L$ does not distinguish
the points on $\Sigma_0$, and $R$ is in the base locus of $L$.
Similarly, the restriction $L|_G=|3\Sigma_0|_G|=|3R|$
has exactly one section on $G$, namely, $3R$. Again it
follows that $L$ does not distinguish the points on $G$.

\smallskip
Away from the closed subset  $Z=\Sigma_0\cup G$, 
the linear system $L$ is in fact very
ample. This can be checked by showing directly that $L$ separates
points and tangent vectors on ${\mathbf P}W-Z$.
 Therefore, $L$ defines a
rational map \[\phi_L:{\proj}W \rightarrow {\mathbf P}(H^0(L)^{\widehat{
\phantom{o}}})={\proj}^N.\]
The map $\phi_L$ is regular on ${\proj}W-R$, embeds ${\proj}W-Z$, and 
contracts $\Sigma_0-R$ and $G-R$ to two distinct
points $p$ and $q$ in ${\proj}^N$. By Bertini's theorem
(cf.~\cite{Hartshorne}), the {\it general} member of $L$ is 
{\it smooth} away from the base locus $R$. 
Let $H$ be a general hyperplane in ${\proj}^N$
not passing through $p$ and $q$. Pulling $H$ back to ${\proj}W$
yields a member $X$
of $L$ not containing $\Sigma_0$ or $G$, and hence
irreducible. 

\medskip
It remains to show that the total space of $X$ has at most
finitely many double point singularities along the curve $R$.
Since the member $3\Sigma_1+G\in|L|$ is smooth along $R$, then
the general member of $|L|$ must be smooth along $R$. Hence our surface $X$
has, in fact, a smooth total space. This concludes the construction
of our maximal bound family of trigonal curves.

\bigskip
\noindent{\bf Example 7.2 ($g\equiv 1(\on{mod}3)$).} Set $g=3e-2$ for
$e\in{\mathbb N}$. Then $e$ satisfies the requirements of our construction:
$e=(g+2)/3$ and $e\equiv g (\on{mod}2)$. For the ruled surface $Y$ we
choose $Y=\proj^1\times\proj^1$. Let $E$ and $\widehat{D}$ be
the following divisors on $Y$:
$E=eB_0+fF_Y\,\,\on{and}\,\,\widehat{D}=3E,$
where $f\in{\mathbb N}$. The vector bundle $W$ on $Y$ is then defined by
$W=\cal{O}_Y\oplus\cal{O}_Y(E)$.
Finally, we indentify the total space of the surface $X$
with an irreducible smooth member of the linear system
$L=|3\zeta+\pi^*\widehat{D}|$ on the 3-fold $\proj W$.

\smallskip
The verification of this construction is similar to the previous example.
Here $L$ is very ample everywhere on $\proj W$ except on the section
$\Sigma_0$, which is contracted to a point under the map $\phi_{L}$.
This example, in somewhat different context, is shown in \cite{small}.

\medskip
\noindent{\bf Remark 7.1} The case of $g\equiv 2(\on{mod}3)$ is 
complicated by the fact that we cannot take $e=(g+1)/3$, for then
$e\not \equiv g(\on{mod}2)$. For example, if $g=5$, then the only
possibility is $e=1$. In the notation of Section~12,
all trigonal curves have lowest Maroni invariant of $1$, and there is
no Maroni locus. For now,
 in this case we have not been able to construct a trigonal
family with singular general member, whose ratio is $36(g+1)/(5g+1)$.

 \bigskip
\section*{\hspace*{1.9mm}8. 
Local Calculation of $\lambda,\delta$ and $\kappa$ in 
the General Case}

\setcounter{section}{8} 
\setcounter{subsection}{0} 
\setcounter{subsubsection}{0} 
\setcounter{lem}{0}
\setcounter{thm}{0}
\setcounter{prop}{0}
\setcounter{defn}{0} 
\setcounter{cor}{0} 
\setcounter{conj}{0} 
\setcounter{claim}{0} 
\setcounter{remark}{0}
\setcounter{equation}{0}
\label{generalcase}

\subsection{Notation and conventions}
\label{Conventions}
In this section we consider the general case of a trigonal
family $X\rightarrow B$. For convenience of notation, we shall assume
that the base curve $B$ intersects transversally and in general points
the boundary divisors of $\overline{\mathfrak{T}}_g$
(cf.~Fig.~\ref{Delta-k,i}). We will call such a base
curve {\it general}, and use this definition throughout
Section~8-10. Since
we work in the rational Picard group of $\overline{\mathfrak{M}}_g$, all arguments
and statements in the remaining cases are shown similarly in 
Sect.~11. From Prop.~6.1, 
we may assume that modulo a base change, our family $X\!\rightarrow\!
B$ fits in the following commutative diagram:

\setlength{\unitlength}{10mm}
\begin{figure}[h]
\begin{picture}(3,4.3)(-1,2)
\put(0,3.9){$\widetilde{X}\,\stackrel{\phi}{\longrightarrow} 
\widetilde{Y}$}
\put(0,5.1){$\widehat{X}\,\stackrel{\widehat{\phi}}{\longrightarrow}
 \widehat{Y}$}
\multiput(0.2,5)(0,-1.2){3}{\vector(0,-1){0.6}}
\put(1.45,3.8){\vector(-1,-2){1}}
\put(0,2.8){$X$}
\put(0,1.6){$B$}
\put(1.55,5){\vector(0,-1){0.6}}
\put(1.55,6.2){\vector(0,-1){0.6}}
\put(1.65,5.8){$\pi$}
\put(1.3,6.3){${\proj}V$}
\put(0.2,5.6){\vector(2,1){1.2}}
\end{picture}
\caption{General base $B$}
\label{general B}\vspace*{-4mm}
\end{figure}

\subsubsection{Relations in ${\rm Pic}_{\mathbb{Q}}\widehat{Y}$ and
${\rm Pic}_{\mathbb{Q}}{\proj}V$}
\label{relations}
The special fibers of of $\widehat{X}$ and of the birationally ruled 
surface $\widehat{Y}$ over $B$ are described in
Fig.~\ref{coef1.fig}--\ref{coef3.fig}. Since each such fiber in $\widehat{Y}$
is a {\it chain} $T$ of rational components, we can fix one of the end 
components to be the {\it root} $R$. We keep the
notation $E^-$ ($E^+$, respectively) for the ancestor (descendants,
respectively)
of a component $E$ in $T$. We also fix a general fiber $F_{\widehat{Y}}\cong
{{\proj}}^1$ of $\widehat{Y}$, and a section $B_{\widehat{Y}}$, which is
the pullback
of the corresponding section $B_0$ in $\widetilde{Y}$ (cf.~(\ref{normalize})). 
The rational Picard group of $\widehat{Y}$ is generated by $F_{\widehat{Y}}$, 
$B_{\widehat{Y}}$ and all non-root components $E$ of the special
fibers of $\widehat{Y}$:
\[\on{Pic}_{\mathbb{Q}}\widehat{Y}=\mathbb{Q}
B_{\widehat{Y}}\bigoplus \mathbb{Q}F_{\widehat{Y}}\!\!\!\bigoplus
_{E-\on{not}\,\on{root}}\!\!\!\mathbb{Q}E.\]
The intersection numbers of these generators are as follows: 
$B_{\widehat{Y}}^2=0=F_{\widehat{Y}}^2,\,\,\,
B_{\widehat{Y}}\cdot F_{\widehat{Y}}=1,\,\,\on{and}\,\,
E\cdot B_{\widehat{Y}}=E\cdot F_{\widehat{Y}}=0.$

\vspace*{-3mm}
\setlength{\unitlength}{7mm}
\begin{figure}[h]
\begin{picture}(3,2.3)(12.6,4.5)
\put(6,5.3){\line(1,0){2}}
\multiput(5,4.6)(0,-1.7){2}{\line(2,1){1.6}}
\multiput(6,3.6)(0,-0.1){2}{\line(1,0){2}}
\put(9.1,4.2){\line(-2,-1){1.6}}
\multiput(4.9,5)(0,-1.8){2}{$E^-$}
\multiput(6.8,5.45)(0,-1.7){2}{$E$}
\put(8.7,4.3){$E^+$}
\end{picture}
\vspace*{9mm}
\caption{$m_{E}$\hspace*{100mm}}
\label{m E}

\begin{picture}(7,4)(3,-2.4)
\put(4.45,2.9){\line(2,1){1.6}}
\put(5.5,3.6){\line(1,0){3}}
\put(6.5,3.4){\line(-1,1){1.2}}
\put(7.5,3.4){\line(1,1){1.2}}
\put(4.5,3.2){$E^-$}
\put(6.8,2.9){$E$}
\put(5.1,4.75){$E^+$}
\put(5.9,4.85){$E^+$}
\put(7.6,4.85){$E^+$}
\put(8.4,4.75){$E^+$}
\put(6.75,3.4){\line(-1,2){0.65}}
\put(7.25,3.4){\line(1,2){0.65}}
\end{picture}\vspace*{-38mm}
\caption{$E^2$}
\label{E^2}

\begin{picture}(3,4)(-1,-1.9)
\multiput(6,5.3)(0,-0.1){2}{\line(1,0){2}}
\multiput(5,4.6)(0,-1.65){2}{\line(2,1){1.6}}
\put(6,3.6){\line(1,0){2}}
\multiput(9.1,4.2)(0,1.8){2}{\line(-2,-1){1.6}}
\put(9.1,4.3){\line(-2,-1){1.6}}
\multiput(4.8,5)(0,-1.7){2}{$R$}
\multiput(6.8,5.45)(0,-1.7){2}{$E^-$}
\multiput(8.7,4.3)(0,1.8){2}{$E$}
\end{picture}
\vspace*{-37.5mm}
\caption{$\theta_E$\hspace*{-90mm}}
\label{theta_E}
\end{figure}
We also set $m_{\!\stackrel{\phantom{.}}{E}}=E\cdot E^-$ (cf.~Fig.~\ref{m E}):
\[m_{\!\stackrel{\phantom{.}}{E}}=\left\{\begin{array}{l}
             0\,\,\on{if}\,\,E=R\,\,\on{root},\\
             1\,\,\on{if}\,\,E\,\,\on{and}\,\,E^-\,\,\on{reduced},\\
             2\,\,\on{if}\,\,E\,\,\on{or}\,\,E^-\,\,\on{nonreduced}.
\end{array}\right.\]
 In this notation, due to the fact that $E\cdot T=E\cdot
F_{\widehat{Y}}=0$, the self-intersection of any $E$ is computed by
(cf.~Fig.~\ref{E^2}): 
\[E^2=-\sum_{\!\stackrel{\scriptstyle{E^{\prime}\not= E}}
{E^{\prime}\cap E\not= \emptyset}} E\cdot
E^{\prime}=-\sum_{\!\stackrel{\scriptstyle{E^{\prime}=E^+}}
{\on{or}\,\,E^{\prime}=E}}m(E^{\prime})\]

\smallskip
 In order to express the dualizing sheaf $K_{\widehat{Y}}$
in terms of the above generators of $\on{Pic}_{\mathbb{Q}}\widehat{Y}$,
for each component $E$ in $\widehat{Y}$ we denote by $\theta_E$ the length
of the path $\stackrel{\longrightarrow}
{RE}$, omitting any nonreduced components except for
$E$ itself. For example, in the two cases in Fig.~\ref{theta_E}
we have $\theta_E=1$ and $\theta_E=2$. Note that $\theta_R=0$.

Considering the ``effective'' blow-ups on $\widetilde{Y}$,
necessary to construct $\widehat{Y}$, we immediately obtain the following
identities (compare with (\ref{K_PV}) and (\ref{K_Yglobal})).
\begin{lem} In $\on{Pic}_{\mathbb{Q}}\widehat{Y}$
and $\on{Pic}_{\mathbb{Q}}{\proj}V$:
\begin{eqnarray*}
\on{(a)}&\!\!\!\!&\!\!
\displaystyle{K_{\widehat{Y}}
\equiv -2B_{\widehat{Y}}+(2g_B-2)F_{\widehat{Y}}+\sum_E \theta_EE},\\
\vspace*{-2mm}
\on{(b)}&\!\!\!\!&\!\!K_{{\proj}V}\equiv -2\zeta -\pi^*
c_1(V)+\pi^*K_{\widehat{Y}},\\
\on{(c)}&\!\!\!\!&\!\!K_{{\proj}V/\!_{\scriptstyle{B}}}
\equiv -2\zeta -\pi^*c_1(V)-2\pi^*B_{\widehat{Y}}+\sum_E \theta_E\pi^*E.
\end{eqnarray*}
\label{Kdivisors}\vspace*{-8mm}
\end{lem}

The hyperplane section $\zeta$ of ${\proj}V$ and
the rank 2 vector bundle $V$ on $\widehat{Y}$ are defined similarly
as in Section~7. 
Thus, in $\on{Pic}_{\mathbb{Q}}{\proj}V$ we have $\widehat{X}\sim 3\zeta +\pi^*D$
for a certain divisor $D$ on $\widehat{Y}$. By analogy with
Lemma~\ref{D=2c_1(V)}, one shows that $D\equiv 2c_1(V)$ in
$\on{Pic}_{\mathbb{Q}}Y$, so that
\begin{equation}
\widehat{X}\equiv 3\zeta+ 2\pi^*c_1(V).
\label{genX}
\end{equation}
Using the above notation for $\on{Pic}_{\mathbb{Q}}\widehat{Y}$ we can write
for some half-integers $c,d,\gamma_{\!\stackrel{\phantom{.}}{E}}$:
\begin{equation}
\displaystyle{c_1(V)\equiv cB_{\widehat{Y}}+dF_{\widehat{Y}}+
\sum_{E}\gamma_{\!\stackrel{\phantom{.}}{E}}E}.
\label{genc_1(V)}
\end{equation}
Here we can assume that $\gamma_{\!\stackrel{\phantom{.}}{R}}=0$
by replacing $R$ with a linear combination of the remaining components $E$
in its chain $T$ (compare with (\ref{define D,c1(V)})).

Finally, we need the top Chern classes of $\widehat{Y}$ and ${\proj}V$
in terms of intersections of known divisors and other known invariants
of the two surfaces (compare with (\ref{c_2(PV)}) and (\ref{c_2(Y)global})).
\begin{lem} In the Chow rings $\mathbb{A}(\widehat{Y})$ and
$\mathbb{A}({\proj}V)$ the following equalities are true:
\begin{eqnarray*}
\on{(a)}&\!\!&\!\!\!\!c_2(\widehat{Y})=c_2(Y)+\sum_{E\not =R} 1=4(1-g_B)+\sum_
{E\not = R} 1, \\
\on{(b)}&\!\!&\!\!\!\!c_2({\proj}V)=
c_2(\widehat{Y})-\pi^*K_{\widehat{Y}}(2\zeta+\pi^*c_1(V)),\\
\on{(c)}&\!\!&\!\!\!\!\displaystyle{
c_2({\proj}V/\!_{\displaystyle{B}})=
-\pi^*K_{\widehat{Y}/B}(2\zeta+\pi^*c_1(V))+\sum_{E\not =R} 1}.
\end{eqnarray*}
\end{lem}
\label{conventions}
\subsubsection{A technical lemma}
\label{technicallemma} In the sequel,
we will work with several functions defined on the set of components 
$\{E\}$ in $\widehat{Y}$. For easier calculations, to
any such function $f$ we associate the {\it difference function} $F$
by setting $F_{\!\stackrel{\phantom{.}}{E}}:=
f_{\!\stackrel{\phantom{.}}{E}}-f_{\!\stackrel{\phantom{.}}
{E^-}}$ for all $E$. Since $R^-$ does not exist, we define
$f_{R^-}=0$ for all roots $R$ in $\widehat{Y}$.
\begin{lem}For any functions $f$ and $h$ defined on the
set of components $\{E\}$ in $\widehat{Y}$, the following identity 
holds true:
\[\sum_E f_{\!\stackrel{\phantom{.}}{E}}E \cdot \sum_E
h_{\!\stackrel{\phantom{.}}{E}}E=-\sum_E
(m\!\cdot\! F\!\cdot\! H)_{\!\stackrel{\phantom{.}}{E}}.\]
\label{technical}\vspace*{-7mm}
\end{lem}
\noindent{\it Proof.} We rewrite the lefthand side as
$\displaystyle{\sum_{E_1\not= E_2}f_{\!\stackrel{\phantom{.}}{E_1}}
h_{\!\stackrel{\phantom{.}}{E_2}}E_1E_2+\sum_E f_{\!\stackrel{\phantom{.}}{E}}h_{\!\stackrel{\phantom{.}}{E}}E^2=}$
\[=\sum_{E}\left(f_{\!\stackrel{\phantom{.}}{E^-}}
h_{\!\stackrel{\phantom{.}}{E}}+f_{\!\stackrel{\phantom{.}}{E}}
h_{\!\stackrel{\phantom{.}}{E^-}}\right)m_{\!\stackrel{\phantom{.}}{E}}-
\sum_E
\left(f_{\!\stackrel{\phantom{.}}{E}}h_{\!\stackrel{\phantom{.}}{E}}+
f_{\!\stackrel{\phantom{.}}{E^-}}
h_{\!\stackrel{\phantom{.}}{E^-}}\right)m_{\!\stackrel{\phantom{.}}{E}}=\]
\[=\sum_E\left(f_{\!\stackrel{\phantom{.}}{E^-}}-
f_{\!\stackrel{\phantom{.}}{E}}\right)\left(h_{\!\stackrel{\phantom{.}}{E}}
-h_{\!\stackrel{\phantom{.}}{E^-}}\right)m_{\!\stackrel{\phantom{.}}{E}}
=\sum_E(m\!\cdot\! F\!\cdot\! H)
_{\!\stackrel{\phantom{.}}{E}}.\,\,\,\qed\]

We have noted that all three functions $m,\theta$ and $\gamma$ are zero on the
roots $R$ in $\widehat{Y}$. Since we shall be working specifically with
these three functions, it makes sense to restrict from now on all sums
$\sum_E$ only to the non-roots $E$ in $\widehat{Y}$. 
With this in mind, in every application
of Lemma~\ref{technical} one must check that
the corresponding functions $f$ and $h$ have the same property:
$f_R=0=h_R$, so that we can restrict the sums 
in Lemma~\ref{technical} also to all {\it non-roots}
$E$ in $\widehat{Y}$. In fact, in all cases this verification
will be obvious as $f$ and
$h$ will  be, for the most part, linear combinations of $\theta$ and $\gamma$.

\smallskip
\noindent{\bf Example 8.1.} From expression (\ref{genc_1(V)}) for $c_1(V)$
as a divisor on $\widehat{Y}$, and Lemma~\ref{technical}:  
\begin{equation} c_1^2(V)=2cd+
\sum_{E}\gamma_{\!\stackrel{\phantom{.}}{E}}E\cdot
\sum_{E}\gamma_{\!\stackrel{\phantom{.}}{E}}E=
2cd-\sum_{E}m_{\!\stackrel{\phantom{.}}{E}}
\Gamma^2_{\!\stackrel{\phantom{.}}{E}}.
\label{c^2_1(V)}
\end{equation}

\subsection{Computation of the invariants $\lambda_{\widehat{X}},
\kappa_{\widehat{X}}$ and $\delta$}
\label{computation}
 The following proposition~\ref{hatlambda} is a generalization
of the corresponding statement in Section~7 (cf.~
Prop.~\ref{lambda_X,kappa_X}). We set $\Gamma_{
\!\stackrel{\phantom{.}}{E}}=\gamma_{\!\stackrel
{\phantom{.}}{E}}-\gamma_{\!\stackrel{\phantom{.}}{E^-}}$ and
$\Theta_{\!\stackrel{\phantom{.}}{E}}
=\theta_{\!\stackrel{\phantom{.}}{E}}
-\theta_{\!\stackrel{\phantom{.}}{E^-}}$ to be the difference functions of
$\gamma$ and $\theta$.

\begin{prop}
The degrees of the invariants $\lambda_{\widehat{X}}$ and 
$\kappa_{\widehat{X}}$
on $\widehat{X}$ are given by
\[\lambda_{\widehat{X}}=d(g+1)-c_2(V)-\frac{1}{4}\sum_E
\left\{m_{\!\stackrel{\phantom{.}}{E}}\cdot\left(2\Gamma^2+2\Gamma\cdot \Theta
+\Theta^2\right)_{\!\stackrel{\phantom{.}}{E}}-1\right\},\]
\vspace*{-6mm}\[\,\,\,\,\kappa_{\widehat{X}}=4dg-3c_2(V)-\sum_Em_{\!\stackrel{\phantom{.}}{E}}\left
(2\Gamma^2+4\Gamma\Theta+3\Theta^2\right)_{\!\stackrel{\phantom{.}}{E}}.\]
\label{hatlambda}\vspace*{-6mm}
\end{prop}

\begin{proof} 
One starts with the Euler characteristic formula
$\lambda_{\widehat{X}}=\chi(\cal{O}_{\widehat{X}})-
\chi(\cal{O}_C)\cdot \chi(\cal{O}_B),$
or the adjunction formula
$\kappa_{\widehat{X}}=\big(\widehat{X}+K_{{\proj}V/\!_{\scriptstyle{B}}})
^2\widehat{X}.$
The rest of the proof is a straight forward calculation, which uses the
equalities given in~\ref{conventions}, and 
is substantially simplified by Lemma~\ref{technical}. \end{proof}

\begin{cor} The degree $\delta$ on the original family $X$ is given by
\[\delta=4d(2g+3)-9c_2(V)-\!\sum_T\mu(T)-\!\sum_{\on{ram}1} 1-\!
\sum_{\on{ram}2}
 3-\!\sum_E\left\{m_{\!\stackrel{\phantom{.}}{E}}
\left(4\Gamma^2+2\Gamma\Theta\right)_{\!\stackrel{\phantom{.}}{E}}-3\right\}\]
\label{hatdelta}\vspace*{-4mm}
\end{cor}
Here $\mu(T)$ stands for the quasi-admissible contribution to
$\kappa_{\widehat{X}}$ of the preimage
$C=\widehat{\phi}^*T$ in $\widehat{X}$, as defined in Lemma~\ref{mu(C)}.

\smallskip
\begin{proof} Since $\lambda=\lambda_{\widehat{X}}$, 
$\kappa=\kappa_{\widehat{X}}+\sum_T \mu(T)+
\sum_{\on{ram}1} 1+\sum_{\on{ram}2}3$, 
and $\delta=12\lambda-\kappa$, the statement immediately
follows from Prop.~\ref{hatlambda}.\end{proof}

\subsection{The arithmetic genus $p_{{E}}$,
 and the invariants $\Gamma_{{E^{\prime}}}$ and $\Theta_{{E^{\prime}}}$}
\label{arithmetic}
\begin{figure}[h]
$$\psdraw{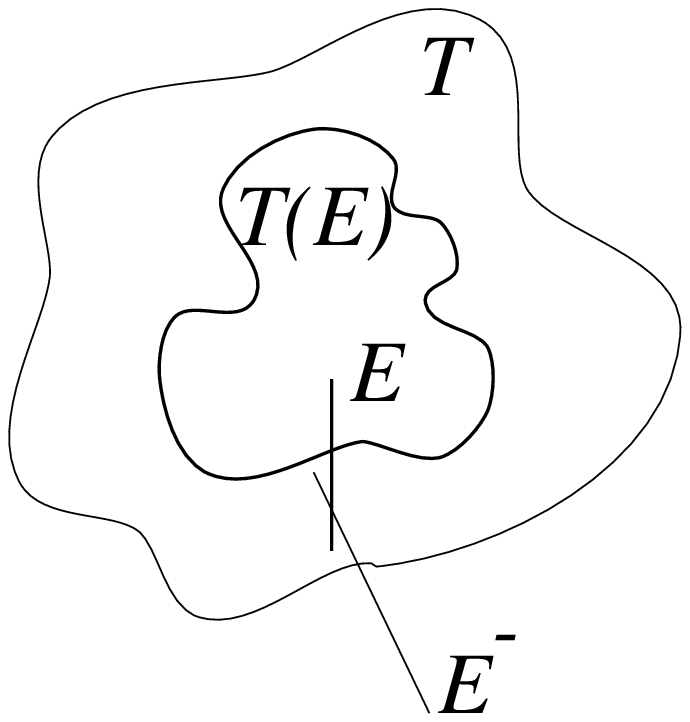}{1.2in}{1.2in}$$
\caption{Subtree $T(E)$}
\label{subtree}
\end{figure}

For a component $E$ in a special fiber $T$ of
$\widehat{Y}$, we define $T(E)$ to be the subtree of $T$ generated 
by the component $E$. In other
words, $T(E)$ is the union of all components $E^{\prime}\in T$ such that
$E^{\prime}\geq E$ (cf. Fig.~\ref{subtree}). For simplicity, we set
$p_{\!\stackrel{\phantom{.}}{E}}:=p_a\big(\phi^*(T(E))\big)$ 
to be the arithmetic genus of the inverse image $\phi^*(T(E))$ in
$\widehat{X}$. It can be easily computed via the following 
 analog of Lemma~\ref{adjunction},
where $T$ consisted of a single component $E=R$.
\begin{lem} For a general base curve $B$ and for any non-root component
$E\in T$:
\begin{equation}
\displaystyle{p_{\!\stackrel{\phantom{.}}{E}}=-m_{\!\stackrel
{\phantom{.}}{E}}\left(\Gamma_E+\frac{3(\Theta_E+1)}{2}\right)+1}.
\label{arithmetic equation}
\end{equation}
\label{arithgenus}\vspace*{-6mm}
\end{lem}
\noindent{\it Proof.} From the adjunction formula for the divisor
$\phi^*(T(E))$ in $\widehat{X}$:
\[2p_{\!\stackrel{\phantom{.}}{E}}-2=\left(K_{\widehat{X}}+\phi^*(T(E))\right)\phi^*(T(E))=
\left((K_{{\proj}V}+\widehat{X})|_{\widehat{X}}+\sum_{E^{\prime}}
\delta_{E^{\prime}}\widehat{\phi}^*E^{\prime}\right)\sum_{E^{\prime}}
\delta_{E^{\prime}}\widehat{\phi}^*E^{\prime}.\]
Here $\delta_{E^{\prime}}=0$ if $E^{\prime}<E$, and $\delta_{E^{\prime}}=1$
otherwise. Thus, the sums above are effectively taken 
over all $E^{\prime}\geq E$. Substituting the expressions
for $K_{{\proj}V}$ and $\widehat{X}$ as divisors in ${{\proj}V}$
from Lemma~\ref{Kdivisors} and (\ref{genX}), we arrive at
\[2p_{\!\stackrel{\phantom{.}}{E}}-2=\sum_{E^{\prime}}\left(2
\gamma_{\!\stackrel{\phantom{.}}{E^{\prime}}}+
3\theta_{\!\stackrel{\phantom{.}}{E^{\prime}}}+
3\delta_{\!\stackrel{\phantom{.}}{E^{\prime}}}
\right)E^{\prime} \sum_{E^{\prime}}
\delta_{\!\stackrel{\phantom{.}}{E^{\prime}}}E^{\prime}.\]
Set $\Delta_{\!\stackrel{\phantom{.}}{E}}
=\delta_{\!\stackrel{\phantom{.}}{E}}
-\delta_{\!\stackrel{\phantom{.}}{E^-}}$, i.e.
$\Delta_{\!\stackrel{\phantom{.}}{E^{\prime}}}=1$ 
only if $E^{\prime}=E$; otherwise,
$\Delta_{\!\stackrel{\phantom{.}}{E^{\prime}}}=0$. 
By Lemma~\ref{technical}, 
\[2p_{\!\stackrel{\phantom{.}}{E}}-2=-\sum_{E^{\prime}} m_
{\!\stackrel{\phantom{.}}{E^{\prime}}}
\left(2\Gamma_{\!\stackrel{\phantom{.}}{E^{\prime}}}+
3\Theta_{\!\stackrel{\phantom{.}}{E^{\prime}}}+3
\Delta_{\!\stackrel{\phantom{.}}{E^{\prime}}}\right)\Delta_{\!\stackrel{\phantom{.}}{E^{\prime}}}\,\,
\Rightarrow\,\,2p_{\!\stackrel{\phantom{.}}{E}}-2=
-m_{\!\stackrel{\phantom{.}}{E}}\left(2\Gamma_{\!\stackrel{\phantom{.}}{E}}
+3\Theta_{\!\stackrel{\phantom{.}}{E}}+3\right).\,\,\,\qed\]

\smallskip
Now we can easily compute the invariants 
$m_{\!\stackrel{\phantom{.}}{E}}$, 
$\Theta_{\!\stackrel{\phantom{.}}{E^{\prime}}}$ and
$\Gamma_{\!\stackrel{\phantom{.}}{E^{\prime}}}$, appearing in the formulas
for $\lambda_{X}$ and $\kappa_{X}$.

\begin{cor} There are three possibilities for  the triple 
$(m_{\!\stackrel{\phantom{.}}{E}}, 
\Theta_{\!\stackrel{\phantom{.}}{E^{\prime}}},
\Gamma_{\!\stackrel{\phantom{.}}{E^{\prime}}})$, depending on
whether the components $E$ and $E^-$ of $T$ are reduced:
\begin{eqnarray*}
\on{(a)}\,\on{if}\,E,E^-\,\on{reduced},\,\on{then}&\!\!\!&\!\!\!\!
m_{\!\stackrel{\phantom{.}}{E}}=1,\,\,\Theta_{\!\stackrel{\phantom{.}}{E}}
=1,\,\,\Gamma_{\!\stackrel{\phantom{.}}{E}}=
-(p_{\!\stackrel{\phantom{.}}{E}}+2).\\
\on{(b)}\, \on{if}\,E\,\,\on{nonreduced}, \,\on{then}&\!\!\!&\!\!\!\!
m_{\!\stackrel{\phantom{.}}{E}}=2,\,\,\Theta_{\!\stackrel{\phantom{.}}{E}}
=1,\,\,\Gamma_{\!\stackrel{\phantom{.}}{E}}=
-({p_{\!\stackrel{\phantom{.}}{E}}+5})/{2}.\\
\on{(c)}\,\on{if}\,E^-\!\on{nonreduced}, \,\on{then}&\!\!\!&\!\!\!\!
m_{\!\stackrel{\phantom{.}}{E}}=2,\,\,\Theta_{\!\stackrel{\phantom{.}}{E}}
=0,\,\,\Gamma_{\!\stackrel{\phantom{.}}{E}}=
-({p_{\!\stackrel{\phantom{.}}{E}}+2})/{2}.
\end{eqnarray*}
\label{constants}\vspace*{-5mm}
\end{cor}
\begin{proof} Note that for the list all possible
special fibers $T$ of $\widehat{Y}$, each component $E$ fits in exactly one
of the three cases above (cf.~Fig.~\ref{coef1.fig}--\ref{coef3.fig}). 
The proof of the statement is 
immediate from the definitions of $m_{\!\stackrel{\phantom{.}}{E}}$ and
$\Theta_{\!\stackrel{\phantom{.}}{E^{\prime}}}$, and
from Lemma~\ref{arithgenus}. \end{proof}

\bigskip\section*{\hspace*{1.9mm}9. 
The Bogomolov Condition $4c_2-c_1^2$ and the $7+6/g$ Bound in
$\overline{\mathfrak{T}}_g$}

\setcounter{section}{9}
\setcounter{subsection}{0} 
\setcounter{subsubsection}{0} 
\setcounter{lem}{0}
\setcounter{thm}{0}
\setcounter{prop}{0}
\setcounter{defn}{0} 
\setcounter{cor}{0} 
\setcounter{conj}{0} 
\setcounter{claim}{0} 
\setcounter{remark}{0}
\setcounter{equation}{0}
\label{Bogomolov1}

With the conventions of Section~8, we state
the main proposition of the section.
\begin{prop} There exists an effective $\mathbb Q$-linear combination
$\cal{E}$ of boundary divisors $\Delta{\mathfrak{T}}_{k,i}$, not containing 
$\Delta{\mathfrak{T}}_0$, such that 
for a general base curve $B$ in $\overline{\mathfrak{T}}_g$:
\[(7g+6)\lambda|_B=g\delta|_B+\cal{E}|_B+\frac{g-3}{2}\left(4c_2(V)-c_1^2(V)
\right).\]
\label{bogomolov1}\vspace*{-10mm}
\end{prop}

 For a shorthand notation, we denote by $\mathfrak{S}$
the difference
\[\mathfrak{S}:=(7g+6)\lambda|_B-g\delta|_B -\frac{g-3}{2}\left(
4c_2(V)-c_1^2(V)\right).\]
Using the results of the previous section, we can write:
\begin{eqnarray*}
\mathfrak{S}&=&
-\frac{1}{4}\sum_E\left\{m_{\!\stackrel{\phantom{.}}{E}}\left(6\Gamma^2+
6(g+2)\Gamma\Theta+(7g+6)\Theta^2\right)_{\!\stackrel{\phantom{.}}{E}}+5g-6
\right\}\\
&&+\sum_T g\mu(T)+\sum_{\on{ram}1} g +\sum_{\on{ram}3} 3g.
\end{eqnarray*}
We defer the
proof of Prop.~\ref{bogomolov1} until the end of this section, when
all of the terms in this sum will be computed.

\subsection{Grouping the contributions of each $\Delta{\mathfrak{T}}_{k,i}$ in 
$\mathfrak{S}$}
\label{Grouping}
Substituting the results of Corollary~\ref{constants}
in the expression for $\mathfrak{S}$, we eliminate 
$m_{\!\stackrel{\phantom{.}}{E}}$, 
$\Theta_{\!\stackrel{\phantom{.}}{E^{\prime}}}$ and
$\Gamma_{\!\stackrel{\phantom{.}}{E^{\prime}}}$:
\begin{eqnarray*}
\mathfrak{S}&\!\!=\!\!&\sum_T g\mu(T)+\sum_{\on{ram}1} g+\sum_{\on{ram}2} 3g+
\frac{1}{4}\sum_{E,E^-\on{red}}\!\!\left(6(2+p_{\!\stackrel{\phantom{.}}{E}})
(g-p_{\!\stackrel{\phantom{.}}{E}})-12g\right)\\
&\!\!-\!\!
&\frac{1}{4}\sum_{E^-\on{nonred}}\!\!\!\!\!\left(3(p_{\!\stackrel{\phantom{.}}{E}}+2)^2+5g-6\right)+\frac{1}{4}\sum_{E\on{nonred}}\!\!\left(3(p_{\!\stackrel{\phantom{.}}{E}}+5)(2g-p_{\!\stackrel{\phantom{.}}{E}}-1)-19g+6\right).
\end{eqnarray*}
For each chain $T$ in $\widehat{Y}$,
the inverse image $\widehat{\phi}^*(T)$ in $\widehat{X}$ is a member 
(or a blow-up of a member) of exactly one boundary divisor
$\Delta{\mathfrak{T}}_{k,i}$. Consequently,
to find the contribution to $\mathfrak{S}$ of a specific $\Delta{\mathfrak{T}}_{k,i}$,
we calculate the sum in $\mathfrak{S}$ corresponding to all types of 
special fibers $\widehat{\phi}^*(T)$.

\begin{figure}[h]
$$\hspace*{-7mm}\psdraw{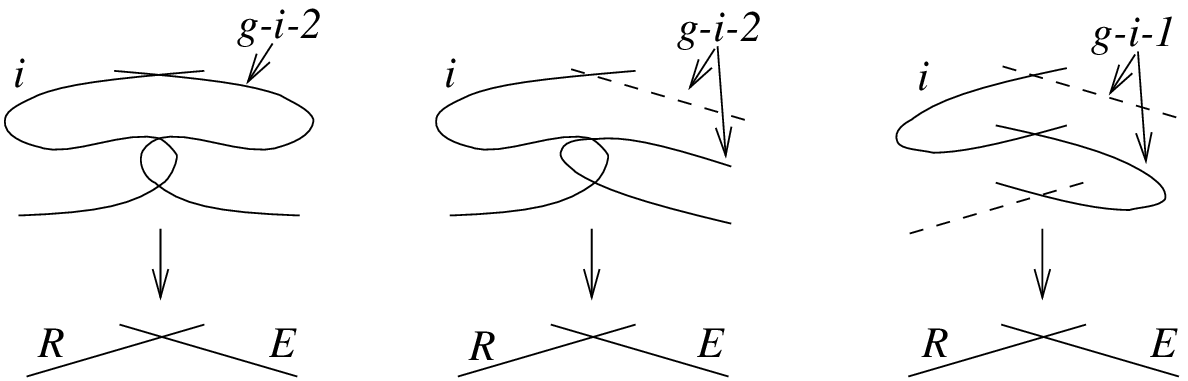}{4.5in}{1.3in}$$
\caption{Coefficients with no ramification}
\label{coef1.fig}
\end{figure}

\subsubsection{Contributions of $\Delta{\mathfrak{T}}_{1,i},
\Delta{\mathfrak{T}}_{2,i}$ and $\Delta{\mathfrak{T}}_{3,i}$} 
\label{contribution1}
Fig.~\ref{coef1.fig} presents the
special fibers corresponding to the boundary divisors 
$\Delta{\mathfrak{T}}_{1,i},\,\,\Delta{\mathfrak{T}}_{2,i},\,\,\Delta{\mathfrak{T}}_{3,i}$. In each of these cases,
there is only one component $E$ in $T$ besides the root
$R=E^-$. Thus, the subchain $T(E)$ in $T$ is trivial -- it consists
only of $E$. Its inverse image $\widehat{\phi}^*(E)$ is connected
for $\Delta{\mathfrak{T}}_{1,i}$, and consists of two connected curves for
$\Delta{\mathfrak{T}}_{2,i}$ and $\Delta{\mathfrak{T}}_{3,i}$. Setting the genus of the inverse
image of $R$ to be $i$, it is easy to see that the genus $p_{\!\stackrel{\phantom{.}}{E}}$ of
$\phi^*(E)$ is $g-i-2$ in the first two cases, and $g-i-1$ in the
third case.(The total genus of the original fiber of $X$,
drawn in full lines, must be $g$.)
Finally, counting the number of ``quasi-admissible''
blow-ups (drawn by  dashed lines), we see that $\mu(T)=0$ for $\Delta{\mathfrak{T}}_{1,i}$,
$\mu(T)=1$ for $\Delta{\mathfrak{T}}_{2,i}$, and $\mu(T)=2$ for
$\Delta{\mathfrak{T}}_{3,i}$ (cf.~Lemma~\ref{mu(C)}).
Note that there are no ramification modifications.

The contribution of each such fiber $\widehat{\phi}^*T$
to the sum $\mathfrak{S}$ is only
one summand of the first type ($E,E^-$reduced), plus the
quasi-admissible adjustment $g\mu(T)$. If $\widehat{\phi}^*T$
corresponds to the boundary divisor $\Delta{\mathfrak{T}}_{k,i}$, we denote this
contribution by $c_{k,i}$. In conclusion,
\[c_{k,i}=
\frac{1}{4}\big(6(2+p_{\!\stackrel{\phantom{.}}{E}})(g-p_{\!\stackrel
{\phantom{.}}{E}})-12g\big)+g\mu(T)\,\,\Rightarrow\,\,
c_{k,i}=\frac{3}{2}(i+2)(g-i)-(4-k)g,\,\,k=1,2,3.\]

\subsubsection{Contributions of $\Delta{\mathfrak{T}}_{4,i}$ and 
$\Delta{\mathfrak{T}}_{5,i}$: ramification index 1}
\label{contribution2}
\begin{figure}[h]
$$\psdraw{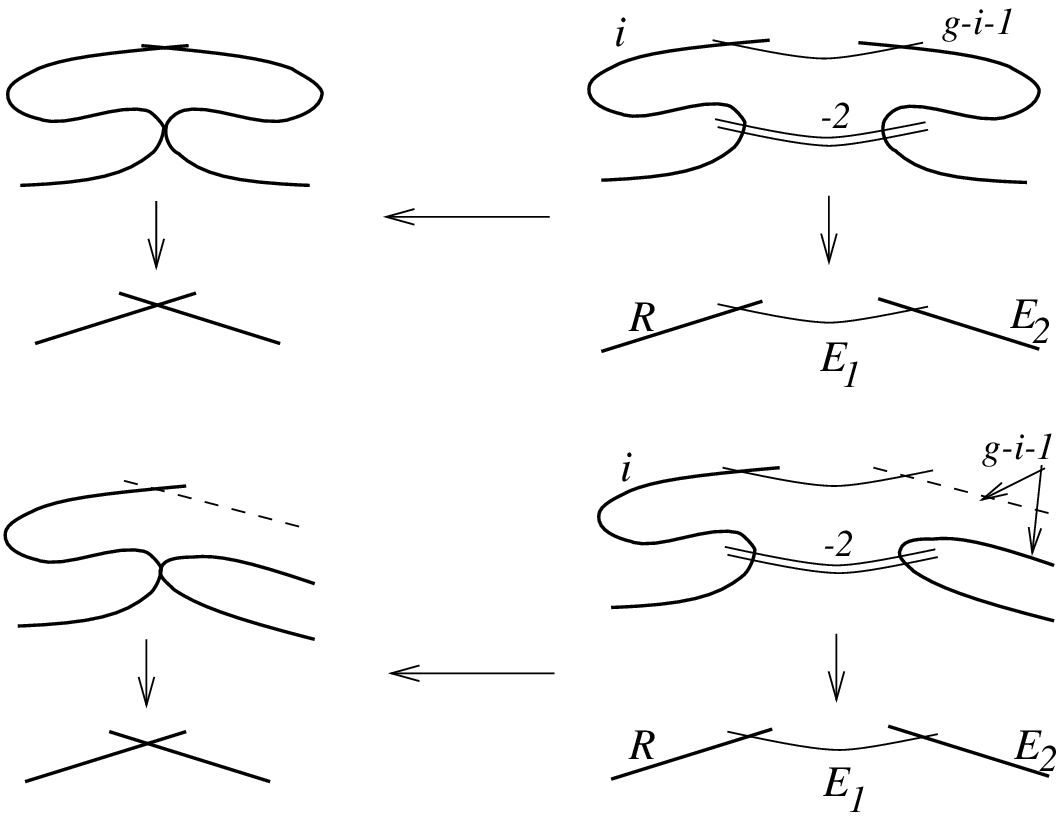}{3.25in}{2in}$$
\caption{Coefficients for ramification index 1}
\label{coef2.fig}
\end{figure}

In each of these cases, 
the fiber $T$ of $\widehat{Y}$ consists of two rational curves $E_1$ and $E_2$,
and the root $R=E_1^-$ (cf.~Fig.~\ref{coef2.fig}). There are
no nonreduced components in $T$, so the contribution to $\mathfrak{S}$
consists of two summands of the first type ($E,E^-$ nonreduced),
plus a quasi-admissible adjustment of $\mu(T)=1$ for $\Delta{\mathfrak{T}}_{5,i}$,
and a ramification adjustment of $g$ in both cases:
\[c_{k,i}=\frac{1}{4}\sum_{j=1,2}
\big(6(2+p_{\!\stackrel{\phantom{.}}{E_j}})(g-p_
{\!\stackrel{\phantom{.}}{E_j}})
-12g\big)+g\mu(T)+g\,\,\on{for}\,\,k=4,5.\]

The arithmetic genus of the nonreduced component of $\widehat{X}$ is
$-2$, and its intersection number with each of the neighboring components is
2. Setting $p_a(\widehat{\phi}^*R)=i$ forces $p_a(\widehat{\phi}^*E_2)= 
g-i-1$. Hence, $p_{\!\stackrel{\phantom{.}}{E_1}}=g-i-1$ and
$p_{\!\stackrel{\phantom{.}}{E_2}}=g-i-2$.
Substituting:
\[c_{k,i}=3(g-i)(i+1)-\frac{7g-3}{2}+g\mu(T),\]
\[\vspace*{-5mm}
c_{4,i}=3(i+1)(g-i)-\frac{7g-3}{2},\,c_{5,i}=3(i+1)(g-i)-\frac{7g-3}{2}+2g.\]

\subsubsection{Contribution of $\Delta{\mathfrak{T}}_{6,i}$: ramification 
index 2}
\label{contribution3}
It remains to consider the case of ramification index 2. 
Here there are four components $E$ besides the root $R$ in the special
fiber $T\subset \widehat{Y}$. Consequently, there
are four summands in $\mathfrak{S}$ corresponding to the $E_i$'s: $E_1$ and
$E_4$ yield summands of the first type ($E,E^-$ reduced),
$E_2$ yields a summand of the second type ($E$ nonreduced),
and $E_3$ yields a summand of the third type ($E^-$ nonreduced).

\begin{figure}[h]
$$\psdraw{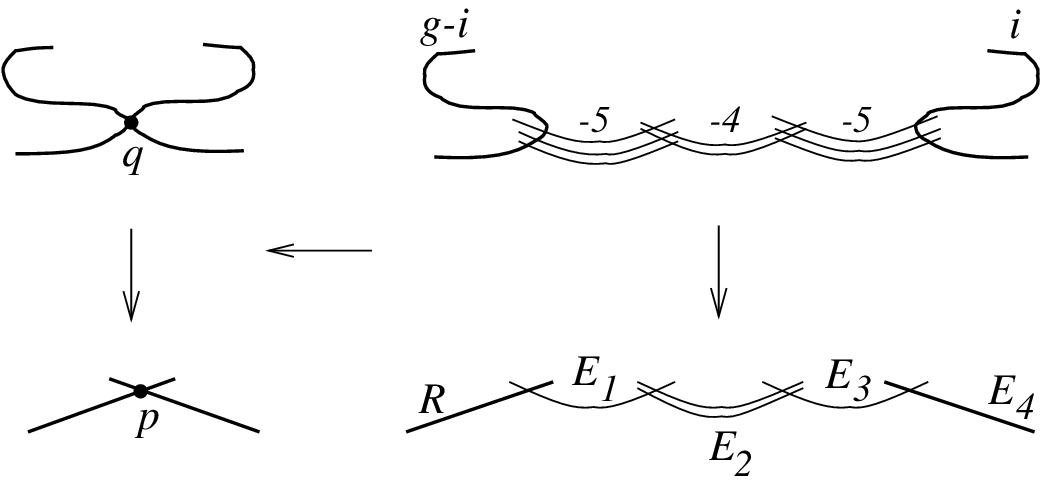}{3.5in}{1.5in}$$
\caption{Coefficients for ramification index 2}
\label{coef3.fig}
\end{figure}
Since $\mu(T)=0$, and the ramification adjustment is $3g$, we obtain for the
contribution of $\Delta{\mathfrak{T}}_{6,i}$ to $\mathfrak{S}$ the following expression:
\begin{eqnarray*}
c_{6,i}\!\!&\!\!=\!\!&\!\!\frac{1}{4}\big(6(2+
p_{\!\stackrel{\phantom{.}}{E_1}})(g-
p_{\!\stackrel{\phantom{.}}{E_1}})-12g\big)+
\frac{1}{4}\big(6(2+p_{\!\stackrel{\phantom{.}}{E_4}})(g-
p_{\!\stackrel{\phantom{.}}{E_4}})-12g\big)+\\
\!\!&\!\!+\!\!&\!\!\frac{1}{4}\left(3(
p_{\!\stackrel{\phantom{.}}{E_2}}+5)(2g-
p_{\!\stackrel{\phantom{.}}{E_2}}-1)-19g+6\right)
-\frac{1}{4}\left(3(p_{\!\stackrel{\phantom{.}}{E_3}}+2)^2+5g-6\right)+3g.
\end{eqnarray*}
The arithmetic genera of the components in $\widehat{X}$ are denoted in
the Fig.~\ref{coef3.fig}. It is easy to see that
 $p_{\!\stackrel{\phantom{.}}{E_4}}
=i$, $p_{\!\stackrel{\phantom{.}}{E_3}}=i-3$,
$p_{\!\stackrel{\phantom{.}}{E_2}}=i-2$, 
$p_{\!\stackrel{\phantom{.}}{E_1}}=i-2$. Finally,
\[c_{6,i}=\frac{9}{2}i(g-i)-\frac{3}{2}(g-1).\]
\subsection{Proof of Proposition~\ref{bogomolov1}} 
\label{Proof}
In the above discussion
we calculated the contributions of the boundary divisors 
$\Delta{\mathfrak{T}}_{k,i}$ to the sum $\mathfrak{S}$, so that
$\mathfrak{S}=\sum_{k,i}c_{k,i}$ with $k=1,...,6$, and the
corresponding limits for the index $i$ (cf.~Prop.~\ref{boundary}).
It is now clear what the divisor
$\cal{E}$ should be. We set
$\cal{E}:=\sum_{k,i}c_{k,i}\Delta{\mathfrak{T}}_{k,i}$, and thus,
$\mathfrak{S}=\cal{E}|_B$,
\[\Rightarrow\,\,\,
(7g+6)\lambda|_B=g\delta|_B+\cal{E}|_B+\frac{g-3}{2}(4c_2(V)-c_1^2(V)).\]
Using the restrictions on the index $i$ for each type of boundary divisor
$\Delta{\mathfrak{T}}_{k,i}$, one can easily deduce that all coefficients $c_{k,i}> 0$.
For instance, when $i=1,...,[g/2]$:
\[c_{6,i}=\frac{9}{2}i(g-i)-\frac{3}{2}(g-1)>\frac{9}{2}1\cdot
(g-1)-\frac{3}{2}(g-1)=3(g-1)>0.\]
In other words, $\cal{E}$ is an effective rational linear
combination of boundary divisors in $\overline{\mathfrak{T}}_g$, which by
construction does not contain $\Delta{\mathfrak{T}}_0.\,\,\,\qed$ 

\subsection{The slope bound $7+6/g$ and a relation 
restricted to the base curve $B$}
\label{slopebound}
Recall that a vector bundle $V$ of rank 2 is {\it Bogomolov
semistable} if $4c_2(V)\geq c^2_1(V)$. 

\begin{prop}[$7+6/g$ bound]
 For a general base curve $B$, if the canonically associated vector bundle $V$ 
is Bogomolov semistable, then the slope of $X/\!_{\displaystyle{B}}$ is 
bounded by
\[\frac{\delta|_B}{\lambda|_B}\leq 7+\frac{6}{g}\cdot\]
\label{7+6/g Bogomolov}
\end{prop}

\vspace*{-5mm}\begin{proof} The statement follows directly from 
Prop.~\ref{bogomolov1}. Indeed, since $\cal{E}$ is effective, then  
$\cal{E}|_B\geq 0$. By hypothesis, $4c_2(V)-c^2_1(V)\geq
0$, and $g\geq 3$. Hence, $(7g+6)\lambda|_B\geq g\delta|_B.$ \end{proof}

\begin{cor} For a general base curve $B$ 
the following relation holds true:
\begin{eqnarray*}
\!(7g+6)\lambda|_B&\!\!\!=\!\!\!&g\delta_0|_B+
\frac{g-3}{2}\left(4c_2(V)-c_1^2(V)\right)\\
&\!\!\!+\!\!\!\!\!\!&
\sum_{i=1}^{[(g-2)/2]}\frac{3}{2}(i+2)(g-i)\delta_{1,i}|_B+
\sum_{i=1}^{g-2}\frac{3}{2}(i+2)(g-i)\delta_{2,i}|_B\\
&\!\!\!+&\sum_{i=1}^{[g/2]}\frac{3}{2}(i+1)(g-i+1)\delta
_{3,i}|_B+\sum_{i=1}^{[(g-1)/2]}
\big(3(i+1)(g-i)-\frac{g-3}{2}\big)\delta_{4,i}|_B\\
&\!\!\!+&\sum_{i=1}^{g-1}\big(3(i+1)(g-i)-\frac{g-3}{2}\big)
\delta_{5,i}|_B+\sum_{i=1}^{[g/2]}
\big(\frac{9}{2}i(g-i)-{\frac{g-3}{2}}\big)\delta_{6,i}|_B.
\end{eqnarray*}
\label{analog1}\vspace*{-7mm}
\end{cor}

\noindent{\it Proof.} This is an immediate consequence of the established
relation in Prop.~\ref{bogomolov1}. We replace $\delta$ by the
linear combination (\ref{divisorrel}) of the boundary 
classes of $\overline{\mathfrak{T}}_g$, and write
\[(7g+6)\lambda=g\delta_0|_B+\sum_{k,i}\widetilde{c}_{k,i}\delta_{k,i}|_B
+\frac{g-3}{2}(4c_2(V)-c_1^2(V)),\]
for some new coefficients $\widetilde{c}_{k,i}$. Recall that
$\on{mult}_{\delta}(\delta_{k,i})$ denotes the {\it multiplicity}
of $\delta_{k,i}$ in $\delta$, so that
$\widetilde{c}_{k,i}=c_{k,i}+\on{mult}_{\delta}(\delta_{k,i})g$.
For example, the coefficient of $\delta_{1,i}$ is
\[\widetilde{c}_{1,i}=\left\{\frac{3}{2}(i+2)(g-i)-3g\right\}+3g=
\frac{3}{2}(i+2)(g-i),\]
or the coefficient of $\delta_{5,i}$ is
\[\widetilde{c}_{5,i}=\left\{3(i+2)(g-i)-\frac{7g-3}{2}+2g\right\}+g=
3(i+1)(g-i)-\frac{g-3}{2}.\,\,\,\qed\]

\medskip\section*{10. Generalized Index Theorem and Upper Bound}

\setcounter{section}{10}
\setcounter{subsection}{0} 
\setcounter{subsubsection}{0} 
\setcounter{lem}{0}
\setcounter{thm}{0}
\setcounter{prop}{0}
\setcounter{defn}{0} 
\setcounter{cor}{0} 
\setcounter{conj}{0} 
\setcounter{claim}{0} 
\setcounter{remark}{0}
\setcounter{equation}{0}
\label{Index1}

\begin{prop}[Index Theorem on $\widehat{X}$] For a general base curve $B$ and
for the rank 2 vector bundle $V$ on $\widehat{Y}$, we have
$9c_2(V)-2c_1^2(V)\geq 0.$
\label{genindex}
\end{prop}
\begin{proof} 
The proof is identical to that of 
Theorem~\ref{indextheorem}. One considers the divisor $\eta$ on 
$\widehat{X}$ defined by
\[\eta:=\left(\zeta+\frac{1}{3}\pi^*c_1(V)\right)\big|_
{\widehat{X}},\]
and shows that $\eta$ kills the pullback of any divisor on $\widehat{Y}$.
In particular, $\eta$ kills an ample divisor on $\widehat{X}$.
By the index theorem on $\widehat{X}$, $\eta^2
\leq 0$. From expression (\ref{genX}),
this can be also written as $9c_2(V)-2c_1^2(V)\geq 0.$ \end{proof}

\medskip
As in Section~7, the index theorem on $\widehat{X}$ suggests to
replace the Bogomolov difference $4c_2(V)-c^2_1(V)$ by another
linear combination of $c_2(V)$ and $c^2_1(V)$, which would behave in a more
``predictable'' way, namely, by $9c_2(V)-2c_1^2(V)$. In the process of doing
so, the only way to eliminate the unnecessary global
terms $d$ and $c$ from a relation among $\lambda|_B$ and $\delta|_B$ 
is to subtract:
$36(g+1)\lambda|_B-(5g+1)\delta|_B.$

\begin{prop} For a general base curve $B$ and an effective rational 
combination $\cal{E}^{\prime}$ of the
boundary divisors $\Delta{\mathfrak{T}}_{k,i}$, not containing $\Delta{\mathfrak{T}}_0$, we have:
\[36(g+1)\lambda|_B=(5g+1)\delta|_B+\cal{E}^{\prime}|_B+(g-3)
\big(9c_2(V)-2c_1^2(V)\big).\]
\label{indexrelation}\vspace*{-5mm}
\end{prop}
Note the apparent similarity between this relation and
Prop.~\ref{bogomolov1}. One may use the latter to prove the former, but
the calculations are not simpler than if one starts from scratch. We will
show a sketch of this proof, leaving the details to the reader, and referring
to the proof of Prop.~\ref{bogomolov1} for comparison.
\begin{proof} We denote by $\mathfrak{S}^{\prime}$ the difference
\[\mathfrak{S}^{\prime}:=
36(g+1)\lambda|_B-(5g+1)\delta|_B -(g-3)\left(9c_2(V)-2c_1^2(V)\right).\]
Substituting for $\delta|_B,\lambda|_B$ and $c_1^2(V)$ the
corresponding identities from Prop.~\ref{hatlambda} and
Example 8.1, and recalling that $c=g+2$ (cf.~Lemma~\ref{adjunction}), 
we write $\mathfrak{S}^{\prime}$ as
\begin{eqnarray*}
\mathfrak{S}^{\prime}&=
&-\sum_E\left\{m_{\!\stackrel{\phantom{.}}{E}}\left(8\Gamma^2+
8(g+2)\Gamma\Theta+9(g+1)\Theta^2\right)_{\!\stackrel{\phantom{.}}{E}}+6(g-1)
\right\}\\
&&+\,(5g+1)\left(\sum_T \mu(T)+\sum_{\on{ram}1}1 +\sum_{\on{ram}3} 3\right).
\end{eqnarray*}
As in Lemma~\ref{bogomolov1}, we group the above summands for every
special fiber in $\widehat{X}$, and correspondingly, for every
chain $T$ in $\widehat{Y}$. Recall Corollary~\ref{constants},
and the computations of the arithmetic genera
$p_{\!\stackrel{\phantom{.}}{E}}$ in the previous section:
\begin{eqnarray*}
\mathfrak{S}^{\prime}
&\!\!=\!\!&(5g+1)\big(\sum_T \mu(T)+\sum_{\on{ram}1} 1+\sum_{\on{ram}2}
3\big)+\!\!\sum_{E,E^-\on{red}}\!\!\left(8(p_{\!\stackrel{\phantom{.}}{E}}+2)
(g-p_{\!\stackrel{\phantom{.}}{E}})-3(5g+1)\right)\\
&\!\!-\!\!\!\!&\sum_{E^-\on{nonred}}\!\!\!\!\!\left(4
(p_{\!\stackrel{\phantom{.}}{E}}+2)^2+6(g-1)\right)\,+
\sum_{E\on{nonred}}\!\!\left(4(p_{\!\stackrel{\phantom{.}}{E}}+5)
(2g-1-p_{\!\stackrel{\phantom{.}}{E}})-12(g-1)\right).
\end{eqnarray*}
With this at hand, it is not hard to calculate the
contributions $d_{k,i}$ of each boundary component $\Delta{\mathfrak{T}}_{k,i}$
to the sum $\mathfrak{S}^{\prime}$:
\[\begin{array}{|l|l|}\hline
\!d_{1,i}\stackrel{\phantom{l}}{=}8(i+2)(g-i)\phantom{+1}\!-3(5g+1)\!&
\!d_{4,i}=16(i+1)(g-i)-2(g-3)-3(5g+1)\\
\!d_{2,i}\stackrel{\phantom{l}}{=}8(i+2)(g-i)\phantom{+1}-\!2(5g+1)\!&
\!d_{5,i}\stackrel{\phantom{l}}{=}16(i+1)(g-i)-2(g-3)-\phantom{3}(5g+1)\\
\!d_{3,i}\stackrel{\phantom{l}}{=}8(i+1)(g-i+1)-(5g+1)\!&
\!d_{6,i}\stackrel{\phantom{l}}{=}24i(g-i)-(5g+1).\phantom{\big)}\\\hline
\end{array}\]
\label{d_{k,i}table}

Let $\cal{E}^{\prime}=\sum_{k,i}d_{k,i}\Delta{\mathfrak{T}}_{k,i}$. Then
$\mathfrak{S}^{\prime}=\cal{E}^{\prime}|_B$, and the desired relation would be
established if $\cal{E}^{\prime}$ is effective. Given the
restrictions on the indices $i$ of the coefficients $d_{k,i}$ in
Prop.~\ref{boundary}, one
easily shows that all $d_{k,i}>0.$ \end{proof}

\begin{prop}[Maximal Bound]
For a general base curve $B$, the slope  satisfies:
\[\frac{\delta}{\lambda}\leq \frac{36(g+1)}{5g+1},\]
with equality if and only all fibers of $X$ are irreducible 
curves, and either $g=3$ or 
the divisor $\eta$ on the total space of ${X}$ is numerically zero.
\label{genmaximal}
\end{prop}

\begin{proof} From the Index Theorem on $\widehat{X}$, it follows
that $9c_2(V)-2c_1^2(V)\geq 0$. Since $\cal{E}^{\prime}$ is
effective, $\cal{E}^{\prime}|_B\geq 0$. Then Prop.~\ref{indexrelation}
implies $36(g+1)\lambda|_B\geq(5g+1)\delta|_B$, with equality exactly
when $9c_2(V)-2c_1^2(V)=0$ and $\cal{E}^{\prime}|_B=0$. The latter
means that $B\cap \Delta{\mathfrak{T}}_{k,i}=\emptyset$ because all coefficients
$d_{k,i}$ of $\cal{E}^{\prime}$ are strictly positive. In other words,
the family $\widehat{X}$ has only {\it irreducible} fibers
($B\cap\Delta{\mathfrak{T}}_{0}\not = \emptyset$). 
This takes us back to Section~7, 
where we presented the global calculation 
on the triple cover $X\rightarrow Y$. There we concluded that
the {\it index condition} $9c_2(V)-2c_1^2(V)=0$ was equivalent to
$\eta\equiv 0$ on $X$($=\widehat{X}$), or the genus $g=3$. \end{proof}

\begin{cor} For a general base curve $B$, 
\begin{eqnarray*}\vspace*{-1mm}
\!36(g+1)\lambda|_B&\!\!\!\!=\!\!\!\!&(5g+1)\delta_0|_B+
(g-3)\left(9c_2(V)-2c_1^2(V)\right)\\
&\!\!\!\!+\!\!\!\!\!&\sum_{i=0}^{[(g-2)/2]}8(i+2)(g-i)\delta_{1,i}|_B+
\sum_{i=1}^{g-2}8(i+2)(g-i)\delta_{2,i}|_B\\
&\!\!\!\!+\!\!\!\!\!&\sum_{i=1}^{[g/2]}8(i+1)(g-i+1)\delta_{3,i}|_B+
\!\!\!\sum_{i=1}^{[(g-1)/2]}\!\!\!\!\!
\big(16(i+1)(g-i)-2(g-3)\big)\delta_{4,i}|_B\\
&\!\!\!\!+\!
\!\!\!&\sum_{i=1}^{g-1}\big(16(i+1)(g-i)-2(g-3)\big)\delta_{5,i}|_B+
\sum_{i=1}^{[g/2]}24i(g-i)\delta_{6,i}|_B.
\end{eqnarray*}
\label{analog2}
\end{cor}
\label{page analog2}

\noindent{\it Proof.} We only need to substitute the known expressions
for the divisors $\cal{E}^{\prime}$ and $\delta$ into 
Prop.~\ref{indexrelation}:
\[36(g+1)\lambda|_B=(5g+1)\delta_0|_B+\sum_{k,i}\big((5g+1)\on{mult}_{\delta}
(\delta_{k,i})+d_{k,i}\big)+(g-3)\big(9c_2(V)-2c_1^2(V)\big).\]
The rest is a simple calculation. For example, the total coefficient 
$\widetilde{d}_{3,i}$ of $\delta_{3,i}$ equals
\begin{eqnarray*}
d_{3,i}+(5g+1)\on{mult}_{\delta}(\delta_{3,i})&=&\{8(i+1)(g-i+1)-(5g+1)\}+
(5g+1)\cdot 1\\
&=&8(i+1)(g-i+1).\,\,\,\qed
\end{eqnarray*}

\medskip 
\section*{11. Extension to an Arbitrary Base $B$}

\setcounter{section}{11}
\setcounter{subsection}{0} 
\setcounter{subsubsection}{0} 
\setcounter{lem}{0}
\setcounter{thm}{0}
\setcounter{prop}{0}
\setcounter{defn}{0} 
\setcounter{cor}{0} 
\setcounter{conj}{0} 
\setcounter{claim}{0} 
\setcounter{remark}{0}
\setcounter{equation}{0}
\label{arbitrary}

We extend now the results of
Sect.~8-10 to arbitrary nonisotrivial
families $X\!\rightarrow \!B$ with smooth trigonal general member. The
essential case is when $B$ is {\it not} tangent to the boundary
$\Delta{\mathfrak{T}}$, from which the remaining cases easily follows.

\subsection{The base curve $B$ not tangent to $\Delta\mathfrak{T}$}
\label{nontangentB}
We now drop the hypothesis of the base curve $B$
intersecting the boundary divisors in general points. Instead, for now we only
assume that the base curve $B$ is not tangent
the boundary $\Delta{\mathfrak{T}}$.
This means that all special fibers of $X$ locally look like the
general ones (cf.~Fig.~\ref{coef1.fig}--\ref{coef3.fig}). Therefore,
from the quasiadmissible cover
$\widetilde{X}\rightarrow \widetilde{Y}$ we can construct an effective cover
$\widehat{\phi}:\widehat{X}\rightarrow \widehat{Y}$ of {\it smooth}
surfaces $\widehat{X}$ and $\widehat{Y}$.
(The {\it smoothness} indicates
that $B$ is {\it not} tangent to any $\Delta{\mathfrak{T}}_{k,i}$. Otherwise,
there would be a higher local multiplicity $xy=t^n$ near a node 
of a special fiber $C_{X}$, $n> 1$. Hence  $\widehat{X}$ would be 
obtained locally via a base change from a smooth surface, but $\widehat{X}$
would have a singular total space.)

Now the special fibers of $\widehat{Y}$ are {\it trees} $T$ (rather than 
just chains) of reduced smooth rational
curves with occasional nonreduced rational components of multiplicity 2.
The latter occur again exactly for each singular point in
$C_{\widetilde{X}}$ of 
ramification index 2 under the quasiadmissible cover $\widetilde{\phi}:
\widetilde{X}\rightarrow\widetilde{Y}$ (cf.~Fig.~\ref{coef3.fig}).

\smallskip
The notation and conventions from Sections 
\ref{conventions} are also valid here. In particular, for any tree $T$, 
we fix one of its end (nonreduced) components to be its root $R$, and we define
as before the functions $m,\theta,\gamma$ on the components $E$ of $T$.
 Moreover, since Lemma
\ref{technical} can be applied also for any tree $T$, the calculations of
$\lambda_{\widehat{X}},\kappa_{\widehat{X}}$ and $\delta$ in 
Prop.~\ref{hatlambda}
and Cor.~\ref{hatdelta} go through without any modifications.

\smallskip Finally, we wish to extend all results of 
Sections~8-10 over the new base $B$. 
The only difference arises in the final calculation of the
coefficients $c_{k,i}$ and $d_{k,i}$. The fiber $C_X$ in $X$, corresponding to
a tree $T$, may now lie 
in the intersection of {\it several} boundary divisors 
$\Delta{\mathfrak{T}}_{k,i}$. Such a trigonal curve $C_X$ is called a {\it special
boundary} curve.
Accordingly, its contribution $c_{\!\stackrel{\phantom{.}}{T}}$ to 
$\mathfrak{S}$ (or $d_{\!\stackrel{\phantom{.}}{T}}$ to $\mathfrak{S}^{\prime}$)
will be distributed among these divisors 
$\Delta{\mathfrak{T}}_{k,i}$'s, rather than
just yielding  a single coefficient $c_{k,i}$ (or $d_{k,i}$) as before.

\begin{figure}[h]
$$\psdraw{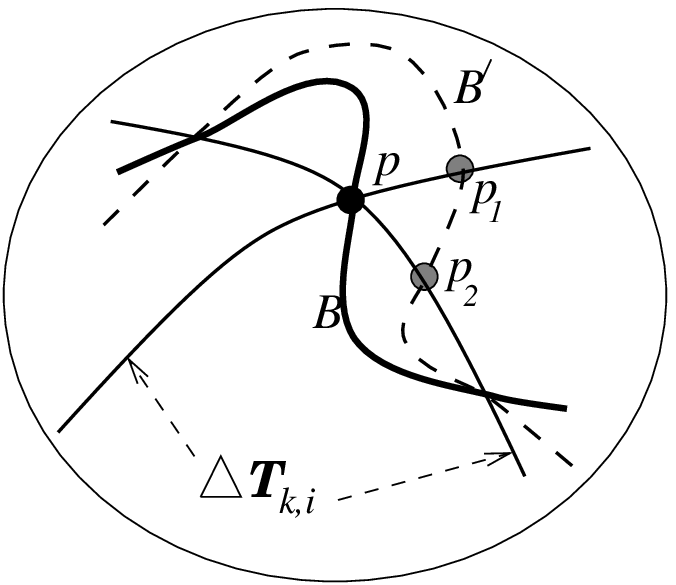}{1.6in}{1.5in}$$
\caption{Moving $B$}
\label{arbitrary1}
\end{figure}
This can be easily resolved. The idea is to replace any special singular
fiber in $\widehat{X}$ by a suitable combination of {\it general} fibers,
without changing the sums $\mathfrak{S}$ and $\mathfrak{S}^{\prime}$. We can
imagine this as ``moving'' the base curve $B$ in $\overline{\mathfrak{T}}_g$
{\it away from} the special singular locus of $\overline{\mathfrak{T}}_g$,
and replacing it with a {\it general} base curve $B^{\prime}$, as defined
in Section~8. For example, in Fig.~\ref{arbitrary1}
the base $B$ passes through a point $p$ in the intersection of
two boundary divisors $\Delta\mathfrak{T}_{k,i}$. Two new general points
$p_1$ and $p_2$, each lying on a single $\Delta\mathfrak{T}_{k,i}$, replace the
special point $p$, and thus $B$ moves to a {\it general} 
 curve $B^{\prime}$.

\begin{lem} Let  $C_X$ be a special boundary curve in
$\overline{\mathfrak{T}}_g$. Denote by 
$\alpha_{k,i}$ the degree of the point $[C_X]$ in the intersection
$\Delta{\mathfrak{T}}_{k,i}\cdot B$. Then the contributions of
$T=\widehat{\phi}(C_{\widehat{X}})$ to 
$\mathfrak{S}$ and to $\mathfrak{S}^{\prime}$ are 
$c_{\!\stackrel{\phantom{.}}{T}}=\sum_{k,i}\alpha_{k,i}c_{k,i}$ and
$d_{\!\stackrel{\phantom{.}}{T}}=\sum_{k,i}\alpha_{k,i}d_{k,i}$, respectively.
\label{contributions d_Tc_T}
\end{lem}
\noindent{\it Proof:} Rewrite $\mathfrak{S}$ and
$\mathfrak{S}^{\prime}$ as sums over the non-root components $E$ of the special
trees $T$:
\begin{eqnarray*}
\mathfrak{S}&\!\!=\!\!&
\sum_{E,E^-\on{red}}F_1(p_{\!\stackrel{\phantom{.}}{E}})+
\sum_{E^-\on{nonred}}\!\!\!\!\!F_2(p_{\!\stackrel{\phantom{.}}{E}})
+\sum_{E\on{nonred}}\!\!F_3(p_{\!\stackrel{\phantom{.}}{E}})
+gH,\\
\mathfrak{S}^{\prime}&\!\!=\!\!&
\sum_{E,E^-\on{red}}G_1(p_{\!\stackrel{\phantom{.}}{E}})+
\sum_{E^-\on{nonred}}\!\!\!\!\!G_2(p_{\!\stackrel{\phantom{.}}{E}})
+\sum_{E\on{nonred}}\!\!G_3(p_{\!\stackrel{\phantom{.}}{E}})+(5g+1)H,
\end{eqnarray*}
where $H=\sum_T\mu(T)+\sum_{\on{ram}1}1+\sum_{\on{ram}2}3$ is the
quasi-admissible and effective adjustment, and the functions $F_i$ and
$G_j$ are quadratic polynomials in $p_{\!\stackrel{\phantom{.}}{E}}$
with linear coefficients in $g$. 
Recall that in these sums each non-root component $E$ appears
exactly once, and $p_{\!\stackrel{\phantom{.}}{E}}$ is the arithmetic
genus of the inverse image $\widehat{\phi}^*(T(E))$ of the
subtree $T(E)$ generated by $E$.

\smallskip
There is a simple way to recognize the boundary divisors 
$\Delta{\mathfrak{T}}_{k,i}$ in which
a special trigonal fiber $C_X$ lies. Consider the 
corresponding ``effective'' fiber $C_{\widehat{X}}=\widehat{\phi}^*T$
in $\widehat{X}$. For any non-root component $E$ in $T$ there are
two possibilities: either $\widehat{\phi}^*E$ and 
$\widehat{\phi}^*{E^-}$ are both reduced, or
$E$ is part of a  chain of length 3 or 5, constructed to
resolve ramifications in the quasi-admissible fiber $C_{\widetilde{X}}$.
 
\smallskip
\begin{figure}[h]
$$\psdraw{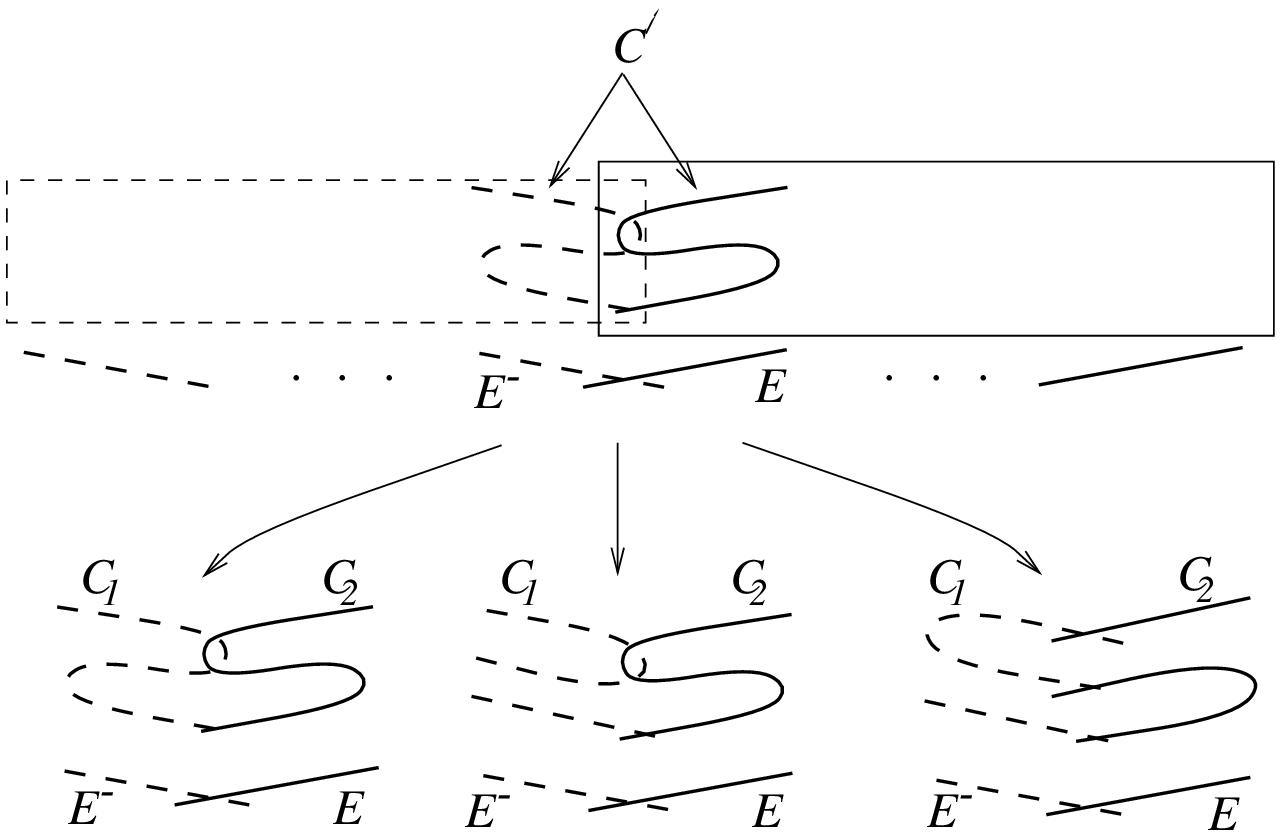}{3.8in}{1.8in}$$ \vspace*{-2mm}
\caption{$E\not\subset$ chain $\rightarrow\,\,
\alpha_{1,i},\alpha_{2,i},\alpha_{3,i}$}
\label{notchain}
\end{figure}
\subsubsection{Contributions to the degrees 
$\alpha_{1,i},\alpha_{2,i},\alpha_{3,i}$}
Consider the first situation, and denote by $C^{\prime}$ the preimage
$\widehat{\phi}^*E\cup \widehat{\phi}^*{E^-}$ in $\widehat{X}$. Thus,
$C^{\prime}$ corresponds
to a general member of $\Delta{\mathfrak{T}}_{1,i},\Delta{\mathfrak{T}}_{2,i},
\Delta{\mathfrak{T}}_{3,i}$, possibly of lower genus (cf.~Fig.~\ref{notchain}).
As part of the fiber $C_{\widehat{X}}$, the curve $C^{\prime}$ is
represented for simplicity by the triple intersection of two {\it smooth}
trigonal curves (in $\Delta{\mathfrak{T}}_{1,i}$), but it could have 
corresponded to any general member of $\Delta{\mathfrak{T}}_{2,i}$ 
or $\Delta{\mathfrak{T}}_{3,i}$. The solid box
encompasses the preimage $\widehat{\phi}^*T(E)$, and the dashed box
encompasses the preimage of the rest, $\widehat{\phi}^*\big(T-T(E)\big)$.
Each of these boxes represents a limit of a quasi-admissible
curve, $C_1$ or $C_2$,
 which is naturally a triple cover of ${\proj}^1$. Thus, we can
{\it ``smoothen''}
each box to such a curve $C_i$. As a result we obtain
a quasiadmissible curve $C_1\cup C_2$ of total genus $g$, which
corresponds to a general member 
of $\Delta{\mathfrak{T}}_{1,i},\Delta{\mathfrak{T}}_{2,i}$ or $\Delta{\mathfrak{T}}_{3,i}$. Depending on
which divisor $\Delta{\mathfrak{T}}_{k,i}$ is evoked, there is a 
corresponding contribution of $1$ to the
coefficient $\alpha_{k,i}$: $[C_X]\in\Delta{\mathfrak{T}}_{k,i}$.

\smallskip
Note that the arithmetic genus of $C_2$ is the previously
defined $p_{\!\stackrel{\phantom{.}}{E}}$. The contribution of $E$
to $\mathfrak{S}$ is $F_1(p_{\!\stackrel{\phantom{.}}{E}})$ plus the possible
quasi-admissible adjustment in $\mu(T)$ needed to obtain $\widetilde
{\phi}^*(E\cup E^-)$. 
In view of the above ``smoothening'', this can be thought of as
the contribution of $C_2$ in the effective curve $C_1\cup C_2$, and by
Prop.~\ref{bogomolov1} this is exactly the coefficient $c_{k,i}$.
The same argument works in the case of $\mathfrak{S}^{\prime}$ from 
Prop.~\ref{indexrelation}. We conclude that $\alpha_{k,i}$ (for $k=1,2,3$)
equals the number of $c_{k,i}$'s  and $d_{k,i}$'s in $\mathfrak{S}$ and
$\mathfrak{S}^{\prime}$, respectively. 

\subsubsection{Contributions to the degrees 
$\alpha_{4,i},\alpha_{5,i},\alpha_{6,i}$} 
\begin{figure}[h]
$$\psdraw{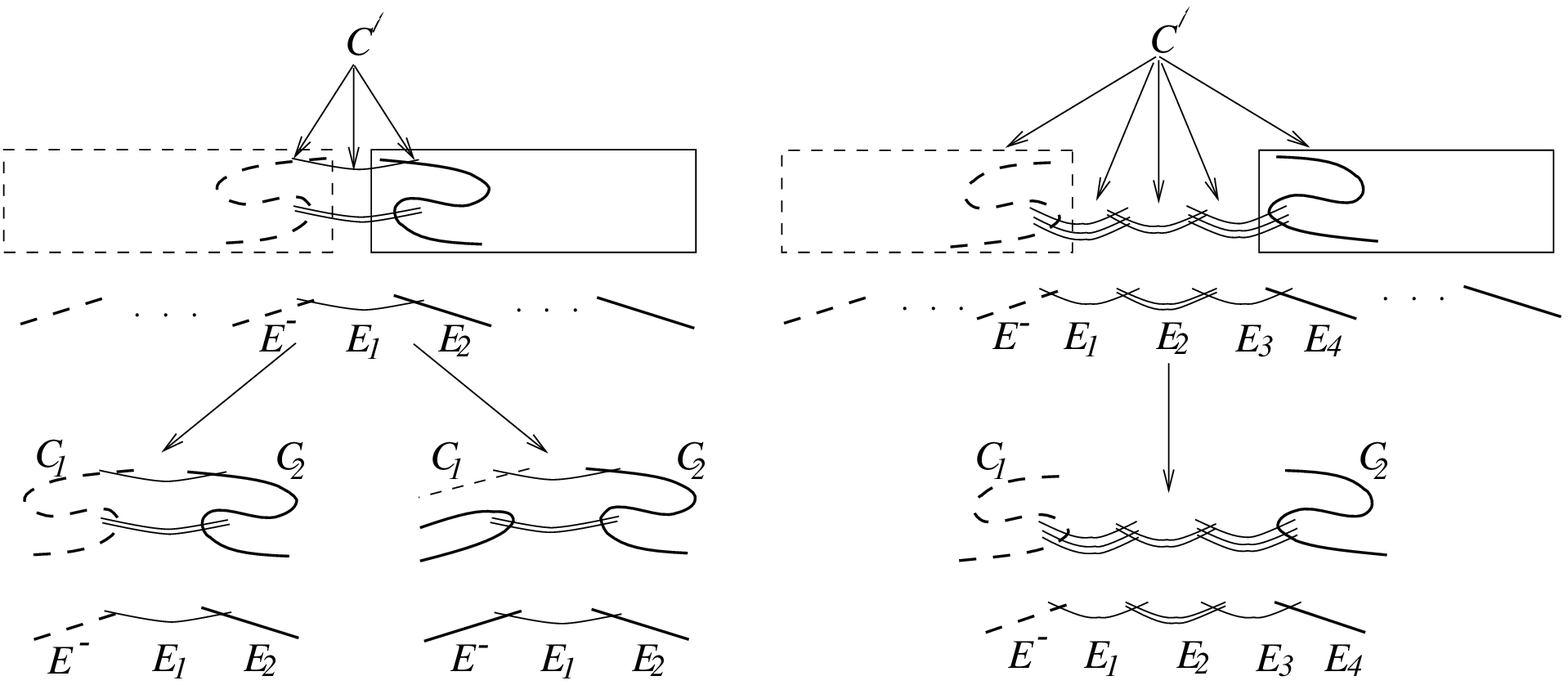}{5.5in}{2in}$$
\caption{$E\subset$ chain $\rightarrow\,\,
\alpha_{4,i},\alpha_{5,i},\alpha_{6,i}$}
\label{chain}
\end{figure}

We treat analogously
the remaining case
when the component $E$ is part of a chain of length 3 or 5.
Here, however, one must consider {\it simultaneously
 all} the components $E$ of $T$ participating in such a chain, and take a
quasi-admissible limit only {\it over the reduced} curves in $C_{\widehat{X}}$.
In Fig.~\ref{chain} one can see all three ramification cases, or
equivalently, the boundary divisors $\Delta{\mathfrak{T}}_{4,i},\Delta{\mathfrak{T}}_{5,i}$
and $\Delta{\mathfrak{T}}_{6,i}$. For simplicity, we have again depicted the
reduced components in $\widehat{X}$ by smooth trigonal curves,
which may not always be true for every tree $T$: they could, for instance,
be singular or reducible, but they will keep the ramification index 1 or 2
at the appropriate points.
\newline 
To see how $c_{k,i}$ and $d_{k,i}$ are obtained,
let us calculate, for example, the contributions of $E_1,E_2,E_3$ and
$E_4$ in the case of $\Delta{\mathfrak{T}}_{6,i}$. The inverse images in
$\widehat{X}$ of $T-T(E_1)$ and $T(E_4)$ are marked by dashed and
solid boxes, respectively.
We {\it ``smoothen''} each box by a smooth trigonal curve, $C_1$ or $C_2$, and
keep the inverse images of $E_1$,$E_2$ and $E_3$. Thus, we obtain
a general member ${C}^{\prime\prime}$ 
of $\Delta{\mathfrak{T}}_{6,i}$. The arithmetic genera, necessary to calculate the
contribution of ${C}^{\prime\prime}$ to $\mathfrak{S}$, are given from
right to left by: \[p_a(C_2)=p_{\!\stackrel{\phantom{.}}{E_4}},\,\,
p_{\!\stackrel{\phantom{.}}{E_3}}=p_{\!\stackrel{\phantom{.}}{E_4}}\!\!-3,\,\,
p_{\!\stackrel{\phantom{.}}{E_2}}=p_{\!\stackrel{\phantom{.}}{E_4}}\!\!-2,\,\,
p_{\!\stackrel{\phantom{.}}{E_1}}=p_{\!\stackrel{\phantom{.}}{E_4}}\!\!-2.\]
As in the proof of Prop.~\ref{bogomolov1}, we substitute these
in the sum $\mathfrak{S}$, and for $i=p_{\!\stackrel{\phantom{.}}{E_4}}$ we obtain
\[F_1(E_1)+F_1(E_4)+F_2(E_2)+F_3(E_3)+3g=\frac{9}{2}
p_{\!\stackrel{\phantom{.}}{E_4}}(g-p_{\!\stackrel{\phantom{.}}{E_4}})-
\frac{3}{2}(g-1)=c_{6,i}.\]

Combining all of the above
results, we conclude that the contributions of any tree $T$ to the sums
$\mathfrak{S}$ and $\mathfrak{S}^{\prime}$ are
$c_{\!\stackrel{\phantom{.}}{T}}=\sum_{k,i}\alpha_{k,i}c_{k,i}\,\,\on{and}\,\,
d_{\!\stackrel{\phantom{.}}{T}}=\sum_{k,i}\alpha_{k,i}d_{k,i}.$\qed 

\smallskip
This allows us to extend all results of
Sect.~\ref{Bogomolov1}--\ref{Index1} to the case of a base curve $B$
meeting transversally the boundary $\Delta{\mathfrak{T}}_g$.

\subsection{Extension to an arbitrary base $B$, not contained in
$\Delta{\mathfrak{T}}_g$}
\label{extension}

If the base curve $B$ happens to be {\it tangent} to a boundary divisor
$\Delta\mathfrak{T}_{k,i}$ at a point $[C_{X}]$, then over some node $p$ of the
corresponding tree $T=\widehat{\phi}(C_{\widehat{X}})$ {\it all}
local analytic multiplicities $m_q$ (cf.~Sect.~\ref{definition})
will be multiplied by the degree 
of tangency of $B$ and $\Delta\mathfrak{T}_{k,i}$. Fig.~\ref{Local
multiplicities} presents a few examples of possible fibers in $\widetilde{X}$:
\begin{figure}[h]
$$\psdraw{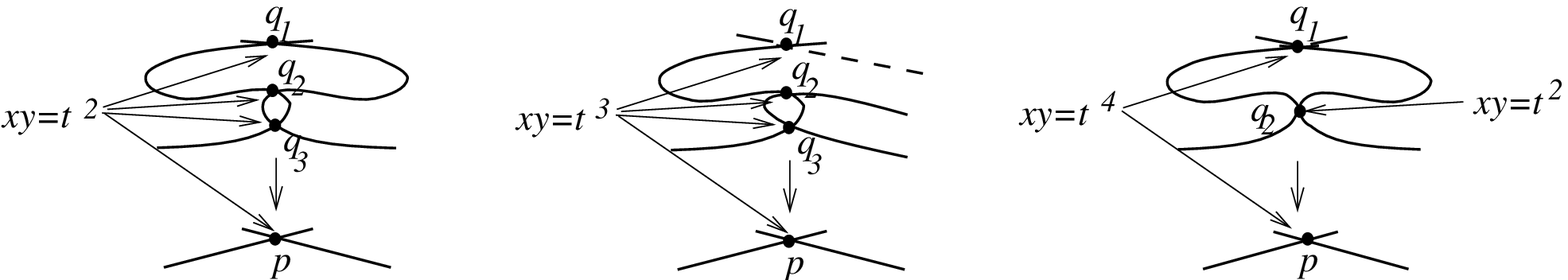}{5.4in}{1.2in}$$
\caption{Local multiplicities}
\label{Local multiplicities}
\end{figure}

\smallskip
In the nonramification cases of $\Delta\mathfrak{T}_{1,i},\Delta\mathfrak{T}_{2,i}$ and
$\Delta\mathfrak{T}_{3,i}$, 
this would force rational double points as singularities on
the total spaces of $\widehat{X}$ and $\widehat{Y}$, 
whereas in the ramification cases of
$\Delta\mathfrak{T}_{4,i},\Delta\mathfrak{T}_{5,i}$ and
$\Delta\mathfrak{T}_{6,i}$, one may arrive at surfaces $\widehat{X}$
and $\widehat{Y}$,  nonnormal  over some nonreduced fibers. 
But in both cases, one can roughly view the
corresponding fibers as being obtained by a base change from the general or
special fibers of
Sect.~8 and Sect. ~11.1. Alternatively,
one can go through the arguments of the paper for the new surfaces
$\widehat{X}$ and $\widehat{Y}$ (normalizing, if necessary), and
notice that all formulas (e.g. Euler characteristic formula for
$\lambda$, Index theorem on $\widehat{X}$, adjunction formula in $\proj V$,
etc.) hold for surfaces with double point singularities.

\smallskip
Thus, in effect, one may replace a given singular fiber $C_X$ by
a bunch of general boundary curves $C$, following
the procedure described in Section~\ref{nontangentB}. Furthermore, if
some of these general curves $C$ are ``multiple'' (i.e. $B$ is tangent to
$\Delta\mathfrak{T}_{k,i}$ at $[C]$), one may in turn replace each $C$
by several ``transversal'' general boundary curves, and refer to the
statements in Sections~8.3 and 9.3. The only notational
difference in this approach will appear in the definition of the invariants
$m,\theta$ and $\gamma$ from Sect.~8: now we have to allow
for them to be {\it rational}, rather than integral, due to possible
rational intersections $E\cdot E^-$. This will be
``compensated'' in the final calculations, which will take into account the
multiplicity of the corresponding fibers, and  roughly speaking,
will ``multiply back'' our invariants $\delta,\lambda$ and $\kappa$
by what they were divided by in the beginning of the calculations.

\bigskip
This concludes the proof of our results for all families of stable curves
$X\rightarrow B$ with  general smooth trigonal members.

\subsection{Statements of the results for any family $X\rightarrow B$}
\label{results}

In the following list of results, 
Theorems~\ref{7+6/g relation2} and \ref{maximal relation2} can be viewed as
local trigonal analogs of Cornalba-Harris's relation in the Picard group of the
hyperelliptic locus $\overline{\mathfrak{I}}_g$ (cf.~Theorem~\ref{CHPic}).
Similarly, Theorem~\ref{maximal bound2} is the analog of the $8+4/g$
maximal bound in the hyperelliptic case (cf.~Theorem~\ref{CHX}).

\begin{thm}[$7+6/g$ relation]
For any family $X\rightarrow B$ of stable curves with smooth
trigonal general member, if 
$V$ is the canonically associated to $X$ vector bundle of rank 2, then
the following relation holds true
\[(7g+6)\lambda|_B=g\delta|_B+\cal{E}|_B+\frac{g-3}{2}\left(4c_2(V)-c_1^2(V)
\right),\] 
where $\cal{E}$ is an effective rational linear combination of
boundary components of $\overline{\mathfrak{T}}_g$, 
not containing $\Delta{\mathfrak{T}}_0$. In particular,
\[(7g+6)\lambda|_B=g\delta_0|_B+\sum_{k,i}\widetilde{c}_{k,i}
\delta_{k,i}|_B+\frac{g-3}{2}\left(4c_2(V)-c_1^2(V)\right),\]
where $\widetilde{c}_{k,i}$ is
quadratic polynomial in $i$ with  linear coefficients in $g$, and it is 
determined by the geometry of $\delta_{k,i}$ (cf.~p.~\pageref{analog1}).
\label{7+6/g relation2}
\end{thm}

\begin{thm}[$7+6/g$ bound]
 For any nonisotrivial family $X\rightarrow B$ of stable curves with smooth
trigonal general member, if the canonically associated to $X$ 
vector bundle $V$ is Bogomolov semistable, then the slope of 
$X/\!_{\displaystyle{B}}$ is bounded from above by
\[\frac{\delta}{\lambda}\leq 7+\frac{6}{g}\cdot\vspace*{-5mm}\]
\label{7+6/g Bogomolov2}
\end{thm}
\begin{thm}[Index relation]
 For any family $X\rightarrow B$ of stable curves with smooth
trigonal general member, 
if $V$ is the canonically associated to $X$ vector bundle of rank 2, then
the following relation holds true
\[36(g+1)\lambda|_B=(5g+1)\delta|_B+\cal{E}^{\prime}|_B+(g-3)
\big(9c_2(V)-2c_1^2(V)\big),\]
where $\cal{E}^{\prime}$ is an effective rational combination of the
boundary divisors $\Delta{\mathfrak{T}}_{k,i}$, not containing 
$\Delta{\mathfrak{T}}_0$. In particular,
\[36(g+1)\lambda|_B=(5g+1)\delta_0|_B+\sum_{k,i}\widetilde{d}_{k,i}
\delta_{k,i}|_B+(g-3)\left(9c_2(V)-2c_1^2(V)\right),\]
where $\widetilde{d}_{k,i}$ is
quadratic polynomial in $i$ with  linear coefficients in $g$, and it is 
determined by the geometry of $\delta_{k,i}$ (cf.~p.~\pageref{page analog2}).
\label{maximal relation2}
\end{thm}

\begin{thm}[Maximal bound]
For any nonisotrivial family $X\rightarrow B$ of stable curves 
with smooth trigonal general member, 
the slope of $X/\!_{\displaystyle{B}}$ satisfies:
\[\frac{\delta}{\lambda}\leq \frac{36(g+1)}{5g+1},\]
with equality if and only all fibers of $X$ are irreducible 
curves, and either $g=3$ or 
the divisor $\eta$ on the total space of ${X}$ is numerically zero.
\label{maximal bound2}
\end{thm}
\label{list of theorems}

\subsection{What happens with the hyperelliptic locus
$\overline{\mathfrak{I}}_g$}

As we promised in Section~\ref{hyperelliptic locus}, we consider the
contribution of the hyperelliptic locus to the above theorems.
For any hyperelliptic curve $C$, we need to blow up a point on $C$
before it starts ``behaving'' like a trigonal curve in the quasi-admissible
and effective covers. Below we have shown what happens to a smooth
or general singular hyperelliptic curve (cf.~Fig.~\ref{hyperboundary}
for the admissible classification of the boundary locus $\Delta\mathfrak{I}_g$).

\begin{figure}[h]
$$\psdraw{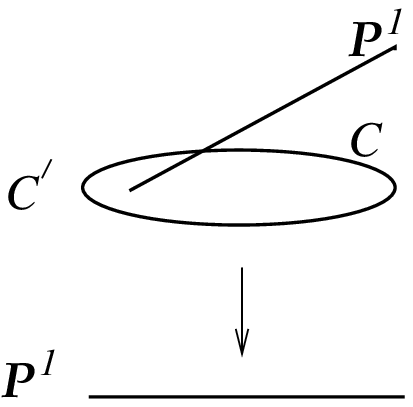}{1in}{1in}$$
\caption{$\overline{\mathfrak{I}}_g\cap\Delta\mathfrak{T}_0$}
\label{smoothhyper}
\end{figure}

\subsubsection{Smooth hyperelliptic curves} 
We blow up $C$ at a point, and thus add a smooth rational component $\proj^1$
to make it a triple cover $C^{\prime}$ (cf.~Fig.~\ref{smoothhyper}).
The quasi-admissible adjustment of $C$ is $\mu(C^{\prime})=1$. From here
on, $C$ will behave essentially like a smooth trigonal curve.
Therefore, in all relations $C$ is going to contribute $g$ or $(5g+1)$,
depending on what $\delta$ is multiplied by.

\subsubsection{Singular hyperelliptic curves in $\Delta\mathfrak{T}_{2,i}$
and $\Delta\mathfrak{T}_{5,i}$}
\begin{figure}[h]
$$\psdraw{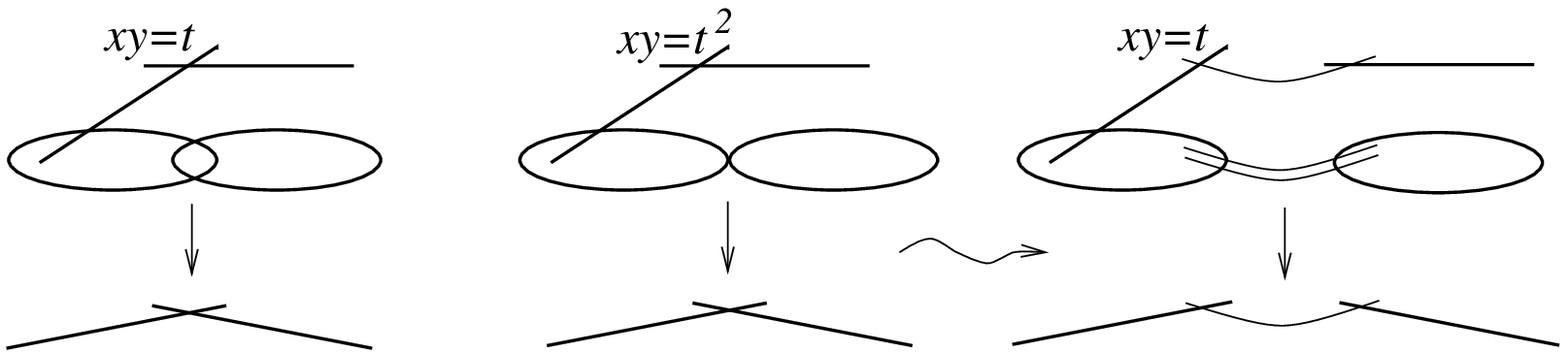}{5in}{1.1in}$$
\caption{$\overline{\mathfrak{I}}_g\cap\Delta\mathfrak{T}_{2,i}$\hspace{22mm}
{\sc Figure 47.} $\overline{\mathfrak{I}}_g\cap\Delta\mathfrak{T}_{5,i}$
\hspace*{25mm}}
\label{singularhyper}
\end{figure}

The necessary effective and quasi-admissible modifications are shown in
Fig.~\ref{singularhyper}--47.

In the first case, there are two hyperelliptic components intersecting
transversally in two points. For the quasi-admissible cover, we need
two ``smooth'' blow-ups, which makes $\mu=2$. From now on, this curve
will behave like a element of $\Delta\mathfrak{T}_{2,i}$, where $\mu_{2,i}=1$. 
Thus, the coefficient
in, say, the maximal bound relation will be: $\widetilde{d}_{2,i}+(5g+1)$,
due to the extra blow-up in $\mu$.

\addtocounter{figure}{1}
In the second case,  two hyperelliptic components meet transversally
in one point, but have a ramification index 1 at this point when viewed as
double covers. Fig.46 presents first the quasi-admissible modification:
as in the case of $\Delta\mathfrak{T}_{5,i}$, the local 
analytic multiplicity between
the two rational components is $2$, which means that we must have made
three ``smooth'' blow-ups and one ``singular'' blow-down. As a result,
$\mu=3$. From here on, this curve behaves exactly as a general member of
$\Delta\mathfrak{T}_{5,i}$.
Recall that $\mu_{5,i}=2$, and the extra $1$ in the hyperelliptic case
accounts for the one extra blow-up. Therefore, the coefficient
of this fiber $C$, say, in the maximal bound relation, will be
$\widetilde{d}_{5,i}+(5g+1)$. 

\smallskip
We conclude that a base curve $B$, passing through the hyperelliptic locus,
will contribute in the results listed in Section~\ref{results}
roughly $g$, or $(5g+1)$, times the number of elements in
$B\cap\overline{\mathfrak{I}}_g$. We cannot write the latter in the form
of a scheme-theoretic intersection, since $\overline{\mathfrak{I}}_g$ is of
a larger codimension in $\overline{\mathfrak{T}}_g$. 
\label{hypercalculations}

\smallskip
One can explain these extra summands in the expressions for $\lambda$ in
the following way. Recall the projection map $pr_1:\overline{\cal{H}}_{3,g}
\rightarrow \overline{\mathfrak{T}}_g$. The exceptional locus of $pr_1$
is the admissible boundary divisor $\Delta{\cal{H}}_{3,0}$, which is
blown down to the codimension 2 hyperelliptic locus $\overline{\mathfrak{I}}_g$
inside $\overline{\mathfrak{T}}_g$. For calculation purposes, it will
be easier to work instead with the space of  minimal quasi-admissible
covers $\overline{\cal Q}_{3,g}$, which replaces $\overline{\cal{H}}_{3,g}$.
The same situation of a blow-down
occurs, where the exceptional divisor in $\overline{\cal Q}_{3,g}$ consists
of reducible curves $C^{\prime}$, as shown in Fig.~\ref{smoothhyper}.

\begin{figure}[h]
$$\psdraw{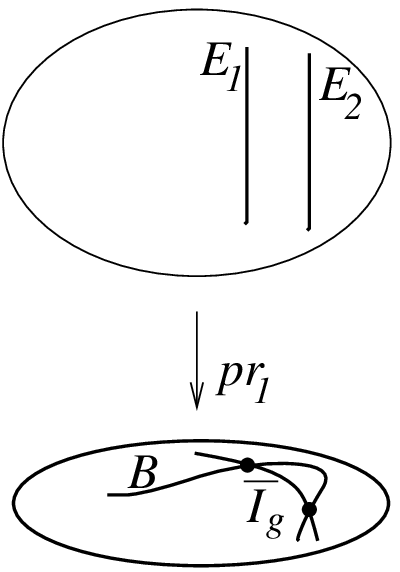}{1.3in}{1.8in}$$
\caption{$B\cap \overline{\mathfrak{I}}_g$}
\label{exceptional.fig}
\end{figure}

Let $D$ be the linear combination of divisors in $\overline{\mathfrak{T}}_g$
given by the restriction $\Delta|_{\overline{\mathfrak{T}}_g}$, and consider
a curve $B\subset \overline{\mathfrak{T}}_g$, intersecting the hyperelliptic
locus in finitely many points. 
By abuse of notation, we 
denote by $pr_1$ the projection from $\overline{\cal Q}_{3,g}$ to
$\overline{\mathfrak{T}}_g$. Then for the intersection $D\cdot B$ we have:
\[D\cdot B=pr_1^*(D)\cdot pr_1^*(B)=pr_1^*(D)\cdot(\overline{B}+\sum
E_j),\]
where $\overline{B}$ is the proper transform of $B$, and the $E_j$'s
are the corresponding exceptional curves above $B$. Note that each $E_j$ is
in fact a line $\proj^1$ representing all possible quasi-admissible covers,
arising from
a hyperelliptic curve $[C]\in B\cap \overline{\mathfrak{I}}_g$. From
Fig.~\ref{smoothhyper}, these are the blow-ups of $C$ at a point, one for
each involution pair $\{p_1,p_2\}\in g^1_2$, and that is \vspace*{3mm}why
$E_j\cong\proj^1$.

The extra summands on p.~\pageref{list of theorems}, induced by
the base curve $B$, are result of the extra intersections $pr_1^*(D)\cdot E_j$
from above. Indeed, the relations, as they stand, compute only 
$pr_1^*(D)\cdot\overline{B}$, the component corresponding to families
with general smooth members. From the calculations
on p.~\pageref{hypercalculations}, we expect that each $pr_1^*(D)\cdot E_j=1$,
and this will account for the extra $1$ apprearing in all $\mu$'s. 

\smallskip
To verify this, we
only needs to show $\delta|_{E_j}=1$. Since we cannot pick out canonically
one point $p_i$ from each hyperelliptic pair $\{p_1,p_2\}$ on $C$, and
thus construct a family of blow-ups at $p_i$ of $C$ over
$E_j\cong \proj^1$, we make a base change of degree two $C\rightarrow E_j$.

\begin{figure}[h]
$$\psdraw{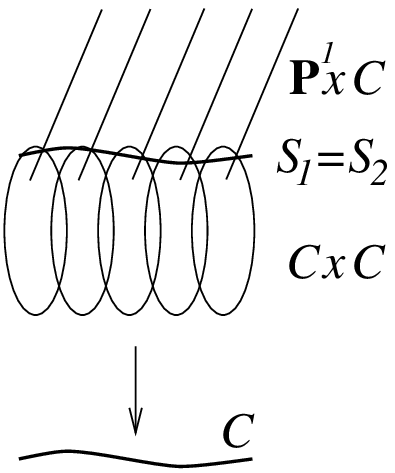}{1.1in}{1.3in}$$
\caption{$\delta|_{C}=2$}
\label{deltahyper.fig}
\end{figure}
We construct a family over $C$, corresponding to {\it all} blow-ups
of $C$ at point $p\in C$. This is simply the products $C\times C$ and
$\proj^1\times C$, identified at two sections $S_i$: $S_1$ is the diagonal on
$C\times C$, and $S_2$ is a trivial section of $\proj^1\times C$ over $C$
(cf.~Fig.~\ref{deltahyper.fig}).
From \cite{CH}, for the base curve $C$ of this family,
the degree $\delta|_C$ is computed as
\[\delta|_C=\delta_{C\times C}+\delta_{\proj^1\times C}+S_1^2+S_2^2=
0+0+2+0=2.\]
Taking into account the base \vspace*{2mm}change $C\rightarrow E_j$, 
$\delta|_{E_j}=1$. 

\smallskip Finally, if we allow for our families to have finitely many
hyperelliptic fibers, we adjust the relation in ~\ref{7+6/g relation2}
by $g\Delta{\cal H}_{3,0}\cdot B$, and the relation in
~\ref{maximal relation2} by $(5g+1)\Delta{\cal H}_{3,0}\cdot B$. The
two bounds in Theorems~\ref{7+6/g Bogomolov2}-\ref{maximal bound2} are
unaffected by the above discussion.

\setcounter{section}{12}
 
\bigskip\section*{12. Interpretation of the Bogomolov Index $4c_2-c_1^2$ 
via the Maroni Divisor}

\setcounter{subsection}{0} 
\setcounter{subsubsection}{0} 
\setcounter{lem}{0}
\setcounter{thm}{0}
\setcounter{prop}{0}
\setcounter{defn}{0} 
\setcounter{cor}{0} 
\setcounter{conj}{0} 
\setcounter{claim}{0} 
\setcounter{remark}{0}
\setcounter{equation}{0}
\label{Bog-Maroni}

\subsection{The Maroni invariant of trigonal curves}
\label{Maroniinvariant}
For any smooth trigonal curve $C$, consider the triple cover
$f:C\rightarrow {{\proj}^1}$. The pushforward $f_*(\cal{O}_{C})$,
as we noted before, is a locally free sheaf of rank 3 on ${{\proj}^1}$, and
hence decomposes into a direct sum of three invertible sheaves on
${{\proj}}^1$:
\[f_*(\cal{O}_{C})=\cal{O}_{{\proj}^1}\oplus \cal{O}_{{\proj}^1}(a)\oplus
\cal{O}_{{\proj}^1}(b).\]
The first summand is trivial due to the split exact sequence 
\[0\rightarrow {V}\rightarrow {\alpha}_*{\cal O}_{C}\stackrel
{\on{tr}}{\rightarrow}{\cal O}_{{\proj}^1}\rightarrow 0,\]
where $V=\cal{O}_{{\proj}^1}(a)\oplus\cal{O}_{{\proj}^1}(b)$.
From GRR, $a+b=g+2$. We have observed in Section~6
that $C$ embeds in the rational
ruled surface ${\proj}V={\mathbf F}_k$, for $k=|b-a|$. 

\medskip
\noindent{\bf Definition 12.1.} The {\it Maroni invariant} of an
irreducible trigonal curve
$C$ is the difference $|b-a|$. The {\it Maroni locus}
in $\overline{\mathfrak{T}}_g$ is the closure of the set of curves
with Maroni invariants $\geq 2$ (cf.~[Ma]).

\begin{lem} For a general trigonal
curve $C$ the vector bundle $V$ is {\it balanced}, i.e.
the integers $a$ and $b$ are equal or 1 apart according to
$g(\on{mod}2)$.
\label{gentrig}
\end{lem}

\begin{proof} Let $a\leq b$. 
The statement follows from a dimension count of
the linear system $L=|3B_0+\frac{g+2}{2}F|$ 
on the ruled surface ${\mathbf F}_{b-a}={\mathbf F}_k$.
Indeed, all trigonal curves with Maroni invariant $(b-a)/2$ are elements
of $L$. If $p:{\mathbf F}_k\rightarrow {\proj}^1$ is the projection map,
the projective  dimension of $L$ equals
\[r(L)=h^0\big(p_*\cal{O}_{{\mathbf F}_k}(3B_0+\textstyle{\frac{g+2}{2}}
F)\big)-1.\]
Denoting by $\widetilde{B}=B_0-\frac{k}{2}F$ the section of
${\mathbf F}_k$ with smallest self-intersection of $-k$, we have
$p_*\cal{O}_{{\mathbf F}_k}(\widetilde{B})\cong \cal{O}_{{\proj}^1}\oplus
\cal{O}_{{\proj}^1}(-k)$. The necessary pushforward from above is:
\begin{equation*}
p_*\cal{O}_{{\mathbf F}_k}(3\widetilde{B}+\textstyle{\frac{g+2+3k}{2}}F)=
\on{Sym}^3(\cal{O}_{{\proj}^1}\!\!\oplus\cal{O}_{{\proj}^1}(-k))
\otimes \cal{O}_{{\proj}^1}({\textstyle{\frac{g+2+3k}{2}}})=
\!\!\!\!\displaystyle{\bigoplus_{j=\pm 1,\pm 3}}
\!\!\!\!\cal{O}_{{\proj}^1}({\textstyle{\frac{g+2+jk}{2}}}).
\end{equation*}
Since an irreducible trigonal curve $C$ lies in $L$, we have
$C\cdot \widetilde{B}\geq 0$, hence $g+2-3k\geq 0$ and $g\equiv k
(\on{mod}2)$. Evaluating the sections of this sum of sheaves, we obtain
$r(L)=2g+7.$

The ruled surface ${\mathbf F}_k$ has automorphisms, inducing automorphisms
of the linear system $L$. We need to mod out these in order to obtain
the dimension of the space of trigonal curves embedded in ${\mathbf F}_k$.
The group $\on{Aut}{\mathbf F}_k$ is a product (not necessarily direct)
of the base automorphisms $\on{Aut}{\proj}^1=\on{PGL}_2$, and
the projective automorphisms of the vector bundle $V$. The latter is an open
set (up to projectivity)
 of the homomorphisms of $V$ into $V$, and hence has the same dimension as:
\[\on{Hom}(V,V)\cong H^0(V\otimes V\,\,\widehat{\phantom{n}})=
H^0\big(\cal{O}_{{\proj}^1}(-k)\oplus\cal{O}_{{\proj}^1}\oplus
\cal{O}_{{\proj}^1}\oplus\cal{O}_{{\proj}^1}(k)\big).\]
For $k>0$, $\on{dim}\on{Aut}V=k+3$, while for $k=0$,
$\on{dim}\on{Aut}V=4$. We conclude that the dimension of
the set of trigonal curves with Maroni invariant $k/2$ is
\[r(L)-\on{dim}\on{Aut}{\mathbf F}_k=
\left\{\begin{array}{l} 2g+1\,\,\on{if}\,\,k=0,\\
                        2g+2-k\,\,\on{if}\,\,k>0.
\end{array}\right.\]

\medskip When $k=0$ or $k=1$, this space corresponds to an open dense set of 
$\overline{\mathfrak{T}}_g$. 
For an even $g$ a general trigonal curve has Maroni invariant $0$
and therefore embeds in ${\mathbf F}_0=
{\proj}^1\times{\proj}^1$, while for an odd $g$
a general trigonal curve has Maroni invariant $1$ and embeds in
${\mathbf F}_1=\on{Bl}_{\on{pt}}(\proj^2)$. 
In both cases, the vector bundle $V$ is balanced. \end{proof}

\begin{cor} For $g$ even, the Maroni locus is a divisor in
$\overline{\mathfrak{T}}_g$ whose general member embeds in ${\mathbf F}_2$.
For $g$ odd, the Maroni locus has codimension 2 in $\overline{\mathfrak{T}}_g$ 
and its general member embeds in ${\mathbf F}_3$.
\label{maronilocus}
\end{cor}

\noindent{\bf Remark 12.1.} It will be useful to identify precisely the
group of authomorphisms of the linear system $L$ for $k=0,1$. We have
$\on{Aut}(\proj^1\!\times \proj^1)\cong PGL_2\times PGL_2\times
{\mathbb Z}/2{\mathbb Z}$. The last factor comes from the
exchange of the fiber and the base of $\proj^1\times \proj^1$ and it
is relevant only for $g=4$: then $L=|3B_0+3F|$. Otherwise,  
for any even $g>4$:
\[\on{Aut}L\cong PGL_2\times PGL_2.\] 
When $g$ is odd, the ruled surface ${\mathbf F}_1$ can be thought of as the
blow-up of $\proj^2$ at the point $q=[0,0,1]$. Any automorphism of
$\on{Bl}_q{\proj^2}$ carries the exceptional divisor $E_q$ of $\mathbf F_1$
to itself, and hence is induced by an automorphism of the plane preserving
the point $q$. The group of such automorphisms of $\proj^2$ 
is the  subgroup of $PGL_3$ corresponding to matrices:
\[\left(\begin{array}{ccc} a_{11} & a_{12} & 0\\ 
                           a_{21} & a_{22} & 0\\
                           a_{31} & a_{32} & a_{33}
        \end{array}\right).\]
Taking into account the discriminant of these matrices, we easily
identify for odd $g$:
\[\on{Aut}L\cong \mathbf A^2\times GL_2.\]
Note that all of the above groups $\on{Aut}L$ have dimension $6$, 
which was claimed already in Lemma~\ref{gentrig}.

\subsection{Generators of Pic$_{\mathbb{Q}}\overline{\mathfrak{T}}_g$}
\label{generators}

\begin{prop} The rational Picard group of 
$\overline{\mathfrak{T}}_g$, $\on{Pic}_{\mathbb{Q}}\overline{\mathfrak{T}}_g$,
is freely generated by the boundary classes
$\delta_0$, $\delta_{k,i}$, and one additional class, which for even genus
$g$ coincides with the Maroni class $\mu$.
\label{genPic}
\end{prop}

\begin{proof} Since a general trigonal curve $C$ embeds in
the ruled surface ${\mathbf F}_k$ ($k=0,1$),
$C$ is a member of the linear system $L=|3B_0+
\frac{g+2}{2}F|$ on ${\mathbf F}_k$. 
Let $U$ be the open set inside ${\proj}L\cong
\proj^{2g+7}$ corresponding to the {\it  smooth trigonal} members of $L$. 
The surjection \[{\mathbb Z}=\on{Pic}\proj^{2g+7}\twoheadrightarrow 
\on{Pic}U \]
has a nontrivial kernel, because the set of singular trigonal curves in
${\mathbf F}_k$ is a divisor in ${\proj}L$.
Hence $\on{Pic}U={\mathbb Z}/n{\mathbb Z}$ for some integer $n\!>\!0$, and 
$\on{Pic}_{\mathbb{Q}}U\!=\!0$.

\medskip
The image of the natural
projection map $p:U\rightarrow \overline{{\mathfrak{T}}}_g$  
is the open dense set $W$ of smooth trigonal curves
with lowest Marone invariant of $0$ or $1$. 
Let $F$ denote the fiber of $p$. From Remark 12.1, 
\[F\cong\left\{\begin{array}{l}
PGL_2\times PGL_2\,\,\on{if}\,\,g-\on{even},g>4;\\
PGL_2\times PGL_2\times {\mathbb Z}/2{\mathbb Z}\,\,\on{if}\,\,g=4;\\
\on{Aut}L\cong {\mathbf A}^2\times GL_2\,\,\on{if}\,\,g-\on{odd}.
               \end{array}\right.\]
Leray spectral sequence or other methods (cf.~~\cite{Gr-Ha,Milne}) yield:
\[H^1(W,f_*{\cal O}^*_U)\hookrightarrow 
H^1(U,{\cal O}^*_U).\]
Pushing  the exponential sequence on $U$ to $W$:
\[0\rightarrow {\mathbb Z}\rightarrow {\cal O}_U \rightarrow {\cal O}^*_U
 \rightarrow 0\,\,
\Rightarrow \,\,0\rightarrow {\mathbb Z}\rightarrow {\cal O}_{W}
 \rightarrow f_*{\cal O}^*_U\rightarrow R^1\!\!f_*{\mathbb Z}.\]
Combining with the exponential sequence on $W$:
\[0\rightarrow {\cal O}^*_{W}
\rightarrow f_*{\cal O}_U^* \rightarrow R^1\!\!f_*{\mathbb Z}\,\,
\Rightarrow\,\,H^1(W,{\cal O}^*_{W})
\stackrel{p^*}{\rightarrow}H^1(U,{\cal O}_X^*),\]
 with $\on{ker}p^*\subset H^0(W,R^1\!\!f_*{\mathbb Z})\subset H^1(F,{\mathbb Z})$.
For even $g$, $H^1(F,{\mathbb Z})$ is torsion (a direct sum of copies of
${\mathbb Z}/2{\mathbb Z}$), but for odd $g$ it is isomorphic to ${\mathbb Z}$.

\smallskip
Hence, for even $g$ we have the natural embedding
$p^*:\on{Pic}_{\mathbb Q}W\hookrightarrow \on{Pic}_{\mathbb Q}U$, 
and in view of $\on{Pic}_{\mathbb{Q}}U=0$, it follows that $\on{Pic}_{\mathbb{Q}}
W=0$. The complement of $W$ in $\overline{\mathfrak{T}}_g$ is
the union of the boundary of $\overline{\mathfrak{T}}_g$ and
the Maroni divisor. Therefore, $\delta_0$, $\delta_{k,i}$ and $\mu$
 generate $\on{Pic}_{\mathbb Q}\overline{\mathfrak{T}}_g$.
Since the class of the Hodge bundle $\lambda$ is {\it not} a linear 
combination of the boundary classes (cf.~p.~\pageref{list of theorems}), 
the boundary divisors are {\it not} sufficient to generate
the rational Picard group of $\overline{\mathfrak{T}}_g$, and 
$\mu$ must be linearly independent of them. We conclude that 
$\delta_0$, $\delta_{k,i}$, and $\mu$ generate freely 
$\on{Pic}_{\mathbb Q}\overline{\mathfrak{T}}_g$ for even genus $g$.

\smallskip
For $g$-odd, 
$p^*:\on{Pic}_{\mathbb{Q}}W\rightarrow \on{Pic}_{\mathbb{Q}}U$
is either an inclusion, or has a kernel with one generator. Since 
the Maroni locus for $g$-odd is not a divisor, an inclusion would
imply as above that $\lambda$ is a linear combination of the boundary
classes, which is not true. Hence, $\on{ker}p^*={\mathbb Q}$ and
$\on{Pic}_{\mathbb{Q}}W$ is generated freely by the boundary
classes $\delta_0$ and $\delta_{k,i}$, and one additional class. \end{proof}

\subsection{The Bogomolov condition and the Maroni divisor}
\label{interpretation}
\begin{prop} For even genus $g$ and a base curve $B$,
not contained in $\Delta{\mathfrak{T}}_g$:
\[(7g+6)\lambda=g\delta_0+
\sum_{k,i}\widehat{c}_{k,i}\delta_{k,i}+2(g-3)\mu,\]
where $\widehat{c}_{k,i}$ are certain polynomial coefficients computed
similarly as $\widetilde{c}_{k,i}$. (cf.~p.~\pageref{analog1})
\label{Maroni2}
\end{prop}

\begin{proof} We set $g=2(m-1)$.
Let us consider for now only families with irreducible trigonal
fibers, i.e. the base curve $B$ intersects only
the boundary component $\Delta{\mathfrak{T}}_0$.

\smallskip
\noindent{\it Case 1.}
If $B$ does not intersect the Maroni divisor $\mu$, then the Maroni
invariant of the fibers in $X$ is constant, and equal to $0$. 
The fibers $C$ of $X$ embed in the projectivization
${\proj}(V|_{F_Y})\cong {\proj}^1\times {\proj}^1$. Since $\on{deg}V|_{F_Y}
=g+2$ and $V$ is balanced, the restriction of $V$ to the fiber $F_Y$ on
the ruled surface $Y$ is
\[V|_{F_Y}=\cal{O}_{{\proj}^1}(m)\oplus 
\cal{O}_{{\proj}^1}(m).\]
Moreover, $V|_{F_Y}$ does not jump as $F_Y$ moves, so that
$V$ can be written as:
\[V\cong h^*M\otimes\cal{O}_{Y}\left(mB_0\right)\]
for some vector bundle $M$ of rank 2 on $B$. But the Bogomolov
index $4c_2(V)-c_1^2(V)$ is independent of twisting $V$ by line
bundles, in particular, by $\cal{O}_{Y}\left(mB_0\right)$, so that
\[4c_2(V)-c_1^2(V)=4c_2(h^*M)-c_1^2(h^*M)=
4c_2(M)-c_1^2(M)=0.\]
The last equality follows from $c_2(M)=0=c_1^2(M)$ for any
bundle on the curve $B$. We conclude that $4c_2(V)-c_1^2(V)=4\mu|_B=0$.

\medskip
\begin{figure}[h]
$$\psdraw{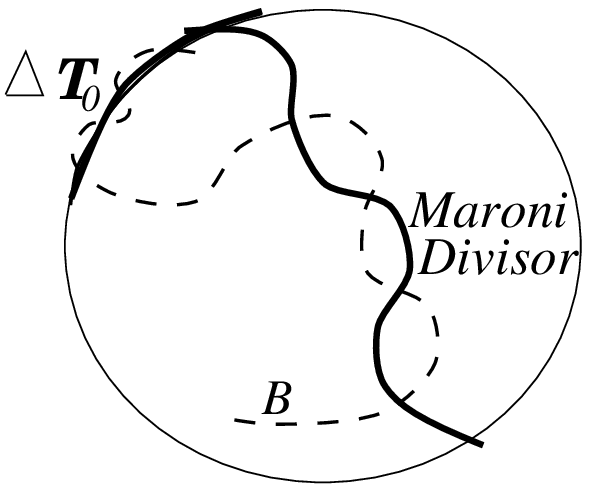}{1.5in}{1.3in}$$
\caption{$B\cap \mu$ in $\overline{\mathfrak{T}}_g$}
\label{intersectBandmu}
\end{figure}
{\it Case 2.} Now let $B$ intersect the Maroni divisor $\mu$
in {\it finitely} many points. Assume also that
these points  are {\it general} in $\mu$, i.e. they correspond to
trigonal curves $C$ embeddable in the ruled surface ${\mathbf F}_2$. 
We twist $V$ by a line bundle $M=
\cal{O}_{Y}\left(mB_0\right)$, and set $\widetilde{V}=V\otimes M$, so that
$\on{deg}\widetilde{V}|_{F_Y}=0$ and
\[\,\,\,\,\left\{\begin{array}{l}
\widetilde{V}|_{F_Y}=
\cal{O}_{{\proj}^1}\oplus \cal{O}_{{\proj}^1}\,\,\,
\phantom{(-1)(1)}\on{when}\,\,F_Y\,\,\on{is\,\,generic,}\\
\widetilde{V}|_{F_Y}=\cal{O}_{{\proj}^1}(-1)\oplus 
\cal{O}_{{\proj}^1}(1)\,\,\,\on{when}\,\,F_Y\,\,\on{is\,\,special}.
\end{array}\right.\]

Then $\widetilde{V}$ is the middle term of a short exact sequence on $Y$
\[0\rightarrow \widetilde{V}^{\prime} \rightarrow \widetilde{V} \rightarrow
\cal{I} \rightarrow 0,\]
where $\widetilde{V}^{\prime}=h^*(h_*\widetilde{V})$ is a vector bundle of
rank 2 on $Y$. In the notation of \cite{Br}, let $W$ be the 
sum of the special fibers of $Y$, and let $Z$ be the union of certain
isolated points on each member of $W$, so that 
$\cal{I}=\cal{I}_{Z\subset W}$ is the ideal sheaf of $Z$ inside $W$.
Note that the number of the special fibers, which comprise $W$,
equals $\on{deg}Z=\mu|_B$.

We can now compute the Chern classes of $\widetilde{V}$:
\[\left\{\begin{array}{l}
c_1(\widetilde{V})=c_1(\widetilde{V}^{\prime})+W=\on{a\,\,sum\,\,of\,\,
fibers\,\,of\,\,Y},\\
c_2(\widetilde{V})=c_2(\widetilde{V}^{\prime})+\on{deg}Z=\on{deg}Z.
\end{array}\right.\]
The last equality follows from the fact that $\widetilde{V}^{\prime}$
is the pull-back of a bundle on the curve $B$, hence of zero higher
Chern classes. We conclude that $c_1^2(\widetilde{V})=0$, and
\[4c_2({V})-c_1^2({V})=
4c_2(\widetilde{V})-c_1^2(\widetilde{V})=4\on{deg}Z=4\mu|_B.\]

Putting the above two cases together, we have for any family with 
irreducible trigonal members, not entirely contained in the Maroni locus:
\begin{equation}
4c_2({V})-c_1^2({V})=4\mu|_B.
\end{equation}
Prop.~\ref{genPic} then implies that $\lambda$ is a linear
combination of the boundary and the Maroni class:
\[(7g+6)\lambda|_{\overline{\mathfrak{T}}_g}=g\delta_0+\sum_{k,i}
\widehat{c}_{k,i}\delta_{k,i}+2(g-3)\mu,\]
where the coefficients $\widehat{c}_{k,i}$ are computed in a similar way, or by
direct computation with families of singular trigonal curves
(cf.~\cite{CH}). \end{proof}

\medskip
We can combine the above results in the following
\begin{thm} For even $g$, $\on{Pic}_{\mathbb{Q}}\overline{\mathfrak{T}}_g$ is
freely generated by all boundary classes $\delta_0$ and $\delta_{k,i}$, and
the Maroni class $\mu$. The class of the Hodge bundle on
$\overline{\mathfrak{T}}_g$ is expressed in terms of these generators as
the following linear combination:
\begin{equation*}
(7g+6)\lambda|_{\overline{\mathfrak{T}}_g}=g\delta_0+
\sum_{k,i}\widehat{c}_{k,i}\delta_{k,i}+2(g-3){\mu}.
\end{equation*}
\label{Pic trigonal}\vspace*{-5mm}
\end{thm}

\noindent{\bf Remark 12.2.} Note that the coefficients
$\widehat{c}_{k,i}$ depend on the specific decriptions of the Maroni 
curves that appear in the boundary divisors $\Delta{\mathfrak{T}}_{k,i}$,
and they are {\it not} always equal to the corresponding
coefficients $\widetilde{c}_{k,i}$ in Theorem \ref{7+6/g relation2}.
Indeed, in the above Proposition, we have shown that
\begin{equation}
4c_2({V})-c_1^2({V})=4\mu|_B+
\sum_{k,i}\alpha_{k,i}\delta_{k,i},
\label{alpha-coef}
\end{equation}
for some $\alpha_{k,i}$, which may be non-zero. Hence,
 $\widehat{c}_{k,i}=\widetilde{c}_{k,i}+\frac{g-3}{2}\alpha_{k,i}$.

\smallskip
For example, consider the case of $\Delta{\mathfrak{T}}_{1,i}$, and let
$C=C_1\cup C_2$ be a general member of it. If $C$ is also Maroni, then
there exists a family $X\rightarrow B$, whose general fiber is
an irreducible Maroni curve, and one of whose special fibers is our $C$.
We can assume, modulo a base change and certain blow-ups not affecting
$C$, that this family fits in the basic construction diagram (cf.~
Fig.~\ref{general B}). Let ${\mathbf R}_1$ and ${\mathbf R}_2$ be
the two ruled surfaces in which $C_1$ and $C_2$ are embedded, and
let $E_1$ and $E_2$ be the projections of $C_1$ and $C_2$ in the
birationally ruled surface $\widehat{Y}$. Then $F=E_1+E_2$ is
a special fiber of $\widehat{Y}$, with self-intersections $E_1^2=E_2^2=-1$.

\smallskip
Now, the general member of $X$, being Maroni, is embedded in 
a ruled surface $\mathbf F_2$ with a section $L$ of self-intersection $-2$.
The union of such $L$'s forms a surface in the 3-fold $\mathbf PV$, whose
closure we denote by $S$. Evidently, $S\cong \widehat{Y}$, at least
outside their special fibers. Let
$S$ intersect ${\mathbf R}_1$ and ${\mathbf R}_2$ in curves $L_1$ and
$L_2$ (over $E_1$ and $E_2$). 
We claim that at least one of ${\mathbf R}_1$ and ${\mathbf R}_2$
is {\it not} isomorphic to $\mathbf F_0=\mathbf P^1\times \mathbf P^1$.
It will suffice to show that $L_1$ or $L_2$ has negative self-intersection.

\smallskip
Indeed, suppose to the contrary that $L_m^2\geq 0$ in ${\mathbf R}_m$ 
($m=1,2$). Note that $S\cdot \mathbf R_m=L_m$ in $\mathbf PV$, so that
\[L_m^2=S|_{\mathbf R_m}\cdot S|_{\mathbf R_m}= S^2\cdot \mathbf R_m\,\,
\Rightarrow\,\,S^2(\mathbf R_1+\mathbf R_2)\geq 0.\]
On the other hand, $\mathbf R_1+\mathbf R_2$ is the fiber of the
projection $\mathbf PV\rightarrow \widehat{Y}$, and as such it is
linearly equivalent to the general fiber $\mathbf F_2$. Hence
\[0\leq S^2\cdot \mathbf F_2=S|_{\mathbf F_2}\cdot S|_{\mathbf F_2}=L^2=-2,\]
a contradiction.
We conclude that if $C=C_1\cup C_2$ is a Maroni curve of boundary type
$\Delta{\mathfrak{T}}_{1,i}$, then either $C_1$ or $C_2$ (or both)
is embedded in a ruled surface $\mathbf F_k$ with $k\geq 1$. This
already distinguishes the cases of odd and even genus $i$.

\smallskip
When $i=g(C_1)$ is even (and hence $j=g(C_2)=g-j-2$ is also even),
the general member of $\Delta{\mathfrak{T}}_{1,i}$ is embedded in a join
of two $\mathbf F_0$'s (each $C_m\subset \mathbf F_0$), and hence
it is {\it not} Maroni. Based on this observation, one can easily find
the coefficient $\alpha_{1,i}$ for $i$-even. To do this, consider
the birationally ruled surface $Y$ which is the blow-up of 
$\mathbf F_0$ at one point. Let again the two components of the
special fiber of $Y$ be $E_1$ and $E_2$, and projectivize the trivial vector
bundle $V=\mathcal{O}_Y\oplus \mathcal{O}_Y$: ${\mathbf P}V=Y\times
\mathbf P^1$. By taking an appropriate linear system in ${\mathbf P}V$,
one obtains a family of trigonal curves $X$, whose fibers are all
irreducible and embedded in $\mathbf F_0$, except for a special
reducible curve $C$ over $E_1\cup E_2$ of the specified above type.
Hence none of $X$'s members are Maroni, and so $\mu|_B=0$. Further,
$4c_2(V)-c_1^2(V)=0$, and $\delta_{1,i}|_B=1$, so that 
equation~(\ref{alpha-coef}) implies $\alpha_{1,i}=0$, and
hence $\widehat{c}_{k,i}=\widetilde{c}_{k,i}$ for $i$-even.

\smallskip
The situation is quite different when the genus $i$ is odd. Then both
components of the general member $C$ of $\Delta{\mathfrak{T}}_{1,i}$
are embedded in $\mathbf F_1$'s, and hence $C$ is potentially Maroni.
One can take further the above general argument of 
intersection theory on ${\mathbf P}V$, and show that the
curves $L_1$ and $L_2$ are in fact both 
sections of negative self-intersection $-1$ in these $\mathbf F_1$'s:
consider the product $S\cdot X \cdot {\mathbf F}_2$ and
its variation over the special fiber of $\widehat{Y}$. But we know
that $L_1$ and $L_2$ intersect, as the fiber of $S$ over
$\widehat{Y}$ is connected.

Thus, the curve $C$ would be Maroni if and only if the two corresponding
ruled surfaces $\mathbf F_1$ are glued along one of their fibers so that
their negative sections intersect on that fiber. (This decsription
can be alternatively derived by considering the degenerations of 
the $g^1_3$'s on the irreducible Maroni curves.)
To find $\alpha_{1,i}$ in this case, we construct a similar example
as above, only changing $V$ to $\mathcal O_Y\oplus \mathcal O_Y(E_1)$.
This, while keeping the general fiber embedded in $\mathbf F_0$, has the 
effect of embedding the special one in a ``Maroni'' gluing of two
$\mathbf F_1$'s. We have $4c_2(V)-c_1^2(V)=-E_1^2=1$, $\mu|_B=1$,
and $\delta_{1,i}|_B=1$, so that equation~(\ref{alpha-coef}) implies 
$\alpha_{1,i}=-3$ for $i$-odd, and hence 
$\widehat{c}_{k,i}=\widetilde{c}_{k,i}-3/2(g-3)$.

\bigskip
\begin{figure}[t]
$$\psdraw{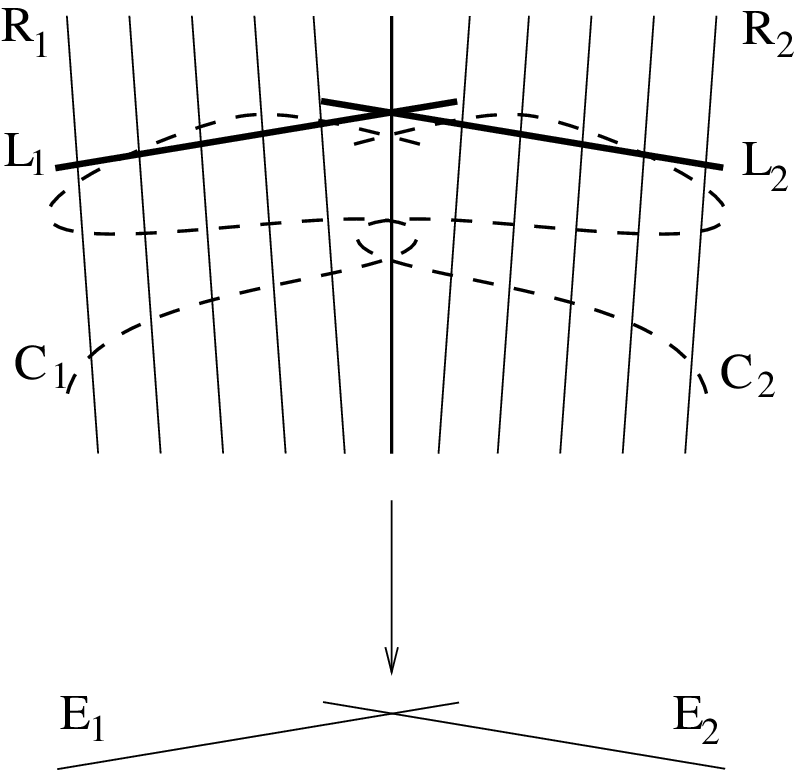}{1.5in}{1.5in}$$
\caption{Maroni curves in $\Delta_{1,i}\mathfrak{T}_g$, $i$-odd}
\label{Maroni-boundary}
\end{figure}

One can similarly compute the remaining coefficients 
$\alpha_{k,i}$, by first figuring out which boundary curves in
$\Delta_{k,i}\mathfrak{T}_g$ are Maroni,  then constructing
an appropriate vector bundle $V$, and finally using
equation~(\ref{alpha-coef}) to compute $\alpha_{k,i}$, and
hence $\widehat{c}_{k,i}$. \qed

\begin{prop} For $g$-even, if the base curve $B$ is not entirely
contained in the Maroni divisor, and the singular members of
$X$ belong only to $\Delta_0\mathfrak{T}_g
\cup \Delta_{1,i}\mathfrak{T}_g$, then the slope of the family $X/_B$
satisfies:
\[\frac{\delta}{\lambda}\leq 7+\frac{6}{g}.\]
\label{Maroni inequality}\vspace*{-5mm}
\end{prop}

\begin{conj} For $g$-even, if the base curve $B$ is not entirely
contained in the Maroni divisor, then the slope of the family $X/_B$
satisfies:
\[\frac{\delta}{\lambda}\leq 7+\frac{6}{g}.\]
\label{Maroni-conj}
\end{conj}

\subsection{The Maroni divisor and the maximal bound}
\label{Maroni-maximal}
Even though for odd genus $g$ the Maroni locus is not large enough
to be a divisor in $\overline{\mathfrak{T}}_g$, we can define a 
{\it generalized Maroni} divisor class by extending the relation from
the $g$-even case.

\medskip
\noindent{\bf Definition 12.2.} For any genus, we define the
{\it generalized Maroni} class $\mu$
in $\on{Pic}_{\mathbb{Q}}\overline{\mathfrak{T}}_g$ by
\[\mu:=\frac{1}{2(g-3)}\big\{(7g+6)\lambda-g\delta_0-
\sum_{k,i}\widehat{c}_{k,i}\delta_{3,i}\big\}.\]

\begin{thm} The maximal bound $36(g+1)/(5g+1)$ is attained
for a trigonal family of curves $X\rightarrow B$ if and only
all fibers of $X$ are irreducible and 
\[\delta_0|_B=-\frac{72(g+1)}{g+2}\mu|_B\]
\label{maximalmaroni}\vspace*{-5mm}
\end{thm}

\begin{proof} The fact that $X$ must have only irreducible fibers
in order to attain the maximum bound is already known from
Theorem~\ref{genmaximal}. This means $\delta_{k,i}|_B=0$ for all
$k,i$. Then, Theorem~\ref{bogomolov1} implies:
\begin{equation}
(7g+6)\lambda|_B=g\delta_0|_B+\frac{g-3}{2}\mu|_B.
\end{equation}
Assume that the maximal bound is attained, i.e. 
$36(g+1)\lambda|_B=(5g+1)\delta_0|_B$. Substituting for
$\lambda|_B$ in the above equation, yields the desired equality.
The converse follows similarly. \end{proof}

\medskip
\noindent{\bf Remark 12.3.} In the $g$-even case, this equality has a 
specific meaning. Since the Maroni class $\mu$ corresponds to an
effective divisor on $\overline{\mathfrak{T}}_g$, the equality (and hence
the maximal bound) is achieved only for base curves $B$ entirely contained
in the Maroni divisor, so that the restriction $\mu|_B$ can be negative.
In fact, in all found examples, the base $B$ is contained in a very small
subloci of the Maroni loci, defined by the highest possible Maroni invariant.

\medskip
\noindent{\bf Remark 12.4.} Theorem~\ref{Pic trigonal}
and Prop.~\ref{maximalmaroni} 
do not have analogs in the hyperelliptic case: there is no additional
Maroni divisor to generate $\on{Pic}_{\mathbb Q}\overline{\mathfrak{I}}_g$
together with the boundary $\Delta\mathfrak{I}_g$.

\medskip
\noindent{\bf Remark 12.5.} When $g=3$, there is no Maroni locus in
$\overline{\mathfrak{T}}_3$ either. Indeed, since an irreducible
 trigonal curve of genus $3$ embeds
only in ruled surfaces ${\mathbf F}_k$ with $k$-odd and
$k\leq (g+2)/3=5/3$, then {\it all} irreducible trigonal curves embed in
${\mathbf F}_1$, and correspondingly they all have the lowest possible
Maroni invariant $k=1$. However, $\on{Pic}_{\mathbb Q}\overline{\mathfrak{T}}_3$
is not generated by the boundary classes of $\overline{\mathfrak{T}}_3$: as
Prop.~\ref{genPic} asserts, in the odd genus case there is always one
additional generating class.

\smallskip
On the other hand, the results on p.~\pageref{list of theorems}
yield apriori {\it two} relations among $\lambda$ and the $\delta_{k,i}$'s.
This would have been a contradiction to the {\it freeness} of the
generators above, unless these two relations are the same. This is in fact
what happens:
\[9\lambda=\delta_0+3\delta_{2,1}+3\delta_{3,1}+4\delta_{4,1}+4\delta_{5,1}
+3\delta_{5,2}+3\delta_{6,1},\]
as restricted to any base curve $B\not\subset\Delta\overline{\mathfrak{T}}_3$.
Note the convenient disappearance of the ``extra''
$(g-3)$--summands in the coefficients of $\delta_{4,i},\delta_{5,i},
\delta_{6,i}$).
Then the maximal and the semistable ratios
 both equal $9$, and are attained for families with
irreducible trigonal members.

\bigskip\section*{13. Further Results and Conjectures}

\setcounter{section}{13}
\setcounter{subsection}{0} 
\setcounter{subsubsection}{0} 
\setcounter{lem}{0}
\setcounter{thm}{0}
\setcounter{prop}{0}
\setcounter{defn}{0} 
\setcounter{cor}{0} 
\setcounter{conj}{0} 
\setcounter{claim}{0} 
\setcounter{remark}{0}
\setcounter{equation}{0}
\label{furtherresults}

\subsection{Results and conjectures for $d$-gonal families, $d\geq 4$}
I have carried out some preliminary research in the $d$-gonal case, and
while the methods and ideas for the trigonal case are
in principle extendable, this appears to be a substantially more
subtle and complex problem.
More precisely, let $\overline{\cal{D}}_d$ be the closure in
$\overline{\mathfrak{M}}_g$ of the stable curves expressible as $d$-sheeted covers 
of ${\proj}^1$. One possible goal is 
to complete the program of describing generators and relations 
for the rational Picard groups $\on{Pic}_{\mathbb{Q}}\overline{\cal{D}}_d$, 
and to find the exact maximal bounds for the slopes of $d$-gonal 
families.

\smallskip
For example, I have obtained the following bound for the slope of a 
general tetragonal family with smooth general member (for odd genus $g$):
\[\frac{\delta}{\lambda}\leq 6\frac{2}{3}+\frac{64}{3(3g+1)}=
\frac{4(5g+7)}{3g+1}.\]

I have also conjectured formulas for the maximal
and general bounds for any $d$-gonal and other families 
of stable curves. Entering these formulas are the {\it Clifford index} of 
curves, {\it Bogomolov semistability} conditions for higher rank bundles,
and some new geometrically described loci in $\overline{\cal{D}}_d$. 
Generalizing the idea of the Maroni locus in the trigonal case, 
these loci are characterized, for example, in the tetragonal case by the
dimensions of the multiples of the $g^1_4$-series. In particular,
there will be another generator of $\on{Pic}_{\mathbb Q}\overline{\mathfrak{T}}_4$
besides the boundary and Maroni divisors. 

\smallskip
In the following I present some of these conjectures
on the upper bounds for $\overline{\cal{D}}_d$. 
We start by comparing all known maximal and general
bounds functions of the genus $g$:
\[\begin{array}{|c|c|c|c|c|c|}
   \hline\hline
   \stackrel{\vspace*{1mm}}{\on{locus\,\, in \,\,}\overline{\mathfrak{M}}_g}
&\on{bound}& g=1& g=2& g=3& g=5\\
   \hline\hline\vspace*{1mm}
    \on{general}\,\overline{\mathfrak{M}}_g&
    \stackrel{\vspace*{1mm}}{\displaystyle{ 6+\frac{12}{g+1\vspace*{1mm}}}}
 & 12 & 10 &9 &8\\
  \hline\vspace*{1mm}
    \on{hyperelliptic}\,\overline{\cal{H}}_g=\overline{\cal{D}}_2&
    \stackrel{\vspace*{1mm}}{\displaystyle{8+ \frac{4}{g\vspace*{1mm}}}} 
  & 12 & 10 & - &-\\
   \hline\vspace*{1mm}
    \on{trigonal}\,\overline{\cal{T}}_g=\overline{\cal{D}}_3
& \stackrel{\vspace*{1mm}}{\displaystyle{ \frac{36(g+1)}{5g+1
\vspace*{1mm}}}} 
& 12 & - & 9 &- \\
   \hline\vspace*{1mm}
    \on{ gen. tetragonal=\overline{\cal{D}}_4}
&\stackrel{\vspace*{1mm}}{\displaystyle{
   \frac{4(5g+7)}{3g+1\vspace*{1mm}}}}& 12 & - & -&8\\
   \hline
\end{array}\]

\medskip
The pattern appearing in this table is clear: the general bound
$6+\displaystyle{12/(g+1)}$ 
coincides with each of the other bounds exactly twice
for some special values of the genus $g$. Evidently, 
$g=1$ is one of these special values, yielding 12 everywhere. (I
owe this observation to Benedict Gross.)
Let $g_d$ be the other genus $g$ for which the
general formula in $\overline{\mathfrak{M}}_g$ and the maximal formula for
$\overline{\cal{D}}_d$ coincide, i.e. $g_2=2$, $g_3=3$, $g_5=5$. 
We notice that for these genera $g_d$
 the  moduli spaces $\overline{\mathfrak{M}}_2,\overline{\mathfrak{M}}_3$
and $\overline{\mathfrak{M}}_5$  consist only of 
hyperelliptic, trigonal or tetragonal curves, respectively.
In general, {\it Brill-Noether} theory (cf.~\cite{ACGH}) asserts  
that  for complete linear series
$g^r_d=g^1_d$ the expected dimension of the variety of $g^1_d$'s on
a smooth curve of genus $g$ is $\rho=g-(r+1)(g-d+r)=2(d-1)-g,$ and hence
the smallest genus $g$ for which
$\overline{\mathfrak{M}}_g=\overline{\cal{D}}_d\supsetneq
\overline{\cal{D}}_{d-1}$
is $g=2d-3$. Thus we set $g_d=2d-3$ for $d\geq 3$ and $g_2=2$. Note that
this coincides with the previously found $g_3=3$ and $g_5=5$.

\begin{conj} If $\cal{F}_d(g)$ is an exact upper bound for the slopes of
families of stable curves with smooth $d$-gonal general member (locus
$\overline{\cal{D}}_d$), then
\begin{eqnarray*}
&&(a)\,\,\cal{F}_d(1)=12.\\
&&(b)\,\,\cal{F}_d(g_d)=6+\displaystyle{\frac{12}{g_d+1}}\cdot
\end{eqnarray*}
\label{conj2}\vspace*{-5mm}
\end{conj}

It is reasonable to expect that the upper bounds for $\overline{\cal{D}}_d$
will be ratios of linear functions of the genus $g$:
$\cal{F}_d(g)=(Ag+B)/(Cg+D)$.
Conjecture~\ref{conj2} then estimates the difference between $\cal{F}_d(g)$
and the general bound for $\overline{\mathfrak{M}}_g$ up to a factor
$f_d=D/C$. 

\begin{conj} The exact upper bounds $\cal{F}_d(g)$ are given by
\[\cal{F}_d(g)=6+\frac{12}{g+1}+6\frac{(1-f_d)(g-g_d)(g-1)}{(g+f_d)(g_d+1)(g+1)},\]\vspace*{-3mm}
or equivalently, \vspace*{-3mm}
\[\cal{F}_d(g)=6+\frac{6}{g+f_d}\left(1+f_d+\frac{1-f_d}{g_d+1}(g-1)\right).\]
\end{conj}

I have a conjecture on how to determine the remaining 
factor $f_d$, which seems to be closely related to the coefficients
of the linear expression in [EMH] 
for the divisor $\overline{\cal{D}}_{\frac{g+1}{2}}$ 
in terms of the Hodge bundle $\lambda$ and
the boundary classes $\delta_i$ on $\overline{\mathfrak{M}}_g$. 
These conjectures are supported by the work of Cornalba-Harris on the
{ hyperelliptic locus} $\overline{\cal{H}}_g=\overline{\cal{D}}_2$,
by the results of this paper
on the { trigonal locus} $\overline{\cal{T}}_g=\overline{\cal{D}}_3$,
and by partial results on the tetragonal locus $\overline{\cal{D}}_4$.

\smallskip In view of Remark 12.5, the equality between the maximal
and semistable trigonal bounds for $g=3$ suggests that a similar situation
might occur for other $d$-gonal families. It is reasonable to expect
two or more ``semistable'' bounds, depending on the number of extra
generators in $\on{Pic}_{\mathbb Q}{\overline{\cal D}}_d$.

\smallskip One of these ``semistable'' bounds relates to 
families obtained as blow-ups of pencils of $d$-gonal curves on
a ruled surface ${\mathbf F}_k$. Example 2.1 yields the maximal bound
$8+4/g$ for hyperelliptic families (no extra generator besides
the boundary classes), and a similar example in the trigonal case
yields the $7+6/g$ semistable bound (one extra generator, the
Maroni locus). We generalize this to any $d$-gonal family of
curves embedded in an arbitrary ruled surface ${\mathbf F}_k$. Invariably, the
slope of $X/\!_{\displaystyle{B}}$ is:
\begin{equation*}
\frac{\delta|_B}{\lambda|_B}=\left(6+\frac{2}{d-1}\right)+\frac{2d}{g}\cdot
\end{equation*}

\medskip
\begin{conj}
 Let $X$ be a family of $d$-gonal curves of genus $g$ whose
base $B$ is not contained in a certain codimension 1 closed subset of
$\overline{\cal D}_d$. Then the slope of $X/\!_{\displaystyle{B}}$ satisfies:
\begin{equation*}
\frac{\delta|_B}{\lambda|_B}\leq \left(6+\frac{2}{d-1}\right)+\frac{2d}{g}\cdot
\end{equation*}
\label{clifford} \vspace*{-5mm}
\end{conj}

Conjectures~\ref{clifford}--4
are modifications of earlier conjectures of Joe Harris.

\subsection{A look at families with special $g^r_d$'s, $r\geq 2$}
The discussion so far was primarily
concerned with the loci $\overline{\cal{D}}_d\subset
\overline{\mathfrak{M}}_g$ corresponding
to linear series $g^1_d$. But all of our problems are well-defined
and quite interesting to solve for curves with series $g^r_d$ of dimension
 $r>1$. Equivalently, we consider the loci $\overline{\cal{D}}^r_d$
of curves mapping with degree $d$ to ${\proj}^r$, $r\geq 1$. 

\medskip
\noindent{\bf Definition 13.1.} The {\it Clifford index} $\mathfrak {c}$
of a smooth curve $C$ is defined as 
\[\mathfrak{c}=\on{min}_L\left\{\on{deg} {L} -2\on{dim}{L}\right\}\]
where $L$ runs over all effective special linear series ${L}$
on $C$.

\medskip
Clifford's theorem implies
${\mathfrak{c}}\geq 0$, with equality if and only if $C$ is
hyperelliptic, i.e. ${L}=g^1_2$ (cf.~\cite{ACGH}). 
On the other hand, ${\mathfrak{c}}=1$
means that there exists a $g^r_d$ on $C$ with $d-2r=1$. From Marten's
Theorem, $\on{dim}W^r_d(C)\leq d-2r-1=0$, where $W^r_d$ is the
variety parametrizing complete linear series on $C$ of degree $d$ 
and dimension at least $r$. Therefore, we must have
$\on{dim}W^r_d=0$. But then Mumford's theorem
asserts that $C$ is either trigonal, or
bi-elliptic, or a smooth plane quintic. The bi-elliptic case would mean
that $W^r_d$ consists of $g^2_6$'s, which contradicts the dimension of
$\on{dim}W^r_d$. In short, $\mathfrak{c}=1$
if and only if $C$ is not hyperelliptic and possesses a $g^1_3$ or a $g^2_5$.

\smallskip
Thus, according to the Clifford index, 
the first case with $r\geq 2$ is the space of plane quintics. Consider a
general pencil of such, and blow up the plane at its 25 base
points. The resulting family $X=\on{Bl}_{25}{\proj^2}\rightarrow\proj^1$
is easily seen to have slope $8=7+6/g$, which corresponds to
the bound in Conjecture~\ref{clifford} with $d-2$ replaced by the Clifford
index $\mathfrak{c}=1$.  Finally, note that for a $d$-gonal curve
$C$ of genus $g$, by definition $\mathfrak{c}\leq d-2$, so that when
$g\gg d$ we may generalize to:

\begin{conj} For a general family $X\rightarrow B$ of genus $g$ stable curves
whose general member has Clifford index $\mathfrak{c}$ and whose base
$B$ is a general curve in $\overline{\cal D}^r_d$, 
the slope of $X/\!_{\displaystyle{B}}$ satisfies:
\begin{equation*}
\frac{\delta_X}{\lambda_X}\leq\left(6+\frac{2}{\mathfrak{c}+1}\right)+
\frac{2\mathfrak{c}+4}{g}\,\,\,\on{for}\,\,\mathfrak{c}<\!<g\cdot
\label{clifford1}
\end{equation*}
\end{conj}

\noindent{\bf Remark 13.1.} It is worth noting that the stratification
of $\overline{\mathfrak{M}}_g$, for which we asked
in the Introduction, is not obtained
via the Clifford index $\mathfrak{c}$. For example, Xiao constructs families
of bi--elliptic curves $C$ with slope $8$ (cf.~\cite{Xiao}), 
which is between the hyperelliptic and the trigonal maximal bounds. 
Since $C$ has a $g^1_4$ as bi--elliptic, this already exceeds the
conjectured maximal bounds for the tetragonal case. This shows
that in some of the above
conjectures we have to exclude the subset of
bi--elliptic curves from the tetragonal locus $\overline{\cal D}_4$, and that
similar modifications might be necessary for the
 other loci $\overline{\cal D}_d$. More precisely, it seems plausible
that the stratification of $\overline{\mathfrak{M}}_g$ according to
successively lower slope bounds is related not just to the existence
of a specific linear series $g^r_d$, but also to the number, dimension and
description of the irreducible components of corresponding varieties $W^r_d$.

\subsection{Other methods via the moduli space 
$\overline{\cal{M}}_{g,n}({\proj}^r,d)$}

The approach in the $g^1_d$-cases is based on a modification 
of the Harris-Mumford's [EHM] {\it Hurwitz scheme
of admissible covers}, which parametrized the
$d$-uple covers of stable pointed rational curves.
However, in the more general situation for linear series with larger 
dimensions $r>1$, such a compactification via
admissible covers does not exist, so we have to look for a
different solution.

\smallskip
Consider moduli spaces of stable maps
$\overline{\cal{M}}_{g,n}({\proj}^r,d)$. They parametrize {\it stable}
maps $(C,p_1,p_2,...,p_n;\mu)$, where $C$ is a projective, connected
nodal curve of arithmetic genus $g$, the $p_i$'s are marked
distinct nonsingular points on $C$, and the map
$\mu:C\rightarrow{\proj}^r$ has image $\mu_*([C])=d[\on{line}]$ and
satisfies certain stability conditions (cf.~\cite{K,KM}). The  space
$\overline{\cal{M}}_{g,n}({\proj}^r,d)$ seems to be the right 
compactification which we need in order to extend our results to
families with $g^r_d$-series on the fibers: the moduli space
of stable maps is somewhat more ``sensitive'' in describing our loci
$\overline{\cal{D}}^r_d$ in terms of their geometry. 

\smallskip
Going back to the $g^1_d$-problems, one can also see the combinatorial
flavor that stands in the background of these questions. 
It is probably not coincidental that the spaces 
$\overline{\cal{M}}_{g,n}({\proj}^r,d)$ are also combinatorially
defined and give rise to many enumerative problems. It will be useful
to understand better the loci $\overline{\cal{D}}^r_d$ via their
connection with the Kontsevich spaces
$\overline{\cal{M}}_{g,n}({\proj}^r,d)$, and ultimately to solve the
remaining questions on $\on{Pic}_{\mathbb{Q}}\overline{\cal{D}}^r_d$
for any $d,r$, as well as related interesting enumerative problems that
will inevitably arise from such considerations.

\section*{14. Appendix: The Hyperelliptic Locus $\overline{\mathfrak{I}}_g$}

\setcounter{section}{14}
\setcounter{subsection}{0} 
\setcounter{subsubsection}{0} 
\setcounter{lem}{0}
\setcounter{thm}{0}
\setcounter{prop}{0}
\setcounter{defn}{0} 
\setcounter{cor}{0} 
\setcounter{conj}{0} 
\setcounter{claim}{0} 
\setcounter{remark}{0}
\setcounter{equation}{0}
\label{hyperelliptic}

In this section we give a proof of Theorems~\ref{theoremCHPic} and~\ref{CHX},
following the same ideas and methods as in the trigonal case. We refer
the reader to previous sections for a detailed proof of certain
statements.

\subsection{Boundary locus of $\overline{\mathfrak{I}}_g$}
\label{hyperellipticboundary}
\begin{figure}[h]
$$\psdraw{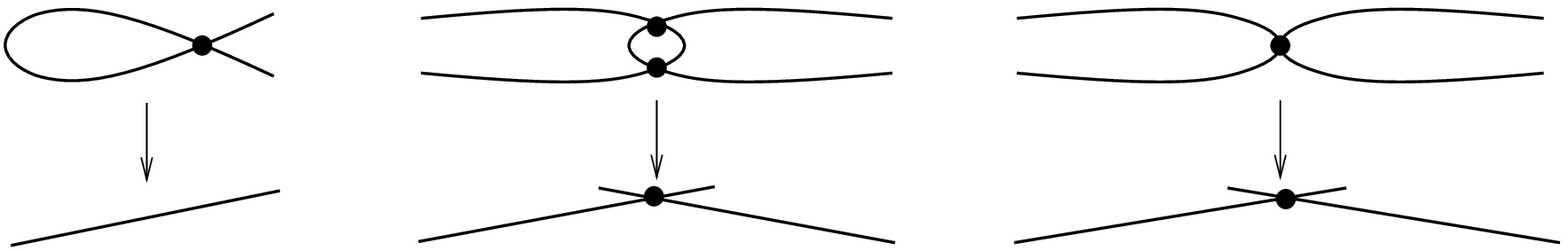}{4.5in}{0.9in}$$
\caption{Boundary of the hyperelliptic locus $\overline{\mathfrak{I}}_g$}
\label{hyperboundary}
\end{figure}

Cornalba-Harris
describe the boundary of $\overline{\mathfrak{I}}_g$ as consisting of
several boundary components, whose general members and indexing are shown
in Fig.~\ref{hyperboundary} (cf.~\cite{CH}). The restriction of the divisor class $\delta$ to $\overline{\mathfrak{I}}_g$
is the following linear combination:
\begin{equation}
\delta\big|_{\overline{\mathfrak{I}}_g}=\delta_0+2\sum_{i=1}^
{[(g-1)/2]}\xi_i+
\sum_{j=1}^{\left[g/2\right]}\delta_j,
\label{boundaryrel}
\end{equation}
where $\xi_i$ and $\delta_i$ are the classes in
$\on{Pic}_{\mathbb{Q}}\overline{\mathfrak{I}}_g$ of the boundary
divisors $\Xi_i$ and $\Delta_j$. 

\newpage
\vspace*{-8mm}$$\Xi_0;\,\,\Xi_i,\,{\scriptstyle {i=1,...,[(g-1)/2]}}
;\,\,\Delta_j,\,{\scriptstyle{j=1,...,[g/2]}}$$

\subsection{Effective covers and embedding for hyperelliptic families}
\label{embeddinghyperelliptic}
In the case of a hyperelliptic family $f:X\rightarrow B$, a minimal
quasi-admissible cover coincides with the original family $X$, because no
blow-ups are necessary to perform on the fibers of $X$: these are
already quasi-admissible double covers. Thus, we have a
degree 2 map $\phi=\widetilde{\phi}:X\rightarrow Y$ for some birationally
ruled surface $Y$ over $B$. As for an effective cover
$\widehat{\phi}:\widehat{X}\rightarrow \widehat{Y}$, only the boundary
divisors $\Delta_i$ require blow-ups (cf.~Fig.~\ref{hyperboundary}). 
This is analogous to the ``ramification index 1''  discussion in 
Fig.~\ref{ram}--\ref{resolve1}. Thus, while in $\widehat{X}$ the special
fibers may have occasional nonreduced rational components of multiplicity
2, the fibers of $\widehat{Y}$ are always 
trees of reduced smooth ${\proj}^1$'s.

\medskip
In the case of a smooth hyperelliptic curve $C$, we consider the
natural double sheeted map $f:C\rightarrow {\proj}^1$. The pushforward
$f_*{\cal{O}_C}$ is a rank 2 vector bundle on ${\proj}^1$, which fits into the
short exact sequence
\begin{equation*}
0\rightarrow {\cal{O}_{{\proj}^1}(g+1)}\rightarrow {f}_*{\cal O}_{C}\stackrel
{\on{tr}}{\rightarrow}{\cal O}_{{\proj}^1}\rightarrow 0.
\end{equation*}
We can embed $C$ in the rational ruled surface 
${\proj}((f_*\cal{O}_C)\,\,\hat{})$.
We generalize this construction
to the effective cover $\widehat{\phi}:\widehat{X}\rightarrow \widehat{Y}$
by setting $V:=({\phi}_*{\cal O}_{X})\,\,\hat{}$. For some line
bundle $E$ on $\widehat{Y}$:
\begin{equation*}
0\rightarrow {E}\rightarrow {\widehat{\phi}}_*{\cal O}_{\widehat{X}}\stackrel
{\on{tr}}{\rightarrow}{\cal O}_{\widehat{Y}}\rightarrow 0.
\end{equation*}
Then $\widehat{X}$ naturally 
embeds in the threefold ${\proj}V$. Let $\pi:{\proj}V\rightarrow
\widehat{Y}$ be the corresponding projection map.

\subsection{The invariants $\lambda,\delta$ and $\kappa$}
\label{Hyperinvariants}
As a divisor in ${\proj}V$, $\widehat{X}\equiv 2\zeta+\pi^*D$,
for some divisor $D$ on $\widehat{Y}$. From the adjunction
formula, $g=\on{deg}c_1(V)|_{F_{\widehat{Y}}}-1=c-1$, where
$c_1(V)=cB_0+dF_Y$.
The arithmetic genus of the inverse image $\widehat{\phi}^*T(E)$ is
given by \[p_{\!\stackrel{\phantom{.}}{E}}=-m_{\!\stackrel{\phantom{.}}{E}}
\left(\Gamma_{\!\stackrel{\phantom{.}}{E}}+
\Theta_{\!\stackrel{\phantom{.}}{E}}\right).\]
 It turns out that these are
the only differences between the set-up of the hyperelliptic and the
trigonal case. The definitions of the functions $m,\theta$ and $\gamma$,
as well as the formulas for $c_1(V), K_{{\proj}V}, c_2({\proj}V)$ and
the congruence $D\equiv 2c_1(V)$ are valid without any modifications. 

\smallskip
As in the trigonal case, it will be sufficient to consider only the
cases when the base curve $B$ intersects {\it transversally}
the boundary divisors of $\overline{\mathfrak{I}}_g$.
But then for all non-root components $E$ in $\widehat{Y}$:
\[m_{\!\stackrel{\phantom{.}}{E}}=1=\Theta_{\!\stackrel{\phantom{.}}{E}}\,\,
\on{and}\,\,
\Gamma_{\!\stackrel{\phantom{.}}{E}}=-(p_{\!\stackrel{\phantom{.}}{E}}+1).\]
We can now easily calculate  the invariants on $X$.

\begin{prop} For any family $f:X\rightarrow B$ of hyperelliptic curves
with smooth general member and a base curve $B$ intersecting transversally
the boundary of $\overline{\mathfrak{I}}_g$:\
\begin{eqnarray*}
\lambda_X&\!\!=\!\!&dg+\frac{1}{2}\sum_{E\not =R}\Gamma_
{\!\stackrel{\phantom{.}}{E}}
(\Gamma_{\!\stackrel{\phantom{.}}{E}}+1),\\
\kappa_X&\!\!=\!\!&4d(g-1)-2\sum_{E\not
=R}(\Gamma_{\!\stackrel{\phantom{.}}{E}}+1)^2+\sum_{\on{ram}1}1,\\
\delta_X&\!\!=\!\!&4d(2g+1)+2\sum_{E\not = R}
(\Gamma_{\!\stackrel{\phantom{.}}{E}}+1)(1-2\Gamma_
{\!\stackrel{\phantom{.}}{E}})+\sum_{\on{ram}1}1.
\end{eqnarray*}
\label{hyperinvariants}
\end{prop}

With this, we are ready to show the linear relations among $\lambda|_B$ and
the boundary restrictions $\delta_i|_B$ and $\xi_i|_B$. It is
evident that in order to cancel the ``global'' term $d$, one
must subtract $(8g+4)\lambda_X|_B-g\delta|_B$, which is the main
idea of the next theorem.

\begin{thm} There exists an effective linear combination $\cal{E}_h$ of
the boundary divisors of $\overline{\mathfrak{I}}_g$, not containing
$\Xi_0$, such that for any family $f:X\rightarrow B$ of hyperelliptic curves
with smooth general member:
\[(8g+4)\lambda_X|_B=g\delta|_B+\cal{E}_h|_B\]
\label{hyperrelation}\vspace*{-8mm}
\end{thm}

\begin{proof} We consider the difference
\begin{eqnarray*}
\mathfrak{S}_h&\!\!=\!\!&(8g+4)\lambda_X|_B-g\delta|_B=
 2\sum_{E\not = R}(1+\Gamma_{\!\stackrel{\phantom{.}}{E}})(g+\Gamma_
{\!\stackrel{\phantom{.}}{E}})+\sum_{\on{ram}1}g\\
&\!\!=\!\!&2\sum_{E\not = R}p_{\!\stackrel{\phantom{.}}{E}}(g-1+p_
{\!\stackrel{\phantom{.}}{E}})+\sum_{\on{ram}1}g.
\end{eqnarray*}
In the hyperelliptic case, as opposed to the trigonal case, there
is only {\it one type} of non-root components $E$, namely, such that
both $E$ and $E^-$ are reduced. That is why there is just one type of summands
in $\mathfrak{S}_h$.

\smallskip
As in Section~\ref{arbitrary}, it is sufficient to calculate the
above sum for general members of $\Xi_{i}$ and
$\Delta_i$, as described in Prop.~\ref{Delta-k,i}, i.e. for a {\it general}
base curve $B$.

\subsubsection{Contribution of the boundary divisors $\Xi_{i}$}
\label{hypercontribution1}
This case is analogous to the case of $\Delta_{3,i}$ (cf.~
Subsection~\ref{contribution1}). The arithmetic genus 
$p_{\!\stackrel{\phantom{.}}{E}}=g-i-1$, and the corresponding summand in
$\mathfrak{S}_h$ is
\[e_i=2p_{\!\stackrel{\phantom{.}}{E}}(g-1+p_
{\!\stackrel{\phantom{.}}{E}})
=2i(g-i-1)>0,\]
where $i=1,...,[(g-1)/2]$.

\subsubsection{Contribution of the boundary divisors $\Delta_{j}$}
\label{hypercontribution2}
Compare this with the contribution of $\Delta_{5,j}$ (subsection
\ref{contribution2}). There are two non-root components $E_1$ and
$E_2$ in the special fiber of $\widehat{Y}$ ($E_1^-=R$),
whose invariants are $p_{\!\stackrel{\phantom{.}}{E_1}}=g-j-1$ and
$p_{\!\stackrel{\phantom{.}}{E_2}}=g-j$. With the
ramification adjustment of $g$, the contribution of $\Delta_j$
to the sum $\mathfrak{S}_h$ is
\begin{equation*}
f_j=
2p_{\!\stackrel{\phantom{.}}{E_1}}(g-1+p_{\!\stackrel{\phantom{.}}{E}})+
p_{\!\stackrel{\phantom{.}}{E_2}}(g-1+p_{\!\stackrel{\phantom{.}}{E}})+g
=4j(g-j)-g>0,
\end{equation*}
where $j=1,...,[g/2]$.

\medskip
Finally,  for the appropriate indices $i$ and $j$
we set $\displaystyle{\cal{E}_h:=\sum_{i>0}e_i\Xi_i+\sum_{j>0}f_j\Delta_j.}$
This is an effective combination of
boundary divisors in $\overline{\mathfrak{I}}_g$, not containing $\Delta_0$
by construction, and satisfying $\mathfrak{S}_h=\cal{E}_h|_B$. \end{proof}

\medskip
Theorem~\ref{hyperrelation} implies immediately the following
\begin{cor}
Let $f:X\rightarrow B$ be a nonisotrivial
family with smooth general member. Then the slope of the family satisfies:
\begin{equation}
\frac{\delta|_B}{\lambda|_B}\leq 8+\frac{4}{g}.
\label{second8+4/g}
\end{equation}
Equality holds if and only if the general fiber of $f$ is hyperelliptic,
and all singular fibers are irreducible.
\end{cor}
It is now straightforward to prove the fundamental
relation in $\on{Pic}_{\mathbb{Q}}\overline{\mathfrak{I}}_g$, shown first
in \cite{CH}. In Theorem~\ref{hyperrelation}, 
we add to the coefficients $e_i$ and $f_j$
the corresponding multiplicities $\on{mult}_{\delta}\xi_i$ and
$\on{mult}_{\delta}\delta_j$: 
\[\widetilde{e}_i=e_i+2\cdot g=2(i+1)(g-i),\,\,\,
\widetilde{f}_j=f_j+1\cdot g=4j(g-j).\] Using the fact that 
$\on{Pic}_{\mathbb{Q}}\overline{\mathfrak{I}}_g$ is generated freely
by the boundary classes $\xi_i$ and $\delta_j$ (see \cite{CH}), 
we obtain
\[(8g+4)\lambda=g\delta_0+\sum_{i>0}\widetilde{e_i}\xi_i+
\sum_{j>0}\widetilde{f_j}\delta_j.\]

\begin{thm} In the Picard group of the hyperelliptic locus,
$\on{Pic}_{\mathbb{Q}}\overline{\mathfrak{I}}_g$, the class of the Hodge bundle
$\lambda$ is expressible in terms of the boundary divisor classes of
$\overline{\mathfrak{I}}_g$ as:
\begin{equation*}
(8g+4)\lambda=g\xi_0+\sum_{i=1}^{[(g-1)/2]}2(i+1)(g-i)\xi_i
+\sum_{j=1}^{[g/2]}4j(g-j)\delta_j.
\end{equation*}
\label{CHPic2}
\end{thm}

\vspace{10mm}
\begin{flushright}\begin{minipage}[r]{2.8in}
{\sc Zvezdelina Stankova-Frenkel}\\
Mathematical Sciences Research Institute\\
1000 Centennial Dr., Berkeley, CA 94720\\
e-mail address: {\tt stankova@msri.org}
\end{minipage}
\end{flushright}
\end{document}